\pdfoutput=1 
\documentclass{amsart}
\usepackage{lipsum}
\usepackage{subfigure}
\usepackage{hyperref}
\usepackage{graphicx,color}
\usepackage{stix} 

\usepackage{amssymb,amsmath,bm,mathtools}
\usepackage{enumerate}
\usepackage{tkz-graph}
\usepackage{ragged2e}
\usepackage{enumitem}
\usepackage{adjustbox}
\usepackage{tipa}
\usepackage{circuitikz} 
\usepackage{scalerel}
\usepackage[margin=1.25in]{geometry}
\usepackage{parskip}
\usepackage{comment}

\usepackage{tikz-cd}
\usepackage{adjustbox}
\usepackage{tikz}
\usepackage[all,knot]{xy}
\usepackage{transparent}
\usepackage{graphicx,import}

\usetikzlibrary{arrows,decorations.pathmorphing}
\usetikzlibrary{backgrounds,positioning}
\usetikzlibrary{fit,petri,shapes.misc}
\usetikzlibrary{shapes.geometric}
\usetikzlibrary{decorations.markings}
\usetikzlibrary{calc}
\usetikzlibrary{braids}

\newlength{\defaultpgflinewidth}
\tikzset{auto}

\tikzset{empty/.style={circle,inner sep=0pt,minimum size=6mm}}
\tikzset{emptyvt/.style={circle,inner sep=0pt,minimum size=0mm}}

\tikzset{plain/.style={circle,draw,very thick,
inner sep=0pt,minimum size=6mm}}

\tikzset{xplain/.style={circle,draw,very thick,
inner sep=0pt,minimum size=8mm}}

\tikzset{smallplain/.style={circle,draw,very thick,
inner sep=0pt,minimum size=4mm}}

\tikzset{xsplain/.style={circle,draw, thick,
inner sep=0pt,minimum size=1.5mm}}

\tikzset{tinyplain/.style={circle,draw, thick,
inner sep=0pt,minimum size=1mm}}

\tikzset{smalldotted/.style={circle,draw,very thick, densely dotted,
inner sep=0pt,minimum size=3mm}}

\tikzset{dottedplain/.style={circle,draw,very thick, densely dotted,
inner sep=0pt,minimum size=6mm}}

\tikzset{normaldot/.style={circle,black,fill=black,inner sep=0pt,minimum size=2mm}}

\tikzset{rectplain/.style={rectangle,draw,very thick,minimum size=6mm}}
\tikzset{bigplain/.style={rectangle,draw,very thick,minimum size=1cm}}

\tikzset{triangular/.style={regular polygon, regular polygon sides=3, draw,very thick,
inner sep=0pt,minimum size=1.2cm}}

\tikzset{arrow/.style={->,thick}}
\tikzset{dashedarrow/.style={->,dashed,thick}}
\tikzset{dottedarrow/.style={->,dotted,thick}}
\tikzset{mapto/.style={|->,thick}}
\tikzset{->-/.style={decoration={markings, mark=at position #1 with {\arrow{>}}},postaction={decorate}}}

\tikzset{implies/.style={thick,double,double equal sign distance,-implies}} 

\tikzset{line/.style={thick}}
\tikzset{dottedline/.style={dotted,thick}}
\tikzset{dashedline/.style={dashed,thick}}

\tikzset{inputleg/.style={<-,thick}}
\tikzset{outputleg/.style={->,thick}}
\tikzset{dottedinput/.style={<-,dotted,thick}}

\usepackage{transparent}
\usepackage{graphicx,import}
\usepackage{caption}

\usepackage{enumitem}

\usepackage{todonotes}
\newcommand\reallywidehat[1]{\arraycolsep=0pt\relax%
\begin{array}{c}
\stretchto{
  \scaleto{
    \scalerel*[\widthof{\ensuremath{#1}}]{\kern-.5pt\bigwedge\kern-.5pt}
    {\rule[-\textheight/2]{1ex}{\textheight}} 
  }{\textheight} %
}{0.5ex}\\           
#1\\                 
\rule{-1ex}{0ex}
\end{array}
}

\newcommand{\bG}{\mathbb{G}}

\newcommand{\bC}{\mathbf{C}}
\newcommand{\bD}{\mathbf{D}}

\newcommand{\bT}{\mathbf{T}}
\newcommand{\bE}{\mathbf{E}}

\newcommand{\bB}{\mathbf{B}}

\newcommand{\calO}{\mathcal{O}}
\newcommand{\calQ}{\mathcal{Q}}
\newcommand{\calC}{\mathcal{C}}

\newcommand{\sSet}{\mathbf{S}}

\newcommand{\Set}{\mathbf{Set}}

\newcommand{\Grpd}{\mathbf{G}}

\newcommand{\Cat}{\mathsf{Cat}}

\newcommand{\Grp}{\mathbf{Gr}}

\newcommand{\kk}{\mathbb{K}}    

\newcommand{\CoRB}{\mathsf{CoRB}}
\newcommand{\PaRB}{\mathsf{PaRB}}

\newcommand{\BV}{\mathsf{BV}}
\newcommand{\magma}{\mathsf{\Omega}}

\newcommand{\PaB}{\mathsf{PaB}}

\newcommand{\calP}{\mathcal{P}}

\newcommand{\Top}{\mathrm{Top}}

\newcommand{\twist}{\tau}

\newcommand{\Cyc}{\mathbf{Cyc}}
\newcommand{\Op}{\mathbf{Op}}

\newcommand{\Prop}{\mathbf{Prop}}

\newcommand{\id}{\operatorname{id}}

\newcommand{\Pro}{\operatorname{Pro}}
\newcommand{\End}{\operatorname{End}}
\newcommand{\Hom}{\operatorname{Hom}}
\newcommand{\Ho}{\operatorname{Ho}}

\newcommand{\ob}{\operatorname{ob}}
\newcommand{\Aut}{\operatorname{Aut}}

\newcommand{\Fun}{\operatorname{Fun}}
\newcommand{\ev}{\operatorname{ev}}
\newcommand{\coev}{\operatorname{coev}}

\newcommand{\Map}{\operatorname{Map}}
\newcommand{\Alg}{\operatorname{Alg}}

\newcommand{\Env}{\operatorname{Env}}
\newcommand{\cs}{\textnormal{B}}
\newcommand{\nerve}{\textnormal{N}}

\newcommand{\gt}{\mathsf{GT}}
\newcommand{\GT}{\widehat{\mathsf{GT}}}

\newcommand{\Br}{\mathsf{B}} 
\newcommand{\PB}{\mathsf{PB}} 

\newcommand{\RB}{\mathsf{RB}} 
\newcommand{\PRB}{\mathsf{PRB}} 

\newcommand{\Free}{\mathsf{F}} 

\newcommand{\Conf}{\mathsf{Conf}}

\newtheorem{theorem}{Theorem}[section]
\newtheorem{thm}[theorem]{Theorem}
\newtheorem{prop}[theorem]{Proposition}
\newtheorem{lemma}[theorem]{Lemma}
\newtheorem{cor}[theorem]{Corollary}
\newtheorem*{thm*}{Theorem}
\newtheorem*{prop*}{Proposition}

\theoremstyle{definition}
\newtheorem{definition}[theorem]{Definition}
\newtheorem{example}[theorem]{Example}
\newtheorem{remark}[theorem]{Remark}

\numberwithin{equation}{section}

\let\oldtocsection=\tocsection
\let\oldtocsubsection=\tocsubsection
\renewcommand{\tocsection}[2]{\hspace{0em}\oldtocsection{#1}{#2}}
\renewcommand{\tocsubsection}[2]{\hspace{1em}\oldtocsubsection{#1}{#2}}
\DeclareRobustCommand{\SkipTocEntry}[5]{}
\renewcommand{\paragraph}[1]{\textbf{{#1}.}\hspace{5pt}}

\title[GT Symmetries of Cyclic Operads and Tangles]{Grothendieck-Teichm\"uller Symmetries of Cyclic Operads and Tangles}

\author[M. Robertson]{Marcy Robertson}
\address{School of Mathematics and Statistics \\ The University of Melbourne \\ Melbourne, Victoria, Australia}
\email{marcy.robertson@unimelb.edu.au}

\author[C. Singh]{Chandan Singh}
\address{School of Mathematics and Statistics \\ The University of Melbourne \\ Melbourne, Victoria, Australia}
\email{chandan.singh1@unimelb.edu.au}

\begin{document}

\begin{abstract}
We identify the profinite Grothendieck--Teichmüller group $\GT$ with the group of homotopy automorphisms of the profinite completion of the cyclic operad of parenthesised ribbon braids. We transport this action to a profinite cyclic $\infty$-operad of framed configuration spaces. Using the metric prop associated to the cyclic ribbon-braid operad, we also obtain a $\GT$-action on a category of framed unoriented profinite tangles, and compare its prounipotent analogue with the action of Kassel and Turaev. Finally, the rational cyclic Grothendieck--Teichmüller action gives an alternative proof of the rational formality of the cyclic framed little-disks operad.
\end{abstract}

\maketitle

Let $\mathcal M_{0,n+1}$ denote the moduli space of smooth
genus-zero curves with $n+1$ marked points. Since these moduli spaces
are defined over $\mathbb Q$, their geometric étale fundamental
groupoids carry actions of the absolute Galois group
$\operatorname{Gal}(\overline{\mathbb Q}/\mathbb Q)$. The first
nontrivial example is
$\mathcal M_{0,4}\cong
\mathbb P^1\setminus\{0,1,\infty\}$, and Belyi's theorem implies that
the induced outer Galois action on its geometric étale fundamental
group is faithful~\cite{belyi}. Grothendieck proposed in
\cite{grothendieck_esquisse} that these actions should be considered
together, as an action on the tower formed by all moduli spaces of
curves and the natural maps between them.

The genus-zero part of this tower has a concrete topological model.
The topological fundamental groups of the spaces
$\mathcal M_{0,n+1}$ are pure mapping class groups of punctured
spheres and are closely related to pure braid groups. Their geometric
étale fundamental groups are the corresponding profinite
completions; see
Examples~\ref{example: etale fundamental groups} and
\ref{example: Oda's theorem}. The maps between the moduli spaces
become operations relating braid groups with different numbers of
strands. Thus the genus-zero tower may be modelled by a tower of braid
groups equipped with these operations.

The Grothendieck--Teichmüller groups arise as symmetries of
completions of this braid-group model. Drinfeld introduced their
profinite and proalgebraic versions and constructed their actions on
completed braid groups~\cite{Drin}. Ihara's study of the Galois
action on profinite braid groups gives an injection
\[
\operatorname{Gal}(\overline{\mathbb Q}/\mathbb Q)
\hookrightarrow\GT;
\]
see \cite{ihara}. Lochak and Schneps made this tower perspective
precise for braid groups. They identify $\GT$ with the automorphism
group of a tower built from profinite completions of the Artin braid
groups, equipped with natural inclusion and strand-doubling
homomorphisms, and obtain an analogous pro-$\ell$ statement up to
inner automorphisms~\cite{ls}.

The compatibility between braid groups with different numbers of
strands also admits a natural operadic formulation. The operad
$\PaB$ consists of parenthesised braids, with operadic composition
given by braid insertion. Building on Bar--Natan's categorical
description of associators~\cite{BN}, Fresse identifies the
proalgebraic Grothendieck--Teichmüller group $\gt(\kk)$ with the
group of object-fixing automorphisms of the Malcev completion
$\PaB_{\kk}$; see \cite[Theorem~11.1.7]{FresseBook1}.

Horel proved the profinite analogue of Fresse's result, identifying
$\GT$ with the group of homotopy automorphisms of the profinite
completion of $\PaB$~\cite{Horel_profinite_groupoids}. Building on
this result, Boavida de Brito, Horel, and the first author proved the
corresponding identification for
$\PaRB$~\cite{Boavida-Horel-Robertson}. Here, the ribbon analogue
$\PaRB$ consists of parenthesised ribbon braids and records a framing
of each strand. Following \cite{campos2019configuration}, we write
$\PaRB^{\mathrm{cyc}}$ for $\PaRB$ equipped with the cyclic structure
introduced there. In Section~\ref{subsec: cyclic PaRB}, we describe
this structure explicitly in terms of the generators of $\PaRB$ and
verify its compatibility with the defining relations and operadic
composition.

This description allows us to determine how the cyclic structure
interacts with automorphisms after profinite completion. Although
preserving the cyclic structure is a priori an additional condition,
we prove that every object-fixing automorphism of
$\widehat{\PaRB}$ does so automatically. Consequently,
\[
\Aut_{\Op}(\widehat{\PaRB})
\cong
\Aut_{\Cyc}(\widehat{\PaRB}^{\mathrm{cyc}})
\cong
\GT,
\]
where $\Aut_{\Op}$ and $\Aut_{\Cyc}$ denote the group of object-fixing automorphisms, of
operads and cyclic operads respectively. The corresponding
statement for endomorphism monoids is
Theorem~\ref{thm: GT monoid iso to cyclic endomorphisms}.

To apply this result to topological models, we replace strict
automorphisms by homotopy automorphisms. For this purpose, we work
with a homotopy-coherent model of cyclic operads. We recall the Segal
model for cyclic $\infty$-operads in
Section~\ref{sec: cyclic infinity operad} and carry out the
homotopy-automorphism calculation in
Section~\ref{subsec: homotopy automorphisms of PaRB}. The result is
the following; see Corollary~\ref{cor: GT as cyclic operad}.

\begin{thm*}
There is a natural isomorphism
\[
\GT
\cong
\pi_0\mathbb R\Aut_{\Cyc_\infty(\widehat{\Grpd})}
\bigl(
\nerve\widehat{\PaRB}^{\mathrm{cyc}}
\bigr).
\]
\end{thm*}

Here $\widehat{\Grpd}$ denotes the category of profinite groupoids,
$\Cyc_\infty(\widehat{\Grpd})$ denotes the homotopy theory of cyclic
$\infty$-operads in profinite groupoids, and $\nerve$ is the cyclic
dendroidal nerve. The notation $\mathbb R\Aut$ denotes the derived
space of self-equivalences, so its set of connected components is the
group of homotopy automorphisms.

We next return to the geometric origin of the cyclic structure on
$\PaRB$. The operad $\PaRB$ is a groupoid model for the framed little
disks operad; see \cite[Chapter~1]{Wahl_Thesis}. Budney proved that
the framed little disks operad admits a cyclic
model~\cite{Budney}. One such model is the cyclic operad
$\mathsf S_0$ of genus-zero surfaces with parametrised boundary
components, recalled in Example~\ref{ex: Budney}. Its cyclic
structure permutes all boundary components, removing the distinction
between inputs and output.

A second cyclic model is the framed Fulton--MacPherson operad
$\mathsf{FFM}_2$, whose cyclic structure was described by Campos,
Idrissi, and Willwacher~\cite{campos2019configuration}. In planar
coordinates, the transposition exchanging an input with the output
is represented by the Möbius involution $z\mapsto-1/z$. Passing to
the fundamental groupoid induces the cyclic structure on $\PaRB$.
In Section~\ref{subsec: cyclic PaRB}, we verify this action explicitly
on the generators of $\PaRB$.

This geometric description motivates realising the preceding
calculation of automorphisms on a topological model. The homotopical
methods discussed above are needed for two related reasons. Framed
configuration spaces do not themselves assemble into an operad,
since the relevant composition maps are defined only up to coherent
homotopy. Moreover, profinite completion of spaces does not in
general preserve the finite products appearing in operadic
composition. We therefore work with cyclic $\infty$-operads,
modelled by cyclic dendroidal Segal spaces.

In Section~\ref{subsec:framed_conf_operad}, we adapt Hackney's
construction~\cite{hackney_config} to model framed configuration
spaces by a semi-dendroidal space $X^{\mathrm{fr}}$. We exhibit a
zigzag of weak equivalences between $X^{\mathrm{fr}}$ and the
underlying non-unital rooted part of $\nerve\mathsf S_0$, and compare
$\nerve\mathsf S_0$ with the classifying-space nerve of
$\PaRB^{\mathrm{cyc}}$. Combined with our calculation of the homotopy
automorphisms of $\widehat{\PaRB}^{\mathrm{cyc}}$, this gives the
following result
(Theorem~\ref{thm: GT acts on widehat csN PaRB}).

\begin{thm*}
There is a natural isomorphism
\[
\GT
\cong
\pi_0\mathbb R\Aut_{\Cyc_\infty(\widehat{\sSet})}
\left(
\bigl(\nerve\mathsf S_0\bigr)^{\wedge}
\right).
\]
Consequently, the injection
$\operatorname{Gal}(\overline{\mathbb Q}/\mathbb Q)\hookrightarrow
\GT$ induces a faithful action by homotopy automorphisms of the
profinite cyclic $\infty$-operad
$\bigl(\nerve\mathsf S_0\bigr)^{\wedge}$.
\end{thm*}

Cyclic operads relate operations with different choices of
output. One of their original motivations is that a cyclic
$\calO$-algebra is an $\calO$-algebra equipped with an invariant
nondegenerate symmetric pairing, and hence with a form of
self-duality~\cite{gk_cyc}. Hinich and Vaintrob encode this
structure in the metric prop $\Pi(\calO)$~\cite{hv_cyclic},
which adjoins evaluation and coevaluation operations compatible
with the cyclic structure. We recall their construction in
Section~\ref{subsec: envelope of cyclic operad}.

For parenthesised ribbon braids, this relationship has a
geometric interpretation in terms of tangles. Algebras over
$\PaRB$ are balanced monoidal categories
\cite[Theorem~1.4.7]{Wahl_Thesis}. Their associator, braiding,
and twist are represented by parenthesised ribbon braids.
Cups and caps supply the additional duality present in ribbon
categories. The category $\bT$ of framed oriented tangles
is the ribbon monoidal category freely generated by one
object~\cite{shum_tortile_1994}. Forgetting strand orientations
makes the generating object self-dual, with the cap and cup
providing its evaluation and coevaluation. We denote the
unoriented category by $\bT'$ and its parenthesised version
by $q\bT'$. The objects of $q\bT'$ are parenthesised words
in a single generating strand. See
\cite[Remark~4.1 and Section~8.2]{mousaaid2021affinization}.

In Section~\ref{subsec:tangles}, we use ribbon duality to
extend the $\PaRB$-action on $q\bT'$ to an action of
$\Pi(\PaRB^{\mathrm{cyc}})$. We describe this extension
through set-valued functors of the input and output words,
so that the cyclic action can exchange inputs and outputs.
We then pass to the profinite setting. Using Furusho's
construction of framed profinite tangles
\cite[Appendix~B]{furusho_galois_2017}, we define the
topological category $q\widehat{\bT}'$ and extend the
construction using topologically enriched functors
(Section~\ref{subsec: profinite tangles}).

An element of $\GT$ transforms the parenthesised profinite
ribbon-braid operations. Compatibility with the cyclic
structure extends this transformation to
$\Pi(\widehat{\PaRB}^{\mathrm{cyc}})$.
In Section~\ref{subsec:GT_action_profinite_tangles}, we
use precomposition to obtain a right $\GT$-action on the
algebra structures carried by $q\widehat{\bT}'$.
The underlying topological category remains fixed, while
its algebra structure is transported. This is the sense
of the action in the following result
(Corollary~\ref{cor:Galois_action_qT}).

\begin{thm*}
The group $\GT$ acts on parenthesised profinite tangles
$q\widehat{\bT}'$. Consequently, the absolute Galois group
$\operatorname{Gal}(\overline{\mathbb Q}/\mathbb Q)$
also acts on $q\widehat{\bT}'$.
\end{thm*}

Over a field $\kk$ of characteristic zero, we show that this
transport is realized by ribbon monoidal equivalences between
the original and transported ribbon structures on the completed
$\kk$-linear tangle category
(Section~\ref{subsec:GT action on turaev category}).
In Proposition~\ref{prop: comparison with KT action}, we
identify the resulting $\gt(\kk)$-action with the construction
of Kassel and Turaev \cite[Appendix~D]{KT_Tangles}.

Our final application concerns the rational cyclic formality of the
framed little disks operad. This was proved by Campos, Idrissi, and
Willwacher in the category of dg Hopf cooperads
\cite[Theorem~3.1]{campos2019configuration}; see also the related
earlier work of Giansiracusa and Salvatore over the real numbers
\cite{Salvatore_2012_cyclic_formality}. In
Section~\ref{subsec: formality}, we give an alternative proof using
the cyclic Grothendieck--Teichmüller action and the
grading-automorphism strategy of
Petersen~\cite{petersen_formality}.

\begin{thm*}
The cyclic operad $\mathsf{FD}_2^{\mathrm{cyc}}$ is rationally
formal.
\end{thm*}

\subsection*{Organisation of the paper}

Sections~\ref{sec: categorical background} and
\ref{sec: Operads, Cyclic operads} recall the categorical,
homotopical, and operadic background, including cyclic operads,
symmetric monoidal envelopes, and metric props.
Section~\ref{sec: cyclic infinity operad} develops the required
homotopy theory of cyclic operads and cyclic $\infty$-operads and
constructs the cyclic $\infty$-operad of framed configurations.

In Section~\ref{sec: cyclic operad of ribbon braids}, we describe the
cyclic structure on $\PaRB^{\mathrm{cyc}}$ explicitly and verify its
compatibility with the defining relations and operadic composition.
Section~\ref{sec: GT} identifies its strict and homotopy
endomorphisms after profinite completion with the profinite
Grothendieck--Teichmüller monoid and transports this description to
the cyclic $\infty$-operad of framed configurations.

Section~\ref{sec: Galois action on tangles} identifies the self-dual
ribbon-braid construction with parenthesised profinite tangles and
constructs the resulting Galois action.
Section~\ref{sec: Turaev's tangle category and rational formality
of the little disks operad} treats prounipotent completions, compares
the induced action with that of Kassel and Turaev, and proves the
rational cyclic formality of the framed little disks operad.

\subsection*{Acknowledgements}

The authors thank Pedro Boavida de Brito, Geoffroy Horel, and Philip
Hackney for many helpful discussions on related material over the
years, and Luci Bonatto and Chris Rogers for conversations on related
projects that directly influenced the development of this work. We would like to thank the anonymous referee for suggesting the crucial cofibrant replacement fix in Section~\ref{subsec: cofibrant cyclic PaRB}.
C.S.\ acknowledges helpful discussions during the Masterclass:
Infinity Operads and Applications to Geometry held in Copenhagen in
August 2025. C.S.\ acknowledges the support of an
Australian Government Research Training Program (RTP) Scholarship
and the Dr Albert Shimmins Fund.  M.R.\ and C.S.\ acknowledge the support of the
Australian Research Council through Future
Fellowship~FT210100256. 
\tableofcontents
\section{Categorical Background}\label{sec: categorical background}
In this first section, we aim to provide a short, self-contained overview of the categorical and homotopical language we use throughout the paper.  In Section~\ref{subsec: tensor categories}, we review the various types of monoidal categories that will be used in the later parts. For a more detailed review, we suggest \cite{kassel_quantum_1995} or \cite{Etingof2015TensorCategories}. In Section~\ref{subsec: Homotopical Preliminaries}, we briefly review the background on Quillen model categories and $\infty$-categories that we will need to discuss maps between monoidal categories and operads up to homotopy. 

\subsection{Monoidal Categories}\label{subsec: tensor categories} 
A \emph{monoidal category} $(\bE, \otimes, \mathbb{1}, \alpha, \lambda, \rho)$ consists of a category $\bE$ equipped with a tensor product $\otimes$, a unit object $\mathbb{1}$, and natural isomorphisms consisting of an associator \[\alpha_{X,Y,Z} : (X \otimes Y) \otimes Z \to X \otimes (Y \otimes Z)\] and left and right unit maps \[\lambda_X : \mathbb{1} \otimes X \to X \quad \text{and
} \quad \rho_X : X \otimes \mathbb{1} \to X,\] subject to coherence conditions expressed via commutative diagrams (the pentagon and triangle identities). A monoidal category is called \emph{strict} if the associativity and unit isomorphisms are identity morphisms.

With the exception of the braid and tangle categories introduced in
Section~\ref{sec: Galois action on tangles}, we assume that our
monoidal categories are cocomplete and that the monoidal product
$\otimes$ preserves small colimits in each variable. A category $\bE$ is simplicially enriched if, for every pair of
objects $X,Y\in\bE$, it has a simplicial mapping space
$\Map_{\bE}(X,Y)$, together with simplicial composition and units.
We will use simplicially enriched categories that are also tensored
and cotensored over simplicial sets. Examples include the categories
of simplicial sets, $\sSet$, and groupoids, $\Grpd$, with the
enrichments recalled below, as well as their profinite counterparts
used later in the paper. To lighten notation, we will often omit the
ambient category and write simply $\Hom(X,Y)$ or $\Map(X,Y)$ when no
confusion can arise.

In Section~\ref{sec: Turaev's tangle category and rational formality of the little disks operad} we will work with $\kk$-linear categories, where $\kk$ is a field of characteristic zero. A monoidal category $(\bE,\otimes,\mathbb{1})$ is \emph{$\kk$-linear} if each hom-set $\Hom_{\bE}(X,Y)$ is a
$\kk$-module and composition and tensor product are bilinear. A functor between $\kk$-linear categories is said to be \emph{$\kk$-linear} if it induces $\kk$-linear maps on hom-sets. 

The composition of morphisms is written in the usual order: for $f\colon X\to Y$ and $g\colon Y\to Z$ we write $g\cdot f\colon X\to Z$.

A \emph{braiding} on a strict monoidal category
$(\bE,\otimes,\mathbb{1})$ is a family of natural isomorphisms
\[
c_{X,Y}\colon X\otimes Y \longrightarrow Y\otimes X
\]
such that, for all objects $X,Y,X',Y',Z\in\bE$ and morphisms
$f\colon X\to X'$ and $g\colon Y\to Y'$, the following identities hold:
\begin{equation}
\label{eq:naturality}
c_{X',Y'}\cdot(f\otimes g)
=
(g\otimes f)\cdot c_{X,Y},
\end{equation}
\begin{equation}
\label{eq:hexagon1}
c_{X,Y\otimes Z}
=
(\id_Y\otimes c_{X,Z})
\cdot
(c_{X,Y}\otimes\id_Z);
\end{equation}
\begin{equation}
\label{eq:hexagon2}
c_{X\otimes Y,Z}
=
(c_{X,Z}\otimes\id_Y)
\cdot
(\id_X\otimes c_{Y,Z}).
\end{equation}

A \emph{strict braided monoidal category} is a strict monoidal category
equipped with a braiding. A \emph{strict symmetric monoidal category}
is a strict braided monoidal category whose braiding satisfies
\[
c_{Y,X}\cdot c_{X,Y}=\id_{X\otimes Y}
\]
for all objects $X,Y\in\bE$.
A \emph{balancing} on a strict braided monoidal category $\bE$ is a natural automorphism $\Theta\colon\id_{\bE}\Rightarrow\id_{\bE}$ such that $\Theta_{\mathbb{1}}=\id_{\mathbb{1}}$ and, for all objects $X,Y\in\bE$,
\begin{equation}
\label{eq:balancing}
\Theta_{X\otimes Y}
=
c_{Y,X}\cdot c_{X,Y}\cdot(\Theta_X\otimes\Theta_Y).
\end{equation}
A \emph{balanced monoidal category} is a braided monoidal category equipped with a balancing.

A \emph{strong monoidal functor} between monoidal categories $(\bE,\otimes,\mathbb{1}_{\bE})$ and $(\mathbf{D},\otimes,\mathbb{1}_{\mathbf{D}})$ is a functor $F\colon\bE\to\mathbf{D}$ equipped with natural isomorphisms
$\varphi_{X,Y}\colon F(X)\otimes F(Y)\xrightarrow{\cong}F(X\otimes Y)$ and
$\varphi_0\colon\mathbb{1}_{\mathbf{D}}\xrightarrow{\cong}F(\mathbb{1}_{\bE})$,
satisfying the usual associativity and unit coherence conditions.

A strong monoidal functor between braided monoidal categories is \emph{braided monoidal} if
\[
\varphi_{Y,X}\cdot c_{F(X),F(Y)}
=
F(c_{X,Y})\cdot\varphi_{X,Y}
\]
for all objects $X,Y\in\bE$. A functor between braided, symmetric, or balanced monoidal categories is understood to be a strong monoidal functor preserving the corresponding additional structure.

\begin{definition}
Let $\bE=(\bE,\otimes,\mathbb{1})$ be a strict monoidal category.
A \emph{left dual} of an object $X\in\bE$ consists of an object $X^*$
together with morphisms
\[
\coev_X\colon\mathbb{1}\longrightarrow X\otimes X^*,
\qquad
\ev_X\colon X^*\otimes X\longrightarrow\mathbb{1},
\]
satisfying the zigzag identities
\begin{equation}
\label{eq: left dual zig zag equations}
(\id_X\otimes\ev_X)\cdot(\coev_X\otimes\id_X)=\id_X,
\qquad
(\ev_X\otimes\id_{X^*})\cdot(\id_{X^*}\otimes\coev_X)=\id_{X^*}.
\end{equation}
We say that $\bE$ \emph{has left duals} if every object of $\bE$ has
a chosen left dual.
\end{definition}

A left dual, if it exists, is unique up to a unique isomorphism
compatible with the evaluation and coevaluation morphisms. In a
strict monoidal category with left duals, the \emph{transpose} of a
morphism $f\colon X\to Y$ is the morphism $f^*\colon Y^*\to X^*$
defined by
\begin{equation}
\label{eq: transpose of a morphism}
f^*
\coloneqq
(\ev_Y\otimes\id_{X^*})
\cdot
(\id_{Y^*}\otimes f\otimes\id_{X^*})
\cdot
(\id_{Y^*}\otimes\coev_X).
\end{equation}

\begin{definition}
\label{def:ribbon category}
A \emph{strict ribbon category} is a strict balanced monoidal category
$(\bE,\otimes,\mathbb{1},c,\Theta)$ with chosen left duals such that
\[
\Theta_{X^*}=(\Theta_X)^*
\]
for every object $X\in\bE$. A \emph{ribbon monoidal functor} is a strong monoidal functor that
preserves the braiding and twist.
\end{definition}

\begin{example}
Let $\kk=\mathbb{C}(q)$, where $q$ is generic. The category of
finite-dimensional type-$1$ representations of
$U_q(\mathfrak{sl}_2)$ is a ribbon category. Its braiding is induced
by the universal $R$-matrix, its twist by the ribbon element, and
its duality by contragredient representations.
\end{example}

Let $\bE=(\bE,\otimes,\mathbb{1},c)$ be a strict symmetric monoidal
category. A \emph{non-degenerate symmetric form} on an object
$X\in\bE$ consists of morphisms
\[
\ev_X\colon X\otimes X\longrightarrow\mathbb{1}
\qquad \text{and} \qquad
\coev_X\colon\mathbb{1}\longrightarrow X\otimes X
\]
satisfying the zig-zag identities
\begin{equation}
\label{eq:self-dual-zig-zag}
(\id_X\otimes\ev_X)\cdot(\coev_X\otimes\id_X)=\id_X
\qquad \text{and} \qquad
(\ev_X\otimes\id_X)\cdot(\id_X\otimes\coev_X)=\id_X,
\end{equation}
and the symmetry conditions
\begin{equation}
\label{eq:symmetric-form}
\ev_X\cdot c_{X,X}=\ev_X
\qquad \text{and} \qquad
c_{X,X}\cdot\coev_X=\coev_X.
\end{equation}
The morphisms $\ev_X$ and $\coev_X$ are also called the \emph{cap}
and \emph{cup}, respectively.

Giving a non-degenerate form on $X$ is equivalent to giving a
self-duality $X\cong X^*$. Under this correspondence, the form is
symmetric precisely when its evaluation morphism
$\ev_X\colon X\otimes X\to\mathbb{1}$ satisfies
$\ev_X\cdot c_{X,X}=\ev_X$. See
\cite[Sections~2.1.1 and~3.1.2]{hv_cyclic} and
\cite{lukacs2010cyclic} for further discussion.
\begin{example}
In the cartesian symmetric monoidal category $(\Cat,\times,*)$, the
only dualizable object, up to isomorphism, is the terminal category
$*$; see \cite[Remark~2.15]{lukacs2010cyclic}. Consequently, the only
small category admitting a non-degenerate symmetric form is $*$,
up to isomorphism.
\end{example}

The preceding definitions extend to non-strict monoidal categories
by inserting the appropriate associators and unitors. For example,
let $(\bE,\otimes,\mathbb{1},\alpha,\lambda,\rho)$ be a symmetric
monoidal category and suppose that an object $X\in\bE$ is equipped
with morphisms
\[
\coev_X\colon\mathbb{1}\longrightarrow X\otimes X
\qquad \text{and} \qquad
\ev_X\colon X\otimes X\longrightarrow\mathbb{1}.
\]
The corresponding snake identities are
\begin{equation}
\label{eq:non-strict-snake-right}
\rho_X
\cdot(\id_X\otimes\ev_X)
\cdot\alpha_{X,X,X}
\cdot(\coev_X\otimes\id_X)
\cdot\lambda_X^{-1}
=
\id_X
\end{equation}
and
\begin{equation}
\label{eq:non-strict-snake-left}
\lambda_X
\cdot(\ev_X\otimes\id_X)
\cdot\alpha_{X,X,X}^{-1}
\cdot(\id_X\otimes\coev_X)
\cdot\rho_X^{-1}
=
\id_X.
\end{equation}
By Mac~Lane's Coherence Theorem, every monoidal category is monoidally equivalent to a strict monoidal category \cite{MacLane_1963}. Analogous strictification results hold for braided, balanced, and ribbon monoidal categories; see \cite{Joyal_Street_Braided_Tensor_Categories,shum_tortile_1994}.

\subsection{Homotopical preliminaries}
\label{subsec: Homotopical Preliminaries}

To discuss homotopy classes of maps between operads and cyclic operads, we use the language of Quillen model categories and $\infty$-categories. We refer to \cite{hirsch,balchin} for model categories and use relative categories as presentations of $\infty$-categories, following \cite{barwick_kan} and \cite[Chapter~13]{balchin}.

\begin{definition}\label{def: relative category} A \emph{relative category} is a pair $(\bE,\mathscr W)$, where $\bE$ is a category and $\mathscr W$ is a wide subcategory whose morphisms are called \emph{weak equivalences}. A functor of relative categories $F\colon(\bE,\mathscr W)\to(\bE',\mathscr W')$ is a functor $F\colon\bE\to\bE'$ such that $F(\mathscr W)\subseteq\mathscr W'$.
\end{definition}

Any model category $\bE$ determines a relative category $(\bE,\mathscr W)$ by taking $\mathscr W$ to be its class of weak equivalences. The homotopy category $\Ho(\bE)$ is obtained by formally inverting the morphisms in $\mathscr W$.

A relative category also presents an $\infty$-category; see \cite{barwick_kan}. We distinguish a relative category $\bE$ from the $\infty$-category $\mathcal E$ that it presents. A concrete model for the mapping spaces of $\mathcal E$ is provided by the Dwyer--Kan simplicial localisation $\mathcal L(\bE,\mathscr W)$. This is a simplicially enriched category with the same objects as $\bE$, equipped with a natural functor $\bE\to\mathcal L(\bE,\mathscr W)$. For $X,Y\in\ob(\bE)$, we write $\mathbb R\Map_{\bE}(X,Y)$ for its simplicial mapping space. This is the \emph{derived mapping space}, or Dwyer--Kan homotopy function complex, from $X$ to $Y$, and it satisfies
\[
\Hom_{\Ho(\bE)}(X,Y) \cong \pi_0\mathbb R\Map_{\bE}(X,Y).
\]

Thus the connected components of $\mathbb R\Map_{\bE}(X,Y)$ are the morphisms from $X$ to $Y$ in the homotopy category, while its higher simplices retain the higher homotopical information encoded by the relative category. In particular, a weak equivalence $X\to X'$ induces weak equivalences
\[
\mathbb R\Map_{\bE}(X',Y)
\xrightarrow{\simeq}
\mathbb R\Map_{\bE}(X,Y)
\qquad \text{and} \qquad
\mathbb R\Map_{\bE}(Y,X)
\xrightarrow{\simeq}
\mathbb R\Map_{\bE}(Y,X')
\]
for every $Y\in\bE$; see \cite{DWYER1980267}.

Although derived mapping spaces are generally difficult to compute, they admit a particularly simple description in a simplicial model
category.

\begin{thm}[{\cite[Corollary~4.7]{dk}}]
\label{thm: identifying mapping spaces in a simplicial model category} Let $\bE$ be a simplicial model category. If $X$ is cofibrant and $Y$ is fibrant in $\bE$, then there is a natural zigzag of weak equivalences
\[
\Map_{\bE}(X,Y) \simeq \mathbb R\Map_{\bE}(X,Y).
\]
\end{thm}

Given simplicial model categories $\bE$ and $\bD$, a simplicial Quillen adjunction
\[
\begin{tikzcd}
\bE \arrow[r, shift left=1, "F"] &
\arrow[l, shift left=1, "G"]
\bD
\end{tikzcd}
\]
induces an adjunction between the $\infty$-categories presented by their underlying relative categories, and hence a homotopy adjunction between their Dwyer--Kan simplicial localisations. Consequently, for $X\in\bE$ and $Y\in\bD$, there is a natural weak equivalence
\[
\mathbb R\Map_{\bD}(\mathbb L F(X),Y)
\simeq
\mathbb R\Map_{\bE}(X,\mathbb R G(Y)),
\]
where $\mathbb L F$ and $\mathbb R G$ denote the corresponding derived functors.

As a particular example, let $\sSet$ denote the category of simplicial sets with the Kan--Quillen model structure, and let $\Grpd$ denote the category of groupoids with the model structure in which weak equivalences are equivalences of categories, cofibrations are functors injective on objects, and fibrations are isofibrations; see \cite{And,Bous_model_cat_groupoid}. Every groupoid is both cofibrant and fibrant in this model structure.

The fundamental groupoid functor is left adjoint to the classifying space functor:
\begin{equation}
\label{classifying space adjunction}
\begin{tikzcd}
\Pi_1\colon\sSet \arrow[r, shift left=1] &
\arrow[l, shift left=1]
\Grpd\colon\cs .
\end{tikzcd}
\end{equation}
We equip $\Grpd$ with the corresponding simplicial enrichment, whose mapping spaces are $\Map_{\Grpd}(\bE,\bD)\coloneqq\Map_{\sSet}(\cs\bE,\cs\bD).$ The adjunction \eqref{classifying space adjunction} is a simplicial Quillen adjunction whose right adjoint is homotopically fully faithful. Thus, for every pair of groupoids $\bE$ and $\bD$, the natural map
\[
\mathbb R\Map_{\Grpd}(\bE,\bD)
\longrightarrow
\mathbb R\Map_{\sSet}(\cs\bE,\cs\bD)
\]
is a weak equivalence.

\section{Operads and cyclic operads}
\label{sec: Operads, Cyclic operads}

We review operads and cyclic operads in
Sections~\ref{subsec: operads} and~\ref{subsec: cyclic operads},
respectively. In Section~\ref{subsec: Envelop and Props}, we recall
the strict symmetric monoidal category associated to an operad. We
then introduce the metric prop associated to a cyclic operad in
Section~\ref{subsec: envelope of cyclic operad}.

\subsection{Operads}
\label{subsec: operads}

Throughout this section, let $\bE=(\bE,\otimes,\mathbb{1})$ be a closed symmetric monoidal category admitting all small colimits, and suppose that $\otimes$ preserves colimits separately in each variable. Let $\Sigma_n=\Aut(\{1,\ldots,n\})$, with $\Sigma_0$ the trivial group.

A \emph{symmetric sequence} in $\bE$ is a collection $\calO=\{\calO(n)\}_{n\geq 0}$ in which each $\calO(n)$ is equipped with a right action of $\Sigma_n$. Thus, for every $\sigma\in\Sigma_n$, there is an automorphism
\[
\sigma^*\colon\calO(n)\longrightarrow\calO(n)
\]
such that $\id^*=\id$ and $\sigma^*\cdot\tau^*=(\tau\sigma)^*$ for all $\sigma,\tau\in\Sigma_n$. Here, permutations are composed in the usual functional order: $(\tau\sigma)(i)=\tau(\sigma(i))$.

\begin{definition}
\label{def: operad}
An \emph{operad} in $\bE$ consists of a symmetric sequence
$\calO=\{\calO(n)\}_{n\geq 0}$ together with a unit morphism
$\eta\colon\mathbb{1}\to\calO(1)$ and partial composition morphisms
\[
\circ_i\colon\calO(n)\otimes\calO(m)
\longrightarrow\calO(n+m-1),
\qquad 1\leq i\leq n,
\]
satisfying the following axioms. We denote the unit operation determined
by $\eta$ by $1\in\calO(1)$ and state the axioms using generalized
elements.

For $x\in\calO(n)$, $y\in\calO(m)$, and $z\in\calO(\ell)$, nested
substitutions satisfy
\begin{equation}
\label{eq:operad-associativity-nested}
(x\circ_i y)\circ_{i+j-1}z
=
x\circ_i(y\circ_j z),
\qquad 1\leq i\leq n,\quad 1\leq j\leq m,
\end{equation}
while substitutions into distinct inputs satisfy
\begin{equation}
\label{eq:operad-associativity-disjoint}
(x\circ_i y)\circ_{k+m-1}z
=
(x\circ_k z)\circ_i y,
\qquad 1\leq i<k\leq n.
\end{equation}
The unit satisfies
\begin{equation}
\label{eq:operad-unit}
1\circ_1x=x 
\qquad \text{and} \qquad
x\circ_i1=x
\quad\text{for }1\leq i\leq n.
\end{equation}
Finally, for $\sigma\in\Sigma_n$ and $\tau\in\Sigma_m$, equivariance
is expressed by
\begin{equation}
\label{eq:operad-equivariance}
(\sigma^*x)\circ_i(\tau^*y)
=
(\sigma\circ_i\tau)^*
\bigl(x\circ_{\sigma(i)}y\bigr).
\end{equation}
Here $\sigma\circ_i\tau\in\Sigma_{n+m-1}$ is the block permutation
obtained by substituting the ordered $m$-element block acted on by
$\tau$ into the $i$th input of $\sigma$.
\end{definition}

A morphism of operads $f\colon\calO\to\calP$ is a morphism of the
underlying symmetric sequences that preserves the unit and partial
composition morphisms. We write $\Op(\bE)$ for the category of
operads in $\bE$. See \cite[Chapter~1]{FresseBook1} or
\cite{mms_operads_topology_and_phys} for further details.

\begin{remark}
Although we allow operads to have nullary operations, most of the
examples considered in this paper have an initial object in arity
zero. In particular, we use the positive-arity versions of the little
disks and framed little disks operads, with
$\mathsf D_2(0)=\mathsf{FD}_2(0)=\varnothing$.
\end{remark}

\begin{figure}[!ht]
\centering
\resizebox{0.65\textwidth}{!}{%
\begin{circuitikz}
\tikzstyle{every node}=[font=\small]
\draw [short] (6,17.5) -- (6,17);
\draw [short] (6,17.5) -- (5,18.5)node[pos=1,above, fill=white]{1};
\draw [short] (6,17.5) -- (5.5,18.5)node[pos=1,above, fill=white]{2};
\draw [short] (6,17.5) -- (6.5,18.5)node[pos=1,above, fill=white]{4};
\draw [short] (6,17.5) -- (7,18.5)node[pos=1,above, fill=white]{5};
\draw [short] (6,17.5) -- (6,18.5)node[pos=1,above, fill=white]{3};
\draw [short] (9.5,17.5) -- (9.5,17);
\draw [short] (9.5,17.5) -- (8.75,18.5)node[pos=1,above, fill=white]{1};
\draw [short] (9.5,17.5) -- (10.25,18.5)node[pos=1,above, fill=white]{4};
\draw [short] (9.5,17.5) -- (9.25,18.5)node[pos=1,above, fill=white]{2};
\draw [short] (9.5,17.5) -- (9.75,18.5)node[pos=1,above, fill=white]{3};
\draw [short] (13.25,17.25) -- (13.25,16.75);
\draw [short] (13.25,17.25) -- (12.25,18.25)node[pos=1,above, fill=white]{1};
\draw [short] (13.25,17.25) -- (12.75,18.25)node[pos=1,above, fill=white]{2};
\draw [short] (13.25,17.25) -- (13.25,18.25)node[pos=1,above, fill=white]{3};
\draw [short] (13.25,17.25) -- (13.75,18.25);
\draw [short] (13.25,17.25) -- (14.25,18.25)node[pos=1,above, fill=white]{8};
\draw [short] (13.75,18.25) -- (13.75,18.75);
\draw [short] (13.75,18.75) -- (13.25,19.25)node[pos=1,above, fill=white]{4};
\draw [short] (13.75,18.75) -- (13.5,19.25)node[pos=1,above, fill=white]{5};
\draw [short] (13.75,18.75) -- (14,19.25)node[pos=1,above, fill=white]{6};
\draw [short] (13.75,18.75) -- (14.25,19.25)node[pos=1,above, fill=white]{7};
\node [font=\large] at (8,17.75) {$\circ_4$};
\node [font=\large] at (11,17.75) {$=$};
\end{circuitikz}
}%
\caption{Partial composition represented by grafting rooted trees: $\circ_4\colon\calO(5)\otimes\calO(4)\longrightarrow\calO(8)$.}
\label{fig:operadic compostion}
\end{figure}
\begin{example}
\label{ex: endomorphism operad}
Fix an object $X\in\bE$. Since $\bE$ is closed, the objects
\[
\End(X)(n)\coloneqq
\underline{\Hom}_{\bE}(X^{\otimes n},X),
\qquad n\geq 0,
\] form a symmetric sequence in $\bE$. For $\sigma\in\Sigma_n$, the right action is given, in generalized-element notation, by
\begin{equation}
\label{eq:endomorphism-operad-action}
(\sigma^*f)(x_1,\ldots,x_n)
=
f(x_{\sigma^{-1}(1)},\ldots,x_{\sigma^{-1}(n)}).
\end{equation}
The partial compositions are given by substitution:\[(f\circ_i g)(x_1,\ldots,x_{n+m-1}) =f\bigl(x_1,\ldots,x_{i-1},g(x_i,\ldots,x_{i+m-1}),x_{i+m},\ldots,x_{n+m-1}\bigr).\] The identity morphism $\id_X\colon X\to X$ determines the operadic unit. The resulting operad $\End(X)$ is called the \emph{endomorphism operad} of $X$.
\end{example}

An \emph{algebra over an operad} $\calO$ is an object $A\in\bE$ together with a morphism of operads $\rho\colon\calO\to\End(A)$. By the tensor--hom adjunction, this is equivalent to specifying structure morphisms $\rho_n\colon\calO(n)\otimes A^{\otimes n}\to A$ for every $n\geq0$, compatible with the symmetric-group actions, partial compositions, and units. A morphism of $\calO$-algebras is a morphism in $\bE$ compatible with the $\calO$-actions. We write $\Alg_{\calO}(\bE)$ for the category of $\calO$-algebras in $\bE$.

\begin{example}
\label{example: D_2}
The \emph{little $2$-disks operad} $\mathsf{D}_2$ is an operad in
topological spaces with $\mathsf D_2(0)=*$. For $n\geq1$,
its space of $n$-ary operations $\mathsf{D}_2(n)$ consists of embeddings
\[
\coprod_{k=1}^{n}D_k^2\longrightarrow D^2
\]
whose restriction to each disk is of the form $x\mapsto r_kx+a_k$, where $r_k>0$, and whose images are pairwise disjoint and contained in the interior of the unit disk. Operadic composition is given by inserting one configuration of disks into a disk of another configuration.

For every $n\geq1$, the map that sends a configuration of disks to the ordered configuration of their centres is a homotopy equivalence $\mathsf{D}_2(n)\simeq\Conf_n(\mathbb{C})$, where $\Conf_n(\mathbb{C})=\{(z_1,\ldots,z_n)\in\mathbb{C}^n \mid z_i\neq z_j\text{ for }i\neq j\}$.
\end{example}

\begin{example}
\label{example: FD_2}
The \emph{framed little $2$-disks operad} $\mathsf{FD}_2$ is obtained
from $\mathsf{D}_2$ by equipping each little disk with a rotation.
Equivalently, it is the semidirect product
$\mathsf{FD}_2\cong\mathsf{D}_2\rtimes\mathsf{SO}(2)$; see
\cite[Example~1.3.1]{Wahl_Thesis} and
\cite{salvatore2003framed}. In particular,
$\mathsf{FD}_2(n)\cong\mathsf{D}_2(n)\times\mathsf{SO}(2)^n$.

Taking centres and rotation parameters gives a homotopy equivalence
\[
\mathsf{FD}_2(n)\simeq\Conf_n^{\mathrm{fr}}(\mathbb{C})
\coloneqq\Conf_n(\mathbb{C})\times(S^1)^n
\]
for every $n\geq1$. Each factor of $S^1$ records a framing, or equivalently a tangent direction, at the corresponding point.
\end{example}

\subsection{Cyclic operads}
\label{subsec: cyclic operads}

A cyclic operad is an operad whose symmetric-group actions extend
to allow the output of an operation to be permuted with its inputs.
Cyclic structures encode the compatibility conditions for invariant
symmetric forms on cyclic algebras; see
\cite[Section~4]{gk_cyc}.

For $n\geq0$, let $\Sigma_n^+=\Aut(\{0,1,\ldots,n\})$. We regard $\Sigma_n$ as the subgroup of $\Sigma_n^+$ fixing $0$, and write $z_{n+1}=(0\,1\,\cdots\,n)\in\Sigma_n^+$, so that $z_{n+1}(i)=i+1\pmod{n+1}$.

A \emph{cyclic structure} on an operad $\calO$ consists of a right $\Sigma_n^+$-action on $\calO(n)$ for every $n\geq0$, extending its given right $\Sigma_n$-action, such that $z_2^*(1)=1$ and
\begin{equation}
\label{eq:cyclic-composition}
z_{n+m}^*(x\circ_i y)
=
\begin{cases}
z_{n+1}^*(x)\circ_{i-1}y,
    & 2\leq i\leq n,\\
z_{m+1}^*(y)\circ_m z_{n+1}^*(x),
    & i=1,
\end{cases}
\end{equation}
for $x\in\calO(n)$, $y\in\calO(m)$, and $n,m\geq1$.  When $m=0$, the second case of
Equation~\eqref{eq:cyclic-composition} is replaced by $z_n^*(x\circ_1y)=\bigl((z_{n+1}^2)^*x\bigr)\circ_ny$ for $x\in\calO(n)$ and $y\in\calO(0)$. Equation~\eqref{eq:cyclic-composition}
expresses the compatibility between cyclic symmetry and operadic
composition.

\begin{figure}[!ht]
\centering
\resizebox{0.5\textwidth}{!}{%
\begin{circuitikz}
\tikzstyle{every node}=[font=\small]
\draw [short] (5.5,20) -- (6,19)node[pos=0,above, fill=white]{2};
\draw [short] (6,19) -- (6,20)node[pos=1,above, fill=white]{3};
\draw [short] (6,19) -- (6.5,20)node[pos=1,above, fill=white]{4};
\draw [short] (6,19) -- (7,20)node[pos=1,above, fill=white]{5};
\draw [short] (6,19) -- (5,20)node[pos=1,above, fill=white]{1};
\draw [short] (6,19) -- (6,18)node[pos=0.75,below, fill=white]{0};
\draw [short] (10.5,19) -- (10.5,20)node[pos=1,above, fill=white]{3};
\draw [short] (10.5,19) -- (11,20)node[pos=1,above, fill=white]{4};
\draw [short] (10.5,19) -- (10,20)node[pos=1,above, fill=white]{2};
\draw [short] (10.5,19) -- (11.5,20)node[pos=1,above, fill=white]{5};
\draw [short] (10.5,19) .. controls (11,17.25) and (11.5,19) .. (12,20)node[pos=1,above, fill=white]{0};
\draw [short] (10.5,19) .. controls (8.75,21.25) and (9.25,18.75) .. (10.5,18)node[pos=1,right, fill=white]{1};
\node [font=\normalsize] at (8,19) {$\overset{z^*}{\longrightarrow}$};
\end{circuitikz}
}%
\caption{The cyclic action moves the first input to the output and
the former output to the final input.}\label{fig:cyclic action}
\end{figure}

\begin{definition}
\label{def: cyclic operad}
A \emph{cyclic operad} is an operad $\calO=\{\calO(n)\}_{n\geq0}$ equipped with a cyclic structure. A morphism of cyclic operads $f\colon\calO\to\calP$ is a morphism of the underlying operads that commutes with the $\Sigma_n^+$-actions. We denote the category of cyclic operads in $\bE$ by $\Cyc(\bE)$.
\end{definition}

\begin{example} \label{example: cyclic endomorphism}
Let $X\in\bE$ be equipped with a non-degenerate symmetric form $(\ev_X,\coev_X)$. The form induces natural isomorphisms
\[
\underline{\Hom}_{\bE}(X^{\otimes n},X)
\cong
\underline{\Hom}_{\bE}(X^{\otimes(n+1)},\mathbb{1})
\] that send an element $f\colon X^{\otimes n}\to X$ to
\[
\widetilde f
\coloneqq
\ev_X\cdot(f\otimes\id_X)
\colon X^{\otimes(n+1)}\longrightarrow\mathbb{1}.
\]
The natural right action of $\Sigma_n^+$ on $\underline{\Hom}_{\bE}(X^{\otimes(n+1)},\mathbb{1})$, given by permuting the $n+1$ tensor factors, therefore induces a right $\Sigma_n^+$-action on $\End(X)(n)$. These actions equip $\End(X)$ with a cyclic structure. The resulting cyclic operad is called the \emph{cyclic endomorphism operad} of $X$ and is denoted $\End^{\mathrm{cyc}}(X)$.
\end{example}

\begin{definition}
\label{def:cyclic-algebra}
Let $\calO$ be a cyclic operad in $\bE$. A \emph{cyclic
$\calO$-algebra} consists of an object $X\in\bE$ equipped with a
non-degenerate symmetric form and a morphism of cyclic operads
\[
\calO\longrightarrow\End^{\mathrm{cyc}}(X).
\]
Equivalently, it consists of an $\calO$-algebra structure
$\rho_n\colon\calO(n)\otimes X^{\otimes n}\to X$ such that the maps
\begin{equation}
\label{eq:cyclic-algebra-invariance}
\ev_X\cdot(\rho_n\otimes\id_X)
\colon
\calO(n)\otimes X^{\otimes(n+1)}
\longrightarrow
\mathbb{1}
\end{equation}
are invariant under the diagonal $\Sigma_n^+$-actions for every
$n\geq0$.

A morphism of cyclic $\calO$-algebras is a morphism of the underlying
$\calO$-algebras that preserves the chosen forms. We write
$\Alg_{\calO}^{\mathrm{cyc}}(\bE)$ for the category of cyclic
$\calO$-algebras.
\end{definition}

\begin{remark}
The cyclic unit axiom $z_{2}^*(1)=1$ ensures that the invariant bilinear form (in Definition~\ref{def:cyclic-algebra}) associated to an algebra over a cyclic operad is \emph{symmetric}. Therefore, cyclic operads model algebras equipped with an invariant \emph{symmetric} pairing, see Example~\ref{example: cyclic endomorphism}.
\end{remark}

\begin{example}
\label{ex: Budney}
Following \cite[Section~5.2]{cyclic_ribbon}, we recall a topological
cyclic operad $\mathsf{S}_0$ modelling genus-zero surfaces with
parametrised boundary components. This construction is closely
related to Budney's conformal-ball model for the framed little-disks
operad~\cite{Budney}.

For $n\geq 1$, let $\operatorname{Adm}(n)$ be the space of admissible maps

\[
a=\coprod_{j=0}^n a_j\colon
\coprod_{j=0}^n D_j^2\longrightarrow S^2.
\]
Here each $a_j$ is an orientation-preserving smooth embedding, the
images of the interiors of the disks are pairwise disjoint, and the
antipodal point $-a_j(0)$ does not lie in the image of $a_j$. We also
require that, after stereographic projection from $-a_j(0)$ to the
tangent plane at $a_j(0)$, the map $a_j$ is given by a positive
rescaling. The direction determined by $a_j(1,0)$ fixes
the corresponding identification of the tangent plane with
$\mathbb{R}^2$.

The group $\mathrm{SO}(3)$ acts on $\operatorname{Adm}(n)$ by
postcomposition, and we define
\[
\mathsf{S}_0(n)\coloneqq
\operatorname{Adm}(n)/\mathrm{SO}(3),
\]
equipped with the quotient topology induced by the compact-open
topology. The group $\Sigma_n^+$ acts on $\mathsf{S}_0(n)$ by
permuting the labels $0,\ldots,n$.

Given $a\in\mathsf{S}_0(n)$ and $b\in\mathsf{S}_0(m)$, the partial
composition $a\circ_i b\in\mathsf{S}_0(n+m-1)$ is defined as follows.
Using stereographic projection from $b_0(0)$, the complement of the
zeroth disk is identified with a planar disk containing the remaining
disk configuration. This configuration is then inserted into the
$i$th disk of $a$. Together with the $\Sigma_n^+$-actions, these
compositions make $\mathsf{S}_0$ a topological cyclic operad.

There is a natural map of operads
$\mathsf{FD}_2\to\mathsf{S}_0$ obtained by mapping a framed
configuration of disks into the lower hemisphere of $S^2$ using
inverse stereographic projection and adjoining the upper hemisphere
as the disk labelled $0$. This map is an aritywise weak homotopy
equivalence and hence a weak equivalence of operads
\[
\mathsf{FD}_2\xrightarrow{\sim}\mathsf{S}_0.
\]
Thus $\mathsf{S}_0$, with its $\Sigma_n^+$-actions, provides a cyclic
model for the framed little $2$-disks operad. We write
$\mathsf{FD}_2^{\mathrm{cyc}}$ for this cyclic model when we wish to
emphasise its cyclic structure.

We may compare this cyclic structure directly with Budney's
conformal-ball model. In that model, exchanging a chosen input with
the output is implemented by reparametrising the configuration by the
inverse of the corresponding conformal map, followed by a fixed
half-turn of $S^2$. Passing to framed point configurations by
stereographic projection, and normalising the point labelled $1$ to
lie at the origin, identifies the transposition $(01)$ with the
involution $z\longmapsto-1/z$.
\end{example}

\subsection{The envelope of an operad}
\label{subsec: Envelop and Props}

A \emph{prop} in $\bE$ is a strict symmetric monoidal
$\bE$-enriched category whose objects are the natural numbers
$\mathbb{N}=\{0,1,2,\ldots\}$ and whose monoidal product on objects
is addition. The hom-object $\calP(m,n)$ is interpreted as the object
of operations with $m$ inputs and $n$ outputs.

Every operad generates a prop, called its \emph{envelope}, through the universal construction described below.

\begin{definition}
\label{def: envelope}
Let $\calO$ be an operad in $\bE$. Its \emph{envelope} $\Env(\calO)$ is the prop whose hom-objects are
\begin{equation}
\label{eq:envelope}
\Env(\calO)(m,n)
=
\coprod_{m_1+\cdots+m_n=m}
\left(
\calO(m_1)\otimes\cdots\otimes\calO(m_n)
\right)
\underset{\Sigma_{m_1}\times\cdots\times\Sigma_{m_n}}{\otimes}
\mathbb{1}[\Sigma_m].
\end{equation}
Here $\mathbb{1}[\Sigma_m]\coloneqq \coprod_{\sigma\in\Sigma_m}\mathbb{1}$ is the coproduct of copies of the monoidal unit indexed by $\Sigma_m$. The group
$\Sigma_{m_1}\times\cdots\times\Sigma_{m_n}$ acts on the first factor through the right actions on the objects $\calO(m_i)$ and on $\mathbb{1}[\Sigma_m]$ by left multiplication through the canonical block inclusion
\[
\Sigma_{m_1}\times\cdots\times\Sigma_{m_n}
\hookrightarrow\Sigma_m.
\]
The tensor product over
$\Sigma_{m_1}\times\cdots\times\Sigma_{m_n}$ denotes the corresponding coequalizer.
\end{definition}

More explicitly, a generalized element of
$\Env(\calO)(m,n)$ is represented by data
$\bigl((f_1,\ldots,f_n),\sigma\bigr)$, where
$m=m_1+\cdots+m_n$, each $f_i\in\calO(m_i)$, and
$\sigma\in\Sigma_m$. Two such representatives are identified under
the relation
\[
\bigl((\sigma_1^*f_1,\ldots,\sigma_n^*f_n),\sigma\bigr)
\sim
\bigl((f_1,\ldots,f_n),
(\sigma_1\oplus\cdots\oplus\sigma_n)\sigma\bigr)
\]
for $\sigma_i\in\Sigma_{m_i}$. Thus $\sigma$ records the ordering of
the $m$ inputs, while $f_i$ labels the operation associated to the
$i$th output.

Equivalently, writing $\underline m=\{1,\ldots,m\}$, with
$\underline0=\varnothing$, there is an isomorphism
\[
\Env(\calO)(m,n)
\cong
\coprod_{\phi\colon\underline m\to\underline n}
\bigotimes_{j=1}^n
\calO\bigl(|\phi^{-1}(j)|\bigr).
\]
Here the ordering inherited from $\underline m$ is used to identify each fibre with an ordered set, and empty fibres are labelled by nullary operations in $\calO(0)$. In particular, $\Env(\calO)(0,0)=\mathbb{1}$, while $\Env(\calO)(m,0)$ is initial for $m>0$.

Categorical composition is induced by operadic substitution. Namely, if a morphism from $m$ to $n$ is represented by operations
$f_1,\ldots,f_n$ and a morphism from $n$ to $p$ is represented by operations $g_1,\ldots,g_p$, their composite is obtained by substituting the appropriate operations $f_i$ into each $g_j$, using the permutation data to order the resulting inputs.

The monoidal product is given by concatenation. If $\bigl((f_1,\ldots,f_n),\sigma\bigr)$ represents a morphism from $m$ to $n$ and $\bigl((g_1,\ldots,g_q),\tau\bigr)$ represents a morphism from $p$ to $q$, their monoidal product is represented by
\[
\bigl((f_1,\ldots,f_n,g_1,\ldots,g_q),\sigma\oplus\tau\bigr)
\in\Env(\calO)(m+p,n+q),
\]
where $\sigma\oplus\tau$ is the image of $(\sigma,\tau)$ under the canonical block inclusion $\Sigma_m\times\Sigma_p\hookrightarrow\Sigma_{m+p}$.

The envelope construction defines a left adjoint $\Env\colon\Op(\bE)\rightleftarrows\Prop(\bE):u$ to the functor $u$, which sends a prop $\calP$ to its underlying operad $u(\calP)(n)=\calP(n,1)$. See \cite[Section~2.6.1]{hv_cyclic}, \cite[Proposition~11]{hr15}, or \cite{phillips2022operads} for further details.

An \emph{algebra over a prop} $\calP$ is a strong symmetric monoidal functor $F\colon\calP\to\bE$. Since $F(n)\cong F(1)^{\otimes n}$, such a functor determines a collection of operations on the object $F(1)$. We write $\Alg_{\bE}(\calP)$ for the category of $\calP$-algebras and monoidal natural transformations.

\begin{prop} \label{prop: algebras over operads and props are the same}
For every operad $\calO$ in $\bE$, there is a natural equivalence of categories
\[
\Alg_{\bE}(\Env(\calO))
\simeq
\Alg_{\bE}(\calO).
\]
\end{prop}

\subsection{The envelope of a cyclic operad}
\label{subsec: envelope of cyclic operad}

The envelope construction can be extended to encode invariant non-degenerate symmetric forms. The resulting prop is called the \emph{metric prop} of a cyclic operad \cite[Section~3.1.3]{hv_cyclic}.

\begin{definition}
\label{def: metric prop}
Let $\calO$ be a cyclic operad. The \emph{metric prop} associated
to $\calO$, denoted $\Pi(\calO)$, is obtained from $\Env(\calO)$
by adjoining morphisms $\ev\colon2\to0$ and $\coev\colon0\to2$,
subject to the following relations.

The evaluation and coevaluation satisfy the snake identities
\begin{equation}
\label{eq:metric-prop-snake}
(\id_1\otimes\ev)\cdot(\coev\otimes\id_1)=\id_1
\qquad\text{and}\qquad
(\ev\otimes\id_1)\cdot(\id_1\otimes\coev)=\id_1,
\end{equation}
and are symmetric:
\begin{equation}
\label{eq:metric-prop-symmetry}
\ev\cdot c_{1,1}=\ev
\qquad\text{and}\qquad
c_{1,1}\cdot\coev=\coev,
\end{equation}
where $c_{1,1}\colon1\otimes1\to1\otimes1$ is the symmetry
in the prop.

Finally, for every $f\in\calO(n)=\Env(\calO)(n,1)$, with $n\geq1$,
the cyclic action is encoded by
\begin{equation}
\label{eq:metric-prop-cyclic-relation}
z_{n+1}^*(f)
=(\id_1\otimes\ev)
 \cdot(\id_1\otimes f\otimes\id_1)
 \cdot(\coev\otimes\id_n).
\end{equation}
Thus the right-hand side is the composite
\[
n
\xrightarrow{\coev\otimes\id_n}
1\otimes1\otimes n
\xrightarrow{\id_1\otimes f\otimes\id_1}
1\otimes1\otimes1
\xrightarrow{\id_1\otimes\ev}
1.
\]
The second leg of the coevaluation enters the first input of $f$,
followed by the first $n-1$ original inputs. The output of $f$ is
paired with the last original input, while the first leg of the
coevaluation becomes the new output.
\end{definition}

An algebra over $\Pi(\calO)$ is a strong symmetric monoidal functor $F\colon\Pi(\calO)\to\bE$. The morphisms $F(\ev)$ and $F(\coev)$ equip $F(1)$ with a non-degenerate symmetric form, while Equation~\eqref{eq:metric-prop-cyclic-relation} expresses the invariance of this form under the $\calO$-algebra structure. We write $\Alg_{\bE}(\Pi(\calO))$ for the category of such algebras.

\begin{lemma}[{\cite[Lemma~3.1.4]{hv_cyclic}}]
\label{lem:algebras-over-metric-prop}
For every cyclic operad $\calO$ in $\bE$, there is a natural equivalence of categories
\[
\Alg_{\bE}(\Pi(\calO))
\simeq
\Alg_{\bE}^{\mathrm{cyc}}(\calO),
\]
where $\Alg_{\bE}^{\mathrm{cyc}}(\calO)$ denotes the category of $\calO$-algebras equipped with an invariant non-degenerate symmetric form.
\end{lemma}

\begin{remark}
\label{rem:Hinich-Vaintrob-cyclic-convention}
Our convention for the cyclic generator differs from that of Hinich and Vaintrob~\cite{hv_cyclic}. We use the forward cycle $z_{n+1}(i)=i+1\pmod{n+1}$, whereas they use the inverse cycle
$\tau$, defined by $\tau(0)=n$ and $\tau(i)=i-1$ for $i>0$. Thus $\tau=z_{n+1}^{-1}$, and their cyclic invariance relation is expressed in our notation using $(z_{n+1}^{-1})^*$.
\end{remark}

\section{Cyclic \texorpdfstring{$\infty$}{infinity}-operads}
\label{sec: cyclic infinity operad}

To study homotopy automorphisms of profinitely completed operads and
cyclic operads, we work with Segal models for $\infty$-operads and
cyclic $\infty$-operads \cite{cisinski2013dendroidal,hry1,doherty2025models}. These models
apply to operads and cyclic operads in both groupoids and spaces and
allow us to compare their derived mapping spaces after profinite
completion.

We begin by recalling the combinatorics of rooted and unrooted trees
and the constructions of free operads and free cyclic operads. We
then introduce the corresponding Segal models for $\infty$-operads
and cyclic $\infty$-operads. Finally, we construct the cyclic
$\infty$-operad of framed configurations used below.

\subsection{Trees and the homotopy theory of (cyclic) operads}
\label{subsec:trees}

A \emph{graph} $G$ consists of a set of vertices $V(G)$, a set of
half-edges $H(G)$, an incidence map $s\colon H(G)\to V(G)$, and an
involution $i\colon H(G)\to H(G)$. The orbits of $i$ are the
\emph{edges} of $G$. An orbit $\{h,i(h)\}$ with $i(h)\neq h$ is an
\emph{internal edge}, while a fixed point of $i$ is an
\emph{external edge}, or \emph{boundary edge}. We write
$\mathrm{iE}(G)$ and $\partial G$ for the sets of internal and
boundary edges, respectively.

\begin{definition}
A \emph{tree} is a connected graph whose geometric realization is
contractible. A \emph{rooted tree} is a tree $T$ equipped with a
distinguished boundary edge $r(T)\in\partial T$, called its
\emph{root}. The remaining boundary edges are its \emph{leaves}, and
we write
\[
l(T)=\partial T\setminus\{r(T)\}.
\]
An \emph{unrooted tree} is a tree with no distinguished boundary
edge.
\end{definition}

Let $W\subseteq V(T)$ be a nonempty set of vertices such that the
full subgraph spanned by $W$ is connected. The \emph{subtree spanned
by $W$} is obtained by retaining the vertices in $W$ and all
half-edges incident to them, and by cutting every internal edge whose
other endpoint does not lie in $W$. More precisely, its half-edge set
is $s^{-1}(W)$, and its involution is given by
\[
i_W(h)=
\begin{cases}
i(h), & i(h)\in s^{-1}(W),\\
h,    & i(h)\notin s^{-1}(W).
\end{cases}
\]
Thus cutting an edge creates a new boundary edge without changing the
valence of any retained vertex. A tree is called a \emph{subtree} of
$T$ if it is isomorphic to one obtained in this way. We also formally
adjoin the \emph{exceptional tree} $\eta$, consisting of a single edge
with two boundary flags and no vertices.

\begin{definition}
A \emph{planar tree} is a tree endowed with an embedding in the plane,
and hence with a cyclic ordering of the half-edges incident to each
vertex. A \emph{labelled planar tree} is a triple $(T,\sigma,\tau)$,
where $T$ is a planar tree equipped with bijections
\[
\sigma\colon\{1,\dots,|V(T)|\}\xrightarrow{\cong}V(T)
\qquad\text{and}\qquad
\tau\colon\{0,1,\dots,|\partial T|-1\}
\xrightarrow{\cong}\partial T.
\]
Two such triples $(T,\sigma,\tau)$ and $(T',\sigma',\tau')$ are
\emph{isomorphic} if there is an isomorphism of planar trees
$T\to T'$ respecting the vertex and boundary labellings.
\end{definition}
Examples include the exceptional tree $\eta$, which has no vertices,
and the \emph{$n$-corolla} $C_{n+1}$, the unrooted
tree with one vertex and $n+1$ boundary edges labelled by
$\{0,1,\ldots,n\}$. We write $C_n$ for the corresponding rooted
corolla with root labelled by $0$ and leaves labelled by
$\{1,\ldots,n\}$.

\begin{remark}
Our notation distinguishes rooted and unrooted corollas by recording
the number of inputs. The rooted corolla representing an $n$-ary
operation is denoted $C_n$, while the corresponding unrooted corolla
is denoted $C_{n+1}$, since the output is included among its $n+1$
boundary edges. The context will always make clear whether a corolla
is rooted or unrooted.
\end{remark}

Labelled trees encode the operations governing cyclic operads. We
will construct an $\mathbb N$-coloured operad whose algebras are
precisely cyclic operads. We first recall the definition of a
coloured operad.

\begin{definition}
\label{def: coloured operad}
Let $\mathfrak C$ be a nonempty set of colours. A
\emph{$\mathfrak C$-coloured operad} $\mathcal P$ in $\bE$ consists
of objects
$\mathcal P(c_1,\ldots,c_n;c_0)$, for $n\geq0$ and
$c_0,c_1,\ldots,c_n\in\mathfrak C$, together with right
symmetric-group actions
\[
\sigma^*\colon
\mathcal P(c_1,\ldots,c_n;c_0)
\longrightarrow
\mathcal P(c_{\sigma(1)},\ldots,c_{\sigma(n)};c_0),
\]
partial composition morphisms
\[
\circ_i\colon
\mathcal P(c_1,\ldots,c_n;c_0)\otimes
\mathcal P(d_1,\ldots,d_m;c_i)
\longrightarrow
\mathcal P(c_1,\ldots,c_{i-1},
d_1,\ldots,d_m,c_{i+1},\ldots,c_n;c_0),
\]
and unit morphisms
$\mathbb 1\to\mathcal P(c;c)$ for every $c\in\mathfrak C$.
These data satisfy the usual associativity, equivariance, and unit
axioms in their coloured forms. See
\cite[Definition~1.1]{bm_resolutions}.
\end{definition}

\begin{remark}
When $\mathfrak C$ has one element, a $\mathfrak C$-coloured operad
is an operad in the sense of Definition~\ref{def: operad}.
\end{remark}

An \emph{algebra} over a $\mathfrak C$-coloured operad $\mathcal P$
consists of an object $A(c)\in\bE$ for every $c\in\mathfrak C$,
together with structure morphisms
\[
\mathcal P(c_1,\ldots,c_k;c_0)\otimes
A(c_1)\otimes\cdots\otimes A(c_k)
\longrightarrow A(c_0)
\]
satisfying the corresponding associativity, equivariance, and unit
axioms. See \cite[Definition~1.2]{bm_resolutions}. We write
$\Alg_{\bE}(\mathcal P)$ for the category of $\mathcal P$-algebras
in $\bE$.
\begin{example}
\label{ex: operad for cyclic operads}
The \emph{operad governing cyclic operads}, denoted $\mathbf C$, is
the discrete $\mathbb N$-coloured operad defined as follows. An
operation in $\mathbf C(m_1,\ldots,m_k;n)$ is represented by an
unrooted tree $T$ with $n+1$ boundary edges and $k$ vertices. It is
equipped with a boundary labelling
\[
\tau\colon\{0,\ldots,n\}\xrightarrow{\cong}\partial T,
\]
a vertex labelling
\[
\sigma\colon\{1,\ldots,k\}\xrightarrow{\cong}V(T),
\]
and, for every $i$, a bijection
\[
\lambda_i\colon
\{0,\ldots,m_i\}
\xrightarrow{\cong}
\operatorname{nb}(\sigma(i)).
\]
Here $\operatorname{nb}(v)$ is the set of half-edges incident to
$v$. In particular, the vertex $\sigma(i)$ has valence $m_i+1$.
Operations are isomorphism classes of trees equipped with these
labellings. The right action of $\Sigma_k$ is induced by relabelling
the vertices.

Let $T$ and $T'$ represent operations in
$\mathbf C(m_1,\ldots,m_k;n)$ and
$\mathbf C(b_1,\ldots,b_\ell;m_i)$, respectively. Their partial
composite $T\circ_iT'$ is obtained by substituting $T'$ for the
vertex $\sigma(i)$ of $T$. For every $0\leq j\leq m_i$, the boundary
edge of $T'$ labelled by $j$ is joined to the half-edge
$\lambda_i(j)$ incident to $\sigma(i)$. The resulting tree represents
an operation in
\[
\mathbf C
(m_1,\ldots,m_{i-1},b_1,\ldots,b_\ell,
m_{i+1},\ldots,m_k;n).
\]
Its boundary labelling is inherited from $T$. The local half-edge
labellings are inherited from the remaining vertices of $T$ and the
vertices of $T'$. Its vertex labelling is obtained by replacing the
$i$th vertex of $T$ by the ordered list of vertices of $T'$.

For every $n\geq0$, the identity operation in $\mathbf C(n;n)$ is
represented by the unrooted corolla with $n+1$ boundary edges, whose
boundary and local half-edge labellings agree. The exceptional tree
represents a nullary operation in $\mathbf C({};1)$. In a
$\mathbf C$-algebra, this operation determines the unit in arity
one. With these operations, $\mathbf C$-algebras are precisely cyclic
operads.
\end{example}
There is similarly an $\mathbb N$-coloured operad $\mathbf O$
governing operads. An operation in
$\mathbf O(m_1,\ldots,m_k;n)$ is represented by a planar rooted tree
with $n$ labelled leaves and $k$ labelled vertices. The vertex
labelled by $i$ has $m_i$ incoming edges. Operations are isomorphism
classes of trees equipped with these labellings. The right action of
$\Sigma_k$ is induced by relabelling the vertices.

Partial composition is given by substitution at a vertex. To form
$T\circ_iT'$, the root of $T'$ is attached to the outgoing edge of
the vertex labelled by $i$ in $T$. The labelled leaves of $T'$ are
attached to the corresponding incoming edges at that vertex. See
\cite[Example~1.5.6]{bm_resolutions} and
\cite[Definition~2.9]{cacti} for further details.

\begin{lemma}
\label{lemma: cyclic operads are algebras over operads}
Let $\bE$ be a cocomplete closed symmetric monoidal category. Regard
the discrete coloured operads $\mathbf C$ and $\mathbf O$ as coloured
operads in $\bE$ by taking copowers of the monoidal unit. There are
natural isomorphisms
\[
\Alg_{\bE}(\mathbf C)\cong\Cyc(\bE)
\qquad\text{and}\qquad
\Alg_{\bE}(\mathbf O)\cong\Op(\bE).
\]
\end{lemma}

Under suitable hypotheses on $\bE$, the model structure on $\bE$
transfers to categories of operads and cyclic operads. Applying
\cite[Theorem~2.1]{bm_resolutions} to the coloured operads
$\mathbf O$ and $\mathbf C$ gives the following.

\begin{prop}
\label{model category on alg over operad}
Let $\bE$ be a cofibrantly generated symmetric monoidal model
category with cofibrant unit, a symmetric monoidal fibrant
replacement functor, and a cocommutative coalgebra interval. Then
$\Op(\bE)$ and $\Cyc(\bE)$ admit cofibrantly generated model
structures in which a morphism
$f\colon\mathcal P\to\mathcal Q$ is a weak equivalence or fibration
if and only if every map
$f_n\colon\mathcal P(n)\to\mathcal Q(n)$ is a weak equivalence or
fibration, respectively, in $\bE$.
\end{prop}

\begin{remark}
The hypotheses of
Proposition~\ref{model category on alg over operad} hold for
$\sSet$ with the Kan--Quillen model structure and for $\Grpd$ with
its standard model structure \cite{And}. We therefore obtain the
corresponding transferred model structures on operads and cyclic
operads in both categories.
\end{remark}

\subsubsection{Adjunctions between operads and cyclic operads}

Every planar unrooted tree $T$ may be regarded as a rooted tree by
choosing a boundary edge $r\in\partial T$ as its root and orienting
all internal edges towards $r$. Conversely, forgetting the root of a
planar rooted tree produces an unrooted tree. In our diagrams, the
chosen root is the boundary edge labelled by $0$.

Let $\mathbf O$ and $\mathbf C$ denote the coloured operads governing
operads and cyclic operads. See
\cite[Example~1.5.6]{bm_resolutions} and
\cite[Definition~2.9]{cacti}. Forgetting the root defines a morphism
of coloured operads
$u\colon\mathbf O\to\mathbf C$. Restriction along $u$ gives the
forgetful functor
$u^*\colon\Cyc(\bE)\to\Op(\bE)$. It has a left adjoint
$u_!\colon\Op(\bE)\to\Cyc(\bE)$, called the \emph{cyclic envelope}.
Thus there is an adjunction
\[
u_!\colon\Op(\bE)
\rightleftarrows
\Cyc(\bE):u^*.
\]

Informally, $u_!\mathcal O$ freely adjoins the operations needed to
regard any boundary edge of an operation in $\mathcal O$ as its
output, subject to the cyclic-operad relations. An explicit
description is given in
\cite[Section~3.1 and Proposition~3.4]{DCH}. See also
\cite[Section~9]{Ward}.

Whenever the model structures on operads and cyclic operads are
transferred from their underlying collections, $u^*$ preserves
fibrations and trivial fibrations. Hence $u_!\dashv u^*$ is a
Quillen adjunction.
\subsubsection{Change of base for operads and cyclic operads}

The coloured operads $\mathbf O$ and $\mathbf C$ governing operads
and cyclic operads are $\Sigma$-cofibrant. Indeed, the
symmetric-group actions given by relabelling vertices are free. It
follows that a weak symmetric monoidal Quillen adjunction between
suitable symmetric monoidal model categories in which every object
is cofibrant induces Quillen adjunctions on operads and cyclic
operads. See
\cite[Proposition~5.2 and Theorem~5.8]{hry_shrink}.

We apply this result to the simplicial Quillen adjunction
\[
\Pi_1\colon\sSet
\rightleftarrows
\Grpd:\cs.
\]
Every object of $\sSet$ and $\Grpd$ is cofibrant, and both $\Pi_1$
and $\cs$ preserve finite products. The adjunction therefore induces
simplicial Quillen adjunctions
\[
\Pi_1\colon\Op(\sSet)
\rightleftarrows
\Op(\Grpd):\cs
\qquad\text{and}\qquad
\Pi_1\colon\Cyc(\sSet)
\rightleftarrows
\Cyc(\Grpd):\cs.
\]

\begin{lemma}
\label{lemma: B is homotopically fully faithful}
The classifying-space functor
$\cs\colon\Cyc(\Grpd)\to\Cyc(\sSet)$ is homotopically fully
faithful. Thus, for all $\calP,\calQ\in\Cyc(\Grpd)$, the natural map
\[
\mathbb R\Map_{\Cyc(\Grpd)}(\calP,\calQ)
\longrightarrow
\mathbb R\Map_{\Cyc(\sSet)}(\cs\calP,\cs\calQ)
\]
is a weak equivalence.
\end{lemma}

\begin{proof}
It is enough to show that the derived counit
$\mathbb L\Pi_1\,\mathbb R\cs(\calQ)\to\calQ$ is a weak equivalence
for every $\calQ\in\Cyc(\Grpd)$.

Every cyclic operad in groupoids is fibrant. Choose a cofibrant
replacement $(\cs\calQ)_c\to\cs\calQ$ in $\Cyc(\sSet)$. The derived
counit is represented by
\[
\Pi_1(\cs\calQ)_c
\longrightarrow
\Pi_1\cs\calQ
\longrightarrow
\calQ.
\]
The second map is the ordinary counit and is an isomorphism. The
first map is an entrywise equivalence of groupoids. Indeed,
$(\cs\calQ)_c(n)\to\cs\calQ(n)$ is a weak equivalence of simplicial sets,
and $\Pi_1$ preserves this weak equivalence because it is left
Quillen and every simplicial set is cofibrant. Hence the derived
counit is a weak equivalence.
\end{proof}

\subsection{Cofibrant cyclic operads in groupoids}
\label{subsec: cofibrancy in grpd}

We equip $\Cyc(\Grpd)$ with the model structure transferred from
cyclic collections. Weak equivalences and fibrations are therefore
defined entrywise. Since every groupoid is fibrant, every cyclic
operad in groupoids is fibrant. The following is the cyclic analogue
of Horel's description of cofibrant operads in groupoids. See
\cite[Proposition~6.8]{Horel_profinite_groupoids}.

Let $\ob\colon\Grpd\to\Set$ be the object functor. It has left and
right adjoints
\[
\operatorname{Disc}\dashv\ob\dashv\operatorname{Codisc}.
\]
These functors preserve finite products, so the adjunctions extend
entrywise to cyclic operads.

\begin{prop}
\label{prop:cofibrations-cyc}
A morphism $g\colon\mathcal P\to\mathcal Q$ in $\Cyc(\Grpd)$ is a
cofibration if and only if
$\ob(g)\colon\ob(\mathcal P)\to\ob(\mathcal Q)$ has the left lifting
property with respect to every entrywise surjective morphism in
$\Cyc(\Set)$.
\end{prop}

\begin{proof}
Suppose first that $g$ is a cofibration, and let
$q\colon X\to Y$ be an entrywise surjective morphism in
$\Cyc(\Set)$. The morphism $\operatorname{Codisc}(q)$ is entrywise
fully faithful and surjective on objects. It is therefore a trivial
fibration in $\Cyc(\Grpd)$. The lifting property of $g$ and the
adjunction $\ob\dashv\operatorname{Codisc}$ give a lift of
$\ob(g)$ against $q$.

Conversely, suppose that $\ob(g)$ has the stated lifting property.
Let $p\colon\mathcal A\to\mathcal B$ be a trivial fibration in
$\Cyc(\Grpd)$, and consider a commutative square
\[
\begin{tikzcd}
\mathcal P \arrow[r,"f"] \arrow[d,"g"'] &
\mathcal A \arrow[d,"p"] \\
\mathcal Q \arrow[r,"h"'] &
\mathcal B.
\end{tikzcd}
\]
The map $\ob(p)$ is entrywise surjective. Hence there is a morphism
of cyclic operads in sets
$\ell_{\mathrm{ob}}\colon\ob(\mathcal Q)\to\ob(\mathcal A)$ such
that
\[
\ell_{\mathrm{ob}}\cdot\ob(g)=\ob(f)
\qquad\text{and}\qquad
\ob(p)\cdot\ell_{\mathrm{ob}}=\ob(h).
\]

For every $n$, the functor
$p_n\colon\mathcal A(n)\to\mathcal B(n)$ is fully faithful. Given a
morphism $\alpha\colon x\to y$ in $\mathcal Q(n)$, there is therefore
a unique morphism
\[
\ell_n(\alpha)\colon
\ell_{\mathrm{ob}}(x)\longrightarrow\ell_{\mathrm{ob}}(y)
\]
whose image under $p_n$ is $h_n(\alpha)$. Uniqueness shows that these
morphisms preserve identities and composition.

They also preserve operadic composition and the extended
symmetric-group actions. After applying $p$, the required equations
are the corresponding equations for $h$. They then hold in
$\mathcal A$ because each $p_n$ is faithful. The functors $\ell_n$
therefore assemble into a morphism of cyclic operads
$\ell\colon\mathcal Q\to\mathcal A$.

By construction, $p\cdot\ell=h$. The morphisms $\ell\cdot g$ and
$f$ agree on objects, and their images under $p$ agree. The
faithfulness of $p$ gives $\ell\cdot g=f$. Thus $g$ has the left
lifting property with respect to every trivial fibration and is a
cofibration.
\end{proof}

\begin{cor}
\label{cor:free-cyc-cof}
Let $\mathcal P$ be a cyclic operad in groupoids. Suppose that
$\ob(\mathcal P)$ is freely generated as a cyclic operad by a cyclic
collection $X$ and that the $\Sigma_n^+$-action on $X(n)$ is free for
every $n$. Then $\mathcal P$ is cofibrant in $\Cyc(\Grpd)$.
\end{cor}

\begin{proof}
Let $q\colon\mathcal A\to\mathcal B$ be an entrywise surjective
morphism in $\Cyc(\Set)$, and let
$h\colon\ob(\mathcal P)\to\mathcal B$ be a morphism of cyclic
operads.

For each $\Sigma_n^+$-orbit in $X(n)$, choose one generator and a
preimage in $\mathcal A(n)$ of its image under $h$. Since the
$\Sigma_n^+$-action on $X(n)$ is free, these choices extend uniquely
to an equivariant lift
\[
X\longrightarrow\mathcal A
\]
of the restriction of $h$ to $X$. Since $\ob(\mathcal P)$ is freely
generated by $X$, this lift extends uniquely to a morphism
$\ob(\mathcal P)\to\mathcal A$ of cyclic operads lifting $h$.
Proposition~\ref{prop:cofibrations-cyc} now implies that $\mathcal P$
is cofibrant.
\end{proof}

\subsubsection{Free operads generated by a tree}
\label{subsec: free operad gen by tree}

The dendroidal category is built from the coloured operads freely
generated by rooted trees. We retain a planar structure on each tree
in order to specify the ordering of the inputs at its vertices. Only
the set-valued construction is needed here. The dendroidal objects
introduced below will take values in Cartesian categories.

\begin{definition}
\label{def: decoration of tree}
Let $\mathfrak C$ be a nonempty set. A rooted tree $T$ is
\emph{$\mathfrak C$-coloured} if it is equipped with a map
$\kappa\colon E(T)\to\mathfrak C$, where $E(T)$ denotes its set of
edges.

Let $\mathcal P$ be a $\mathfrak C$-coloured symmetric sequence in
$\Set$. A \emph{$\mathcal P$-decoration} of a
$\mathfrak C$-coloured planar rooted tree $T$ assigns to every vertex
$v\in V(T)$ an element
\[
p(v)\in
\mathcal P\bigl(
\kappa(e_1),\ldots,\kappa(e_k);\kappa(e_0)
\bigr),
\]
where $e_1,\ldots,e_k$ are the incoming edges at $v$, ordered by the
planar structure, and $e_0$ is its outgoing edge.
\end{definition}

The free operad $F(\mathcal P)$ consists of isomorphism classes of
$\mathcal P$-decorated, $\mathfrak C$-coloured planar rooted trees.
More precisely,
$F(\mathcal P)(c_1,\ldots,c_n;c_0)$ is represented by such trees
whose root has colour $c_0$ and whose leaves are labelled by
$\{1,\ldots,n\}$, with the leaf labelled by $i$ having colour $c_i$.
One quotients by isomorphisms of decorated trees and by the
symmetric-group actions at the vertices. Operadic composition is
given by grafting, and the symmetric-group action is given by
relabelling the leaves. See
\cite[Section~3 and Theorem~3.2]{bm_resolutions} and
\cite[Definition~2.3]{cacti}.

Now let $T$ be a planar rooted tree and take $\mathfrak C=E(T)$.
Define an $E(T)$-coloured symmetric sequence $X(T)$ as follows. If
$v$ has ordered incoming edges $e_1,\ldots,e_n$ and outgoing edge
$e_0$, then $X(T)(e_1,\ldots,e_n;e_0)$ contains a generator denoted
by $v$. The components obtained by permuting the incoming colours
contain its corresponding symmetric-group translates. All other
components are empty.

\begin{definition}
The \emph{free operad generated by $T$} is the $E(T)$-coloured
operad
\[
\Omega(T)\coloneqq F(X(T)).
\]
\end{definition}

The non-unital operations of $\Omega(T)$ are represented by connected
subtrees of $T$, with the orientation induced by the root. Thus
$\Omega(T)(e_1,\ldots,e_n;e_0)$ contains a unique operation if there
is a subtree with outgoing boundary edge $e_0$ and ordered incoming
boundary edges $e_1,\ldots,e_n$, and is empty otherwise. The unit in
$\Omega(T)(e;e)$ is represented by the exceptional one-edge tree
$\eta_e$. Operadic composition grafts compatible subtrees along
matching boundary edges, while the symmetric-group actions reorder
the incoming boundary edges.
\subsubsection{The free cyclic operad generated by a tree}
\label{sec:free cyclic operad}

Coloured cyclic operads require a slight extension of the usual
notion of a coloured operad. Following
\cite[Section~2.3]{DCH}, we equip the set of colours $\mathfrak C$
with an involution $c\mapsto c^*$, which is allowed to be the
identity. See also
\cite[Definition~3.2]{doherty2025models}.

\begin{definition}
\label{def: coloured cyclic operad}
Let $\mathfrak C$ be a set equipped with an involution. A
\emph{$\mathfrak C$-coloured cyclic operad} $\mathcal P$ in $\bE$
consists of objects $\mathcal P(c_0,\ldots,c_n)$, indexed by
nonempty finite lists of colours, together with structure maps
\[
\sigma^*\colon
\mathcal P(c_0,\ldots,c_n)
\longrightarrow
\mathcal P(c_{\sigma(0)},\ldots,c_{\sigma(n)})
\]
for $\sigma\in\Sigma_n^+$, units
$\id_c\colon\mathbb 1\to\mathcal P(c^*,c)$, and partial
compositions
\[
\circ_i^j\colon
\mathcal P(c_0,\ldots,c_n)\otimes
\mathcal P(d_0,\ldots,d_m)
\longrightarrow
\mathcal P(\underline c\circ_i^j\underline d),
\]
defined whenever $c_i=d_j^*$. Here
\[
\underline c\circ_i^j\underline d
=
(c_0,\ldots,c_{i-1},
d_{j+1},\ldots,d_m,d_0,\ldots,d_{j-1},
c_{i+1},\ldots,c_n).
\]
These data satisfy the commutativity, associativity, equivariance,
and unit axioms of \cite[Definition~2.3]{DCH}.
\end{definition}

We follow the convention of excluding operations indexed by the empty
profile; see \cite[Remark~4.3]{DCH}.

A planar tree $T$ is \emph{$\mathfrak C$-coloured} if it is equipped
with a map $\kappa\colon H(T)\to\mathfrak C$ satisfying
$\kappa(i(h))=\kappa(h)^*$. Let $\mathcal P$ be a $\mathfrak C$-coloured cyclic symmetric
sequence in $\Set$. Thus, for every $\sigma\in\Sigma_n^+$, it has
structure maps
\[
\sigma^*\colon
\mathcal P(c_0,\ldots,c_n)
\longrightarrow
\mathcal P(c_{\sigma(0)},\ldots,c_{\sigma(n)}).
\]
A \emph{$\mathcal P$-decoration} of a
$\mathfrak C$-coloured planar tree assigns to every vertex $v$ an
element
\[
p(v)\in
\mathcal P(\kappa(h_0),\ldots,\kappa(h_n)),
\]
where $h_0,\ldots,h_n$ are the half-edges incident to $v$, listed in
an order compatible with the cyclic order determined by the planar
structure.

\begin{definition}
\label{def: free cyclic operad}
The \emph{free cyclic operad} $F(\mathcal P)$ is the
$\mathfrak C$-coloured cyclic operad represented by
$\mathfrak C$-coloured, $\mathcal P$-decorated planar trees equipped
with boundary labellings.

More precisely, $F(\mathcal P)(c_0,\ldots,c_k)$ is the quotient
\[
\left(
\coprod_{(T,\kappa,\tau)}
\prod_{v\in V(T)}
\mathcal P(\operatorname{nb}(v))
\right)\bigg/\sim,
\]
where
$\tau\colon\{0,\ldots,k\}\xrightarrow{\cong}\partial T$
is a boundary labelling satisfying
$\kappa(\tau(i))=c_i$. The equivalence relation is generated by
isomorphisms of coloured decorated planar trees and by changes in
the chosen ordering of the half-edges at a vertex, compensated by
the corresponding extended symmetric-group action on its
decoration.

The right $\Sigma_k^+$-action is induced by relabelling the boundary
edges. Partial composition joins the boundary edge labelled $i$ in
one tree to the boundary edge labelled $j$ in another. It is defined
when $c_i=d_j^*$, and the remaining boundary edges are labelled
according to $\underline c\circ_i^j\underline d$. Units are
represented by the exceptional trees with no vertices.
\end{definition}

Now let $T$ be a planar unrooted tree and take
$\mathfrak C=H(T)$, with colour involution given by the half-edge
involution. Define an $H(T)$-coloured cyclic symmetric sequence
$X(T)$ as follows. If $v\in V(T)$ has incident half-edges
$\operatorname{nb}(v)=(h_0,\ldots,h_n)$, set
\[
X(T)(h_{\sigma(0)},\ldots,h_{\sigma(n)})
=
\{v_\sigma\}
\]
for every $\sigma\in\Sigma_n^+$, and let every other component be
empty. Its extended symmetric-group structure is given by
$\rho^*(v_\sigma)=v_{\sigma\rho}$.

\begin{definition}
\label{def: free cyclic operad generated by a tree}
The \emph{free cyclic operad generated by $T$} is the
$H(T)$-coloured cyclic operad
\[
\Omega^{\mathrm{cyc}}(T)\coloneqq F(X(T)).
\]
\end{definition}

The non-unital operations of $\Omega^{\mathrm{cyc}}(T)$ are represented by
subtrees of $T$. More precisely,
$\Omega^{\mathrm{cyc}}(T)(h_0,\ldots,h_n)$ contains a unique
operation when the distinct half-edges $h_0,\ldots,h_n$ form the
boundary of a subtree of $T$, and is empty otherwise. Partial
composition joins compatible subtrees. This is the construction
described in
\cite[Definitions~4.6 and~4.8]{doherty2025models}.

\subsection{\texorpdfstring{$\infty$}{infinity}-operads}
\label{subsec: infinity operad}

We model one-coloured $\infty$-operads using a planar presentation of
the dendroidal category. Its objects are planar rooted trees, and its
morphisms are defined by
\[
\Hom_{\Omega}(S,T)
=
\Hom_{\Op(\Set)}\bigl(\Omega(S),\Omega(T)\bigr).
\]
The planar structures specify orderings of the inputs at the vertices,
but morphisms are maps of the associated symmetric coloured operads
and need not preserve these orderings. This is the presentation used
in \cite[p.~235]{hackney_config}, where it is identified with the
category denoted $\Omega'$ in
\cite[Example~2.8]{bm_reedy}. It is equivalent to the usual dendroidal
category, and we denote it by $\Omega$ throughout.

Morphisms in $\Omega$ are generated by isomorphisms, face maps, and
degeneracies. If $e$ is an internal edge of $T$, contracting $e$
gives an inner face map $\partial_e\colon T/e\to T$. An outer vertex
of a tree with at least two vertices is a vertex incident to exactly
one inner edge. Removing such a vertex and its adjacent outer edges
gives an outer face map $\partial_v\colon T/v\to T$. If $v$ is unary,
contracting $v$ and identifying its two incident edges gives a
degeneracy $s_v\colon T\to T/v$. See
\cite{MW1,moerdijk_lectures,HM}.

Let $\mathcal E$ be an $\infty$-category with finite products. A
\emph{dendroidal object} in $\mathcal E$ is a functor
\[
X\colon\Omega^{\mathrm{op}}\longrightarrow\mathcal E.
\]
Thus the $\infty$-category of dendroidal objects in $\mathcal E$ is
$\operatorname{Fun}(\Omega^{\mathrm{op}},\mathcal E)$.

For each vertex $v\in V(T)$, the corolla at $v$ determines a morphism
$C_{|\mathrm{in}(v)|}\to T$ in $\Omega$. Applying a dendroidal object
$X$ and taking the product over all vertices produces the
\emph{Segal map}
\[
X_T\longrightarrow
\prod_{v\in V(T)}X_{C_{|\mathrm{in}(v)|}}.
\]

\begin{definition}
\label{def: infty operad}
Let $\mathcal E$ be an $\infty$-category with finite products. An
\emph{$\infty$-operad in $\mathcal E$} is a dendroidal object
$X\colon\Omega^{\mathrm{op}}\to\mathcal E$ such that $X_\eta\simeq *$ and,
for every tree $T$, the Segal map
\[
X_T\longrightarrow
\prod_{v\in V(T)}X_{C_{|\mathrm{in}(v)|}}
\]
is an equivalence in $\mathcal E$.

We write $\Op_\infty(\mathcal E)$ for the full sub-$\infty$-category
of $\Fun(\Omega^{\mathrm{op}},\mathcal E)$ spanned by these objects.
\end{definition}

The condition $X_\eta\simeq *$ restricts us to one-coloured
$\infty$-operads, since $X_\eta$ represents the object of colours.
The object $X_{C_n}$ represents the $n$-ary operations, while the
Segal condition says that the value on a tree is determined, up to
equivalence, by the operations associated to its vertices.

Suppose now that $\bE$ is a Cartesian model category presenting the
Cartesian $\infty$-category $\mathcal E$. An entrywise fibrant
dendroidal object $X\colon\Omega^{\mathrm{op}}\to\bE$ represents an object of
$\Op_\infty(\mathcal E)$ provided that $X_\eta\to *$ is a weak
equivalence and all its Segal maps are weak equivalences in $\bE$.
Entrywise fibrancy ensures
that the ordinary finite products appearing in the Segal maps
represent the corresponding products in $\mathcal E$.

Equivalently, let $\Op_\infty(\bE)$ denote the relative category whose objects are the entrywise fibrant dendroidal objects satisfying these conditions and whose weak equivalences are the entrywise weak equivalences. Its localization presents $\Op_\infty(\mathcal E)$. We use the model-categorical notation when
working with explicit representatives and the $\infty$-categorical notation when discussing derived mapping spaces.

Every one-coloured operad $\mathcal O$ in $\bE$ determines a
dendroidal object through the fully faithful \emph{dendroidal nerve}
\[
\nerve\colon
\Op(\bE)\longrightarrow
\Fun(\Omega^{\mathrm{op}},\bE).
\]
Its value at a tree is
\[
(\nerve\mathcal O)_T
=
\prod_{v\in V(T)}
\mathcal O\bigl(|\mathrm{in}(v)|\bigr),
\]
with functoriality induced by the operadic structure maps of
$\mathcal O$. The nerve satisfies the Segal condition strictly.
Consequently, if $\mathcal O$ is entrywise fibrant, then
$\nerve\mathcal O$ represents an object of
$\Op_\infty(\mathcal E)$. See \cite{MW1} for the dendroidal nerve and
its full faithfulness.

\begin{remark}
When $\bE$ is $\sSet$ or $\Grpd$, the relative category above may be
presented model-categorically by localizing the projective model
structure on $\Fun(\Omega^{\mathrm{op}},\bE)$ at the Segal maps and at the map
$\varnothing\to\Omega(-,\eta)$. Its fibrant objects are precisely the
entrywise fibrant dendroidal objects satisfying the two conditions in
Definition~\ref{def: infty operad}; see
\cite{cisinski2013dendroidal} for the corresponding construction for
dendroidal spaces.

The model structure on profinite spaces is fibrantly generated, so
this localization argument is not immediately available in that
setting. We will instead define profinite $\infty$-operads directly as
Segal objects in the underlying Cartesian $\infty$-category of
profinite spaces.
\end{remark}

The dendroidal nerve gives strict Segal objects, but naturally
occurring homotopy-coherent operadic structures need not be presented
in this form. Configuration spaces provide a useful non-unital
example.

\begin{definition}
\label{def: semi-dendroidal space}
Let $\Omega_{\mathrm{inj}}\subseteq\Omega$ be the wide subcategory
generated by face maps and isomorphisms. A
\emph{semi-dendroidal space} is a functor
$X\colon\Omega_{\mathrm{inj}}^{\mathrm{op}}\to\sSet$. It satisfies the
\emph{Segal condition} if, for every tree $T$, the canonical map
\[
X_T\longrightarrow
\prod_{v\in V(T)}X_{C_{|\mathrm{in}(v)|}}
\]
is a weak equivalence.
\end{definition}

A semi-dendroidal space satisfying the Segal condition may be regarded
as a non-unital $\infty$-operad: restricting to
$\Omega_{\mathrm{inj}}$ retains the face maps but omits the
degeneracies that encode operadic units.

\begin{example}
\label{ex: infinity operad of configuration spaces}
The ordered configuration spaces $\Conf_n(\mathbb C)$ do not form a
strict operad, although forgetting the radii of the disks gives weak
equivalences
$\mathsf D_2(n)\simeq\Conf_n(\mathbb C)$.

Hackney constructs a semi-dendroidal space $X$ satisfying the Segal
condition, with
$X_{C_n}\simeq\Conf_n(\mathbb C)$. Hence, for every tree $T$, the
Segal map
\[
X_T\longrightarrow
\prod_{v\in V(T)}X_{C_{|\mathrm{in}(v)|}}
\]
is a weak equivalence. The construction begins with a two-coloured
operad of configurations of points and disks. Its shift maps replace
disks by their centres.

The resulting semi-dendroidal space does not extend to a dendroidal
space. In particular, no compatible collection of degeneracy maps
can be added. See
\cite[Corollary~5 and Proposition~6]{hackney_config}.
\end{example}

\subsection{Cyclic \texorpdfstring{$\infty$}{infinity}-operads}
\label{subsec: cyclic infinity operad}

Let $T$ be a planar unrooted tree. As described in
Section~\ref{sec:free cyclic operad}, it generates a free cyclic
operad $\Omega^{\mathrm{cyc}}(T)$.

We use the planar presentation of the \emph{cyclic dendroidal
category}, denoted $\Omega^{\mathrm{cyc}}$. Its objects are planar
unrooted trees, and its morphisms are given by
\[
\Hom_{\Omega^{\mathrm{cyc}}}(S,T)
=
\Hom_{\Cyc(\Set)}
\bigl(
\Omega^{\mathrm{cyc}}(S),
\Omega^{\mathrm{cyc}}(T)
\bigr).
\]
Although the objects carry chosen planar structures, morphisms need
not preserve them. This category is equivalent to the category
$U_{\mathrm{cyc}}$ introduced in \cite{hry1} and to the category
$\Upsilon$ of \cite[Definition~4.8]{doherty2025models}. As in the
rooted case, its morphisms are generated by isomorphisms, inner and
outer face maps, and degeneracies.

Let $\mathcal E$ be an $\infty$-category with finite products. A
\emph{cyclic dendroidal object} in $\mathcal E$ is a functor
\[
X\colon
(\Omega^{\mathrm{cyc}})^{\mathrm{op}}
\longrightarrow
\mathcal E.
\]
For every vertex $v$ of an unrooted tree $T$, the corolla at $v$
determines a morphism
$C_{|\operatorname{nb}(v)|}\to T$ in
$\Omega^{\mathrm{cyc}}$. Thus a vertex with $n+1$ incident
half-edges determines the unrooted corolla $C_{n+1}$ and corresponds
to an $n$-ary cyclic operation. Applying $X$ and taking the product
over the vertices of $T$ gives the \emph{Segal map}
\[
X_T\longrightarrow
\prod_{v\in V(T)}
X_{C_{|\operatorname{nb}(v)|}}.
\]

\begin{definition}
\label{def: cyclic infinity operad}
Let $\mathcal E$ be an $\infty$-category with finite products. A
\emph{cyclic $\infty$-operad in $\mathcal E$} is a cyclic dendroidal
object $X$ such that $X_\eta\simeq *$ and, for every unrooted tree
$T$, the Segal map
\[
X_T\longrightarrow
\prod_{v\in V(T)}
X_{C_{|\operatorname{nb}(v)|}}
\]
is an equivalence in $\mathcal E$.

We write $\Cyc_\infty(\mathcal E)$ for the full
sub-$\infty$-category of
$\operatorname{Fun}((\Omega^{\mathrm{cyc}})^{\mathrm{op}},\mathcal E)$
spanned by these objects.
\end{definition}

The condition $X_\eta\simeq *$ restricts us to one-coloured cyclic
$\infty$-operads. The Segal condition says that the value on an
unrooted tree is determined, up to equivalence, by the operations
associated to its vertices.

Suppose that $\bE$ is a Cartesian model category presenting the
Cartesian $\infty$-category $\mathcal E$. An entrywise fibrant cyclic
dendroidal object
$X\colon(\Omega^{\mathrm{cyc}})^{\mathrm{op}}\to\bE$
represents an object of $\Cyc_\infty(\mathcal E)$ provided that
$X_\eta\to *$ is a weak equivalence and all its Segal maps are weak
equivalences in $\bE$. Entrywise fibrancy ensures that the ordinary finite products
appearing in the Segal maps represent the corresponding products in
$\mathcal E$.

Equivalently, let $\Cyc_\infty(\bE)$ denote the relative category
whose objects are the entrywise fibrant cyclic dendroidal objects
satisfying these conditions and whose weak equivalences are the
entrywise weak equivalences. Its localization presents
$\Cyc_\infty(\mathcal E)$.

Every one-coloured cyclic operad $\mathcal P$ in $\bE$ determines a
cyclic dendroidal object through the fully faithful
\emph{cyclic dendroidal nerve}
\[
\nerve_{\mathrm{cyc}}\colon
\Cyc(\bE)
\longrightarrow
\operatorname{Fun}
\bigl((\Omega^{\mathrm{cyc}})^{\mathrm{op}},\bE\bigr).
\]
Its value at an unrooted tree $T$ is
\[
(\nerve_{\mathrm{cyc}}\mathcal P)_T
=
\prod_{v\in V(T)}
\mathcal P\bigl(|\operatorname{nb}(v)|-1\bigr),
\]
with functoriality induced by the cyclic-operadic structure maps of
$\mathcal P$. When $\bE=\Set$, this agrees with
\[
(\nerve_{\mathrm{cyc}}\mathcal P)_T
\cong
\Hom_{\Cyc(\Set)}
\bigl(\Omega^{\mathrm{cyc}}(T),\mathcal P\bigr).
\]
The cyclic dendroidal nerve satisfies the Segal condition strictly.
Consequently, if $\mathcal P$ is entrywise fibrant, then
$\nerve_{\mathrm{cyc}}\mathcal P$ represents an object of
$\Cyc_\infty(\mathcal E)$.

\begin{remark}
Throughout this paper, we use both the ordinary dendroidal nerve
$\nerve_{\mathrm{op}}$ and the cyclic dendroidal nerve
$\nerve_{\mathrm{cyc}}$. When applied to a cyclic operad, the
unsubscripted notation $\nerve$ always denotes
$\nerve_{\mathrm{cyc}}$. We write $\nerve_{\mathrm{op}}$ explicitly
for the ordinary dendroidal nerve applied to an operad.
\end{remark}

\subsection{Properties of the (cyclic) dendroidal nerve}
\label{subsec: properties of the (cyclic) dendroidal nerve}

Under the hypotheses of
Proposition~\ref{model category on alg over operad}, the categories
$\Op(\bE)$ and $\Cyc(\bE)$ carry transferred model structures in
which weak equivalences and fibrations are defined entrywise. The
relative categories $\Op_\infty(\bE)$ and $\Cyc_\infty(\bE)$ were
defined above using entrywise fibrant dendroidal and cyclic
dendroidal objects satisfying the one-coloured Segal conditions.

The following homotopical full-faithfulness result is the cyclic
analogue of
\cite[Proposition~4.3]{Boavida-Horel-Robertson}. The proof given
there for operads uses the strict Segal description of the
dendroidal nerve and a resolution by free operads. Replacing free
operads by free cyclic operads gives the cyclic statement. The same
argument is carried out for modular operads in
\cite[Theorem~B.2]{BR_modular} and restricts to cyclic operads
through the genus-zero identification of
\cite[Proposition~A.11]{BR_modular}.

\begin{theorem}
\label{thm: dendroidal nerve is homotopically fully faithful}
Let $\bE$ be either $\sSet$, with the Kan--Quillen model structure,
or $\Grpd$, with its standard model structure. For all entrywise
fibrant $\mathcal O,\mathcal P\in\Op(\bE)$, the dendroidal nerve
induces a weak equivalence
\[
\mathbb R\Map_{\Op(\bE)}(\mathcal O,\mathcal P)
\longrightarrow
\mathbb R\Map_{\Op_\infty(\bE)}
\bigl(
\nerve_{\mathrm{op}}\mathcal O,
\nerve_{\mathrm{op}}\mathcal P
\bigr).
\]
Similarly, for all entrywise fibrant
$\mathcal O,\mathcal P\in\Cyc(\bE)$, the cyclic dendroidal nerve
induces a weak equivalence
\[
\mathbb R\Map_{\Cyc(\bE)}(\mathcal O,\mathcal P)
\longrightarrow
\mathbb R\Map_{\Cyc_\infty(\bE)}
\bigl(
\nerve_{\mathrm{cyc}}\mathcal O,
\nerve_{\mathrm{cyc}}\mathcal P
\bigr).
\]
\end{theorem}

For operads and cyclic operads in groupoids, entrywise fibrancy is
automatic. The same applies to the classifying spaces of such
groupoids, since they are Kan complexes.
There is a \emph{root-elision} functor
$f\colon\Omega\to\Omega^{\mathrm{cyc}}$ that sends a rooted tree
$(T,r)$ to its underlying unrooted tree $T$. On morphisms, it is
induced by applying the cyclic-envelope functor to maps between the
free operads generated by rooted trees; see
\cite[Definition~4.11]{doherty2025models}.

Restriction along $f$ defines a functor from cyclic dendroidal
objects to dendroidal objects. Since root elision preserves the
exceptional tree, corollas, and the vertices of a tree, it preserves
the one-coloured and Segal conditions. It therefore restricts to a
functor
$f^*\colon\Cyc_\infty(\mathcal E)\to\Op_\infty(\mathcal E)$.
Whenever $\mathcal E$ admits the required Kan extensions, restriction
fits into an adjoint string
$f_!\dashv f^*\dashv f_*$, where $f_!$ and $f_*$ are left and right
Kan extension.

This construction is compatible with the forgetful functor
$u^*\colon\Cyc(\bE)\to\Op(\bE)$. Recall that $u^*$ has the cyclic
envelope $u_!$ as its left adjoint and the exceptional functor $u_*$
as its right adjoint; see \cite[Section~3]{DCH}. Thus
$u_!\dashv u^*\dashv u_*$.

For every rooted tree $(T,r)$, root elision gives an isomorphism
$u_!\Omega(T,r)\cong\Omega^{\mathrm{cyc}}(T)$; see
\cite[Definition~4.11 and the following discussion]
{doherty2025models}. Consequently, the adjunction
$u_!\dashv u^*$ gives, for every cyclic operad $\mathcal P$, an
isomorphism
\[
\Hom_{\Cyc(\bE)}
\bigl(\Omega^{\mathrm{cyc}}(T),\mathcal P\bigr)
\cong
\Hom_{\Op(\bE)}
\bigl(\Omega(T,r),u^*\mathcal P\bigr).
\]
Equivalently, the cyclic and ordinary dendroidal nerves are related
by
$f^*\nerve_{\mathrm{cyc}}(\mathcal P)
\cong\nerve_{\mathrm{op}}(u^*\mathcal P)$.
Thus the diagram
\[
\begin{tikzcd}
\Cyc(\bE)
    \arrow[r,"\nerve_{\mathrm{cyc}}"]
    \arrow[d,"u^*"']
&
\operatorname{Fun}
\bigl((\Omega^{\mathrm{cyc}})^{\mathrm{op}},\bE\bigr)
    \arrow[d,"f^*"]
\\
\Op(\bE)
    \arrow[r,"\nerve_{\mathrm{op}}"']
&
\operatorname{Fun}(\Omega^{\mathrm{op}},\bE)
\end{tikzcd}
\]
commutes up to natural isomorphism.

\subsection{A cyclic \texorpdfstring{$\infty$}{infinity}-operad of
framed configurations}
\label{subsec:framed_conf_operad}

In this subsection, we construct a topological semi-dendroidal space
$X^{\mathrm{fr}}\colon\Omega_{\mathrm{inj}}^{\mathrm{op}}\to\Top$
whose value at the corolla $C_n$ is weakly equivalent to the framed
configuration space
$\Conf_n^{\mathrm{fr}}(\mathbb C)
=\Conf_n(\mathbb C)\times(S^1)^n$.
Although the construction itself records only rooted operations, we
show that $X^{\mathrm{fr}}$ is weakly equivalent to the underlying
non-unital rooted part of a cyclic $\infty$-operad.

Here and below, we use the topological analogue of a semi-dendroidal
space. Applying the singular-complex functor gives a semi-dendroidal
space in the sense of Definition~\ref{def: semi-dendroidal space}.

To identify the cyclic enhancement geometrically, we use two cyclic
models for the framed little disks operad: the operad $\mathsf S_0$
of parametrised disk configurations in the sphere, recalled in
Example~\ref{ex: Budney}, and the framed Fulton--MacPherson operad
$\mathsf{FFM}_2$. We first compare these models as cyclic
$\infty$-operads and then identify their underlying non-unital rooted
models with $X^{\mathrm{fr}}$.

Our construction is a framed version of Hackney's non-unital
$\infty$-operad of configuration spaces
(Example~\ref{ex: infinity operad of configuration spaces} and
\cite[Example~1]{hackney_config}). It is built from a two-coloured
operad distinguishing disks from points, with shift maps that replace
selected disks by their centres. We equip both disks and points with
framings and require the shift maps to preserve them.

Let $D\subset\mathbb C$ denote the open unit disk.

\begin{definition}
\label{def:operad-with-shifts}
We define a two-coloured topological operad $\mathsf{PD}$ with colour
set $\{1,2\}$. For an ordered list of colours
$\ell=(l_1,\ldots,l_n)$, set
\[
\mathsf{PD}(\ell;2)
=
\begin{cases}
* & \text{if $\ell=(2)$},\\
\varnothing & \text{otherwise}.
\end{cases}
\]
An element $x\in\mathsf{PD}(\ell;1)$ is a configuration
$x=(x_1,\ldots,x_n)$ of framed disks and framed points in $D$. If
$l_j=1$, then $x_j=(r_j,c_j,u_j)$, where $r_j>0$, $c_j\in D$, and
$u_j\in S^1$, represents the affine embedding
$z\mapsto r_ju_jz+c_j$ of $D$ into itself. If $l_j=2$, then
$x_j=(c_j,u_j)\in D\times S^1$ is a framed point.

We require the closure of the image of every disk to lie in $D$. The
closures of the images of distinct disks must be disjoint, the
framed points must have distinct centres, and no framed point may lie
in the closure of the image of a disk. We give
$\mathsf{PD}(\ell;1)$ the resulting subspace topology.

The symmetric-group actions relabel the entries of a configuration.
Let $x\in\mathsf{PD}(\ell;1)$ and
$y\in\mathsf{PD}(\lambda;1)$, with $l_i=1$, and write
$x_i=(r_i,c_i,u_i)$. The composite $x\circ_i y$ replaces the $i$th
disk of $x$ by the image of $y$ under the corresponding affine
embedding. A framed disk or framed point of $y$ is transformed
according to
\[
\begin{aligned}
(r_j',c_j',u_j')
&\longmapsto
\bigl(r_ir_j',\,r_iu_ic_j'+c_i,\,u_iu_j'\bigr),\\
(c_j',u_j')
&\longmapsto
\bigl(r_iu_ic_j'+c_i,\,u_iu_j'\bigr).
\end{aligned}
\]
Substitution into an input of colour $2$ is determined by the unique
operation in $\mathsf{PD}((2);2)$. The units of colours $1$ and $2$
are the identity embedding of $D$ and the unique element of
$\mathsf{PD}((2);2)$, respectively.

For colour strings $\ell\leq\ell'$ of equal length, ordered entrywise
by $1<2$, define a shift map
$s_{\ell,\ell'}\colon\mathsf{PD}(\ell;1)\to
\mathsf{PD}(\ell';1)$ as follows. Components whose colours agree in
$\ell$ and $\ell'$ are unchanged, while a framed disk
$(r_i,c_i,u_i)$ is replaced by its framed centre $(c_i,u_i)$ whenever
$l_i=1$ and $l_i'=2$.
\end{definition}

Each shift map is a weak homotopy equivalence. The argument is the
same as in \cite[Example~1]{hackney_config}: a framed point may be
replaced by a sufficiently small framed disk, and the space of
permissible radii is contractible. The framing is left unchanged
throughout.

The shift maps satisfy $s_{\ell,\ell}=\id$ and
$s_{\ell',\ell''}\cdot s_{\ell,\ell'}=s_{\ell,\ell''}$ whenever
$\ell\leq\ell'\leq\ell''$. They are also compatible with operadic
substitution. If $\ell\leq\ell'$, $\lambda\leq\lambda'$, and
$l_i=l_i'=1$, then
\[
s_{\ell,\ell'}(x)\circ_i s_{\lambda,\lambda'}(y)
=
s_{\ell\circ_i\lambda,\,
  \ell'\circ_i\lambda'}(x\circ_i y).
\]
The condition $l_i=l_i'=1$ ensures that both composites are defined.
Finally, if $\sigma\in\Sigma_n$ and
$\ell\sigma=(l_{\sigma(1)},\ldots,l_{\sigma(n)})$, then
\[
s_{\ell\sigma,\ell'\sigma}(\sigma^*x)
=
\sigma^*s_{\ell,\ell'}(x).
\]
It follows that $\mathsf{PD}$, together with the maps
$s_{\ell,\ell'}$, is an operad with shifts in the sense of
\cite{hackney_config}.

The dendroidal nerve of $\mathsf{PD}$ records colourings of the edges
of a rooted tree together with compatible decorations of its
vertices. For a planar rooted tree $T$, there is an isomorphism
\[
\nerve_{\mathrm{op}}(\mathsf{PD})_T
\cong
\coprod_{\kappa\colon E(T)\to\{1,2\}}
\prod_{v\in V(T)}
\mathsf{PD}\bigl(
\kappa(\mathrm{in}(v));\kappa(\mathrm{out}(v))
\bigr).
\]
Here $\kappa(\mathrm{in}(v))$ is the ordered list of colours of the
incoming edges at $v$, while $\kappa(\mathrm{out}(v))$ is the colour
of its outgoing edge.

For every nontrivial rooted tree $T$, let
$\kappa_T\colon E(T)\to\{1,2\}$ assign colour $2$ to the leaves and
colour $1$ to the root and internal edges.

\begin{definition}
\label{def: infinity operad of framed config}
Set $X^{\mathrm{fr}}_\eta=*$. For every nontrivial rooted tree $T$,
let $X^{\mathrm{fr}}_T$ be the summand of
$\nerve_{\mathrm{op}}(\mathsf{PD})_T$ indexed by $\kappa_T$. Thus
\[
X^{\mathrm{fr}}_T
=
\prod_{v\in V(T)}
\mathsf{PD}(\ell_v;1),
\]
where $\ell_v=\kappa_T(\mathrm{in}(v))$ is the ordered list of colours
of the incoming edges at $v$.

Let $\alpha\colon S\to T$ be a morphism in
$\Omega_{\mathrm{inj}}$, with $S\neq\eta$. To define
$\widehat\alpha\colon X^{\mathrm{fr}}_T\to X^{\mathrm{fr}}_S$, first
apply the restriction map of
$\nerve_{\mathrm{op}}(\mathsf{PD})$ and then apply the shift maps at
the vertices where an edge becomes a leaf. If $S=\eta$, let
$\widehat\alpha$ be the unique map to $X^{\mathrm{fr}}_\eta=*$. The
compatibility of the shift maps with operadic substitution and the
symmetric-group actions implies, as in
\cite[Theorem~3]{hackney_config}, that these maps are functorial.
They therefore define a semi-dendroidal space
$X^{\mathrm{fr}}\colon
\Omega_{\mathrm{inj}}^{\mathrm{op}}\to\Top$.
\end{definition}

This is the framed analogue of Hackney's construction.

\begin{prop}
\label{prop: framed configurations Segal}
The semi-dendroidal space $X^{\mathrm{fr}}$ satisfies the Segal
condition, and
$X^{\mathrm{fr}}_{C_n}\simeq
\Conf_n^{\mathrm{fr}}(\mathbb C)$ for every $n\geq0$.
\end{prop}

\begin{proof}
Let $T$ be a nontrivial rooted tree. Its Segal map
\[
X^{\mathrm{fr}}_T
\longrightarrow
\prod_{v\in V(T)}
X^{\mathrm{fr}}_{C_{|\mathrm{in}(v)|}}
\]
is induced by the corolla inclusions
$C_{|\mathrm{in}(v)|}\to T$. Under $\kappa_T$, the leaves have
colour $2$, while the root and internal edges have colour $1$. For
each vertex $v$, write
$\ell_v=\kappa_T(\mathrm{in}(v))$, and let
$M_{|\mathrm{in}(v)|}$ denote the colour string consisting entirely
of $2$s. We then have
$X^{\mathrm{fr}}_T=\prod_v\mathsf{PD}(\ell_v;1)$ and
$X^{\mathrm{fr}}_{C_{|\mathrm{in}(v)|}}
=\mathsf{PD}(M_{|\mathrm{in}(v)|};1)$.

The corolla inclusion at $v$ replaces the internal input edges
incident to $v$ by leaves. Its structure map therefore applies
$s_{\ell_v,M_{|\mathrm{in}(v)|}}$ at $v$. Consequently, the Segal
map is the product of these shift maps over the vertices of $T$.
Each factor is a weak homotopy equivalence, and the product is finite,
so the Segal map is a weak homotopy equivalence.

For $n\geq1$, the corolla $C_n$ has $n$ input edges of colour $2$
and an output edge of colour $1$. Hence
$X^{\mathrm{fr}}_{C_n}=\mathsf{PD}(M_n;1)$, the space of
configurations of $n$ framed points $(c_j,u_j)\in D\times S^1$ with
distinct centres. It is therefore homeomorphic to
$\Conf_n(D)\times(S^1)^n=\Conf_n^{\mathrm{fr}}(D)$ and, since
$D\cong\mathbb C$, to $\Conf_n^{\mathrm{fr}}(\mathbb C)$. When
$n=0$, $\mathsf{PD}(M_0;1)$ consists of the empty configuration, so
$X^{\mathrm{fr}}_{C_0}\cong*
\cong\Conf_0^{\mathrm{fr}}(\mathbb C)$.
\end{proof}

The full one-coloured suboperad of $\mathsf{PD}$ on colour $1$ is
the framed little disks operad $\mathsf{FD}_2$. We use the version
with $\mathsf{FD}_2(0)=*$, represented by the empty configuration.
Indeed, if $m_n$
denotes the colour string consisting entirely of $1$s, then
$\mathsf{PD}(m_n;1)$ is the space of configurations of $n$ framed
disks in $D$ with pairwise disjoint closed images, and composition is
given by affine substitution. Consequently,
\[
\nerve_{\mathrm{op}}(\mathsf{FD}_2)_T
=
\prod_{v\in V(T)}
\mathsf{PD}(m_{|\mathrm{in}(v)|};1).
\]

For each vertex $v$, the inequality
$m_{|\mathrm{in}(v)|}\leq\ell_v$ holds entrywise. Define
$q_T\colon\nerve_{\mathrm{op}}(\mathsf{FD}_2)_T\to
X^{\mathrm{fr}}_T$ as the product, over the vertices of $T$, of the
shift maps $s_{m_{|\mathrm{in}(v)|},\ell_v}$. Compatibility with
operadic substitution and the symmetric-group actions makes these
maps natural with respect to face maps and isomorphisms. They
therefore assemble into a natural transformation
\[
q\colon
\left.\nerve_{\mathrm{op}}(\mathsf{FD}_2)
\right|_{\Omega_{\mathrm{inj}}}
\longrightarrow
X^{\mathrm{fr}}.
\]

\begin{cor}
\label{cor:nerve-FD2-eq-Xfr}
The natural transformation $q$ is an entrywise weak homotopy
equivalence of semi-dendroidal spaces.
\end{cor}

\begin{proof}
For $T=\eta$, the map $q_T$ is the identity of a point. For every
nontrivial tree $T$, it is a finite product of shift maps and hence a
weak homotopy equivalence.
\end{proof}

To equip framed configurations with a cyclic structure, we return to
the geometric models of the cyclic framed little disks operad
recalled in Example~\ref{ex: Budney}. Two such models will be useful:
the operad $\mathsf S_0$ of parametrised disk configurations in
$S^2$ and the framed Fulton--MacPherson operad $\mathsf{FFM}_2$,
whose cyclic structure is described in
\cite[Section~3.1]{campos2019configuration}. These geometric models
determine the cyclic enhancement of $X^{\mathrm{fr}}$ used below.

Their cyclic actions have compatible geometric descriptions. On
$\mathsf S_0$, the symmetric group permutes all parametrised disks,
including the disk labelled $0$. Choosing that disk as the output and
passing to planar coordinates identifies the transposition $(01)$
with the Möbius involution $z\mapsto-1/z$. The same involution
describes the transposition $(01)$ on $\mathsf{FFM}_2$. The following
proposition makes the equivalence of these cyclic structures precise.

\begin{prop}
\label{prop: cyclic comparison FFM2 S0}
The cyclic dendroidal nerves of $\mathsf{FFM}_2$ and $\mathsf S_0$
are equivalent as cyclic $\infty$-operads.
\end{prop}

\begin{proof}
After forgetting their cyclic structures, both operads are models for
the framed little disks operad. In particular, there are zigzags of
aritywise weak equivalences
\[
u^*\mathsf{FFM}_2
\simeq
\mathsf{FD}_2
\simeq
u^*\mathsf S_0.
\]
The spaces $\mathsf{FD}_2(n)$ are $K(\pi,1)$-spaces. Indeed,
$\mathsf{FD}_2(n)\simeq
\Conf_n(\mathbb C)\times(S^1)^n$, and
$\Conf_n(\mathbb C)$ is a $K(\PB_n,1)$. The same therefore holds for
$\mathsf{FFM}_2(n)$ and $\mathsf S_0(n)$.

We compare their cyclic structures through their fundamental
groupoids. Campos, Idrissi, and Willwacher identify the full
fundamental subgroupoid of $\mathsf{FFM}_2$ on the associahedral
vertices with $\PaRB$ (Definition~\ref{def: PaRB}); see
\cite[Section~3.1]{campos2019configuration}. Here the comparison uses
the version of $\PaRB$ with a terminal arity-zero component,
corresponding to the empty configuration. Under this identification,
the cyclic structure induced by the transposition $(01)$ is
represented, after normalising the point labelled $1$ to lie at the
origin, by the involution $z\mapsto-1/z$. It is therefore the cyclic
structure on $\PaRB^{\mathrm{cyc}}$ described in
Section~\ref{subsec: cyclic PaRB}.

Forgetting parenthesisations gives an equivalence of cyclic operads in
groupoids
\[
\PaRB^{\mathrm{cyc}}
\longrightarrow
\mathsf{RBr}^{\mathrm{cyc}}.
\]
Müller and Woike compare $\mathsf{RBr}^{\mathrm{cyc}}$ with
$\Pi_1\mathsf S_0$ as cyclic operads in groupoids
\cite[Proposition~5.3]{cyclic_ribbon}. As observed above,
the transposition $(01)$ on $\mathsf S_0$, expressed in the planar
chart determined by the output, is given by the same involution
$z\mapsto-1/z$. The underlying symmetric-group actions are preserved
by these comparisons, and $\Sigma_{n+1}$ is generated by $\Sigma_n$
together with $(01)$. Agreement on this transposition therefore
identifies the full cyclic actions.

Since the arity spaces are $K(\pi,1)$-spaces, their classifying
fundamental groupoids recover their weak homotopy types. The preceding
comparisons therefore give a zigzag of entrywise weak equivalences
between the cyclic dendroidal nerves of $\mathsf{FFM}_2$ and
$\mathsf S_0$. Hence these nerves are equivalent as cyclic
$\infty$-operads.
\end{proof}

It remains to compare the underlying rooted part of this cyclic
$\infty$-operad with $X^{\mathrm{fr}}$. Recall from
Example~\ref{ex: Budney} that there is an aritywise weak equivalence
of operads
$\varphi\colon\mathsf{FD}_2\to u^*\mathsf S_0$. For every rooted tree
$T$, the map $\nerve_{\mathrm{op}}(\varphi)_T$ is the product of the
maps $\varphi_{|\mathrm{in}(v)|}$ over the vertices of $T$, and is
therefore a weak homotopy equivalence. Combining this observation
with Corollary~\ref{cor:nerve-FD2-eq-Xfr} gives the following result.

\begin{theorem}
\label{thm: cyclic infinity operad of framed configurations}
The cyclic dendroidal nerve
$\nerve_{\mathrm{cyc}}\mathsf S_0$ is a cyclic
$\infty$-operadic model for framed configuration spaces. More
precisely, there is a zigzag of entrywise weak equivalences of
semi-dendroidal spaces
\[
\left.
f^*\nerve_{\mathrm{cyc}}\mathsf S_0
\right|_{\Omega_{\mathrm{inj}}}
\ \xleftarrow{\ \sim\ }\
\left.
\nerve_{\mathrm{op}}(\mathsf{FD}_2)
\right|_{\Omega_{\mathrm{inj}}}
\ \xrightarrow{\ \sim\ }\
X^{\mathrm{fr}}.
\]
\end{theorem}

\begin{proof}
Compatibility of the ordinary and cyclic dendroidal nerves with root
elision gives
$f^*\nerve_{\mathrm{cyc}}\mathsf S_0
\cong\nerve_{\mathrm{op}}(u^*\mathsf S_0)$.
Under this identification, the left-hand map is the restriction of
$\nerve_{\mathrm{op}}(\varphi)$ to
$\Omega_{\mathrm{inj}}$. The right-hand map is the weak equivalence
$q$ of Corollary~\ref{cor:nerve-FD2-eq-Xfr}.
\end{proof}

\section{A cyclic operad of ribbon braids}
\label{sec: cyclic operad of ribbon braids}

Campos, Idrissi, and Willwacher observed in
\cite[Section~3.1]{campos2019configuration} that the operad of
parenthesised ribbon braids $\PaRB$ admits a cyclic structure. Their
construction is motivated by the involution $z\mapsto-1/z$ on
normalised configurations of framed points in $\mathbb C$. Müller and
Woike subsequently described the corresponding cyclic structure on
the ribbon-braid groupoid and compared it with the cyclic surface
model $\mathsf S_0$; see \cite[Section~5]{cyclic_ribbon} and
Example~\ref{ex: Budney}.

We first recall the necessary background on braids, ribbon braids, and
the operad $\PaRB$ in Sections~\ref{subsec: ribbon braids} and
\ref{subsec:operad PaRB}. In Section~\ref{subsec: cyclic PaRB}, we
describe the cyclic enhancement of $\PaRB$ explicitly and verify its
compatibility with the defining relations and operadic composition.

\subsection{Braids and ribbon braids}
\label{subsec: ribbon braids}

The configuration space of $n$ unordered points in the complex plane,
$\Conf_n(\mathbb C)/\Sigma_n$, has fundamental group the
\emph{braid group} on $n$ strands,
\[
\Br_n=\pi_1\bigl(\Conf_n(\mathbb C)/\Sigma_n,p\bigr).
\]
It has the standard presentation
\[
\Br_n=
\left\langle
\beta_1,\ldots,\beta_{n-1}
\,\middle|\,
\beta_i\beta_j=\beta_j\beta_i\text{ if }|i-j|\geq2,\quad
\beta_i\beta_{i+1}\beta_i
=\beta_{i+1}\beta_i\beta_{i+1}
\right\rangle,
\]
where the second relation holds for $1\leq i\leq n-2$. The generator
$\beta_i$ represents the elementary exchange in which the strand in
position $i$ crosses over the strand in position $i+1$.

For each $n\geq1$, there is a short exact sequence
\[
\begin{tikzcd}
1 \arrow[r]
& \PB_n \arrow[r]
& \Br_n \arrow[r,"\pi"]
& \Sigma_n \arrow[r]
& 1,
\end{tikzcd}
\]
where $\pi$ records the induced permutation of the strands. Its kernel
is the \emph{pure braid group}
$\PB_n=\pi_1(\Conf_n(\mathbb C),p)$. It is generated by the standard
pure braids
\[
x_{ij}
=
\beta_{j-1}\cdots\beta_{i+1}\beta_i^2
\beta_{i+1}^{-1}\cdots\beta_{j-1}^{-1},
\qquad 1\leq i<j\leq n.
\]

Similarly, the configuration space of unordered framed points in the
plane has fundamental group the \emph{ribbon braid group}
\[
\RB_n
=
\pi_1\bigl(
\Conf_n^{\mathrm{fr}}(\mathbb C)/\Sigma_n,p
\bigr).
\]
It is generated by the braid generators
$\beta_1,\ldots,\beta_{n-1}$ and the framing twists
$\twist_1,\ldots,\twist_n$. In addition to the braid relations, these
generators satisfy
\[
\twist_i\twist_j=\twist_j\twist_i,\qquad \beta_i\twist_j=\twist_j\beta_i \quad\text{if }j\notin\{i,i+1\},
\]
and
\[
\beta_i\twist_i=\twist_{i+1}\beta_i,\qquad
\beta_i\twist_{i+1}=\twist_i\beta_i.
\]
Consequently, $\RB_n\cong\Br_n\ltimes\mathbb Z^n$, where the action of $\Br_n$ on $\mathbb Z^n$ factors through $\Br_n\to\Sigma_n$ and permutes the framing coordinates.

The kernel of the projection $\RB_n\to\Sigma_n$ is the \emph{pure ribbon braid group} $\PRB_n = \pi_1\bigl(\Conf_n^{\mathrm{fr}}(\mathbb C),p\bigr) \cong \PB_n\times\mathbb Z^n.$
 
\subsection{The operad of parenthesised ribbon braids}
\label{subsec:operad PaRB}

The operad of \emph{coloured ribbon braids}, $\CoRB$, is an operad in
groupoids whose automorphism groups are the pure ribbon braid groups
$\PRB_n$. Its classifying space provides a model for the
positive-arity framed little disks operad: there is a weak equivalence
of operads $\cs\CoRB\simeq\mathsf{FD}_2$; see
\cite[Proposition~6.8]{Boavida-Horel-Robertson} and
\cite[Chapter~1]{Wahl_Thesis}. The operad of
\emph{parenthesised ribbon braids}, $\PaRB$, is a cofibrant replacement
of $\CoRB$ and hence gives a homotopically well-behaved groupoid model
for $\mathsf{FD}_2$.

\begin{definition}
Set $\CoRB(0)=\varnothing$, where $\varnothing$ denotes the empty
groupoid. For each $n\geq1$, define the groupoid $\CoRB(n)$ as follows:
\begin{itemize}
\item the objects of $\CoRB(n)$ are the elements of $\Sigma_n$;
\item for $\sigma_1,\sigma_2\in\Sigma_n$, the set
$\Hom_{\CoRB(n)}(\sigma_1,\sigma_2)$ consists of the ribbon braids
whose underlying permutation is $\sigma_2\sigma_1^{-1}$.
\end{itemize}
If $r\colon\sigma_1\to\sigma_2$ and
$s\colon\sigma_2\to\sigma_3$, their categorical composite is denoted
$s\cdot r$ and is given by concatenation of ribbon braids.

Using the group isomorphism
$\RB_n\cong\Br_n\ltimes\mathbb Z^n$, a morphism
$r\in\Hom_{\CoRB(n)}(\sigma_1,\sigma_2)$ may be written as
$r=(\gamma,[\tau_1,\ldots,\tau_n])$, where
$\gamma\in\Br_n$ has underlying permutation
$\pi(\gamma)=\sigma_2\sigma_1^{-1}$ and
$\tau_i\in\mathbb Z$ records the number of twists of the $i$th ribbon
strand.
\end{definition}

The symmetric group $\Sigma_n$ acts naturally on $\CoRB(n)$, and the
collection $\CoRB=\{\CoRB(n)\}_{n\geq0}$ forms a symmetric operad in
groupoids. The partial composition functor
\[
\circ_i\colon
\CoRB(n)\times\CoRB(m)\longrightarrow\CoRB(n+m-1)
\]
is given on objects by the usual operadic composition of permutations.
Thus $\sigma\circ_i\rho$ is obtained by inserting
$\rho\in\Sigma_m$ into the $i$th position of
$\sigma\in\Sigma_n$ and reindexing. See
\cite[Section~5.2]{FresseBook1} or
\cite[Definition~6.1]{Boavida-Horel-Robertson}.

On morphisms, operadic composition inserts a ribbon braid on $m$
strands into the $i$th strand of a ribbon braid on $n$ strands. For
$m\geq1$, let
$\mathbf R_m=(\beta_1\cdots\beta_{m-1})^m\in\PB_m$ denote the full
twist, with $\mathbf R_1=1$. Given
$r_1=(\gamma_1,[\tau_1,\ldots,\tau_n])\in\CoRB(n)$ and
$r_2=(\gamma_2,[\tau'_1,\ldots,\tau'_m])\in\CoRB(m)$, their operadic
composite is
\[
r_1\circ_i r_2 = \Bigl( \gamma_1\circ_i \bigl((\mathbf R_m)^{\tau_i}\cdot\gamma_2\bigr), [\tau_1,\ldots,\tau_{i-1}, \tau_i+\tau'_1,\ldots,\tau_i+\tau'_m, \tau_{i+1},\ldots,\tau_n] \Bigr).
\]
Here $(\mathbf R_m)^{\tau_i}$ denotes the $\tau_i$th integer power of
the full twist under categorical composition, regarded as an
automorphism in $\CoRB(m)$. See
\cite[Definition~6.4]{Boavida-Horel-Robertson} and
\cite[Section~1.5]{Wahl_Thesis}.

\begin{definition}
\label{def:magma_operad}
We define the sets $\magma(n)$ as follows. Set
$\magma(0)=\varnothing$, and let $\magma(1)$ consist of the exceptional
one-edge tree $\lvert$, representing the operadic unit. For $n\geq2$,
let $\magma(n)$ be the set of isomorphism classes of planar rooted
binary trees with $n$ input edges, equipped with a boundary labelling
$\ell_T\colon\{0,1,\ldots,n\}\xrightarrow{\cong}\partial(T)$ such that
$\ell_T(0)$ is the root.

These sets form an operad in sets, called the \emph{magma operad}.
Operadic composition is given by grafting the root of one tree onto a
labelled input edge of another.

Equivalently, $\magma$ is the free operad generated by the
$\Sigma_2$-set consisting of the two rooted binary corollas
\[
\mu=\mu(1,2;0) 
\qquad \text{and} \qquad
\mu_{21}:=(12)^*\mu=\mu(2,1;0).
\]
These correspond to the two possible orders of the input edges after the root of the binary corolla has been fixed.
\end{definition}

\begin{remark}
\label{rem: permutation conventions PaRB}
We follow the geometric labelling convention of \cite[Section~3.1]{campos2019configuration}. Permutations of boundary labels are written in cycle notation and are composed from right to left. The symmetric actions are written on the right, so
\[
(\sigma\rho)^*=\rho^*\sigma^*.
\]
By contrast, a string $i_1\cdots i_n$ records an ordered list of labels. If $x$ is an object or morphism whose inputs are labelled $1,\ldots,n$, we write $x_{i_1\cdots i_n}$ for the result of relabelling its inputs, in order, by $i_1,\ldots,i_n$.

Thus $\alpha_{231}$ denotes the associator $\alpha_{123}$ with its input labels reordered as $2,3,1$; the string $231$ is not cycle notation. This differs from the one-line permutation convention used in \cite{Boavida-Horel-Robertson}, where the same relabelling is written $(231)^*\alpha$.
\end{remark}

Elements of $\magma(n)$ may equivalently be described as fully
parenthesised words in the symbols $\{1,\ldots,n\}$, in which every
symbol occurs exactly once. We use typewriter font for these words to
distinguish their parentheses from cycle notation. Under this
identification, $\mu$ and $\mu_{21}$ correspond to
$\mathtt{(12)}$ and $\mathtt{(21)}$, respectively. For example,
\[
\mathtt{((12)3)4},
\qquad
\mathtt{(32)(14)},
\qquad
\mathtt{1((23)4)}
\]
are elements of $\magma(4)$.

For each $n\geq1$, there is a natural map $\omega_n\colon\magma(n)\to\Sigma_n$ which forgets the tree and remembers only the permutation encoded by the ordering of its leaf labels. In terms of parenthesised words, one removes the parentheses and regards the resulting word as the ordered list of images of $1,\ldots,n$. For example, removing the parentheses from $\mathtt{(32)(14)}$ gives $3214$, and hence $\omega_4(\mathtt{(32)(14)})=(13)$ in cycle notation. The maps $\omega_n$ are compatible with the symmetric-group actions and operadic composition, and therefore assemble into an operad morphism $\omega\colon\magma\to\ob(\CoRB)$.

\begin{definition}
\label{def: PaRB}
Set $\PaRB(0)=\varnothing$. For each $n\geq1$, define the groupoid
$\PaRB(n)$ as follows:
\begin{itemize}
\item the objects of $\PaRB(n)$ are the labelled trees $T\in\magma(n)$;
\item for $T_1,T_2\in\magma(n)$, set
\[
\Hom_{\PaRB(n)}(T_1,T_2) = \Hom_{\CoRB(n)}\bigl(\omega_n(T_1),\omega_n(T_2)\bigr).
\]
Thus a morphism $T_1\to T_2$ is a ribbon braid whose underlying permutation is $\omega_n(T_2)\omega_n(T_1)^{-1}$.
\end{itemize}

Categorical composition is inherited from $\CoRB(n)$. The collection $\PaRB=\{\PaRB(n)\}_{n\geq0}$ forms an operad in groupoids. On objects, operadic composition is given by grafting labelled rooted trees, as in $\magma$, and on morphisms it is given by insertion of ribbon braids, as in $\CoRB$.
\end{definition}

\begin{figure}[htbp]
    \centering
    \includegraphics[width=0.5\linewidth]{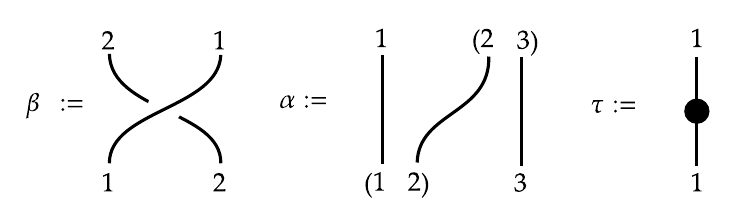}
  \caption{The generating braiding $\beta$, associator $\alpha$, and
twist $\tau$ in $\PaRB$. Morphisms are drawn from the bottom of the page
to the top.}
    \label{fig:Generating morphisms of PaRB}
\end{figure}

The object operad $\ob(\PaRB)=\magma$ is generated by the binary corolla $\mu$ and its symmetric translate $\mu_{21}=(12)^*\mu$. Morphisms in $\PaRB$ are generated under operadic and categorical composition by
\[
\beta_{12}\in \Hom_{\PaRB(2)}(\mu,\mu_{21}),
\qquad
\twist\in\Hom_{\PaRB(1)}(\lvert,\lvert),
\]
together with the associator
\[
\alpha_{123}\in \Hom_{\PaRB(3)}\bigl(\mu\circ_1\mu,\mu\circ_2\mu\bigr)
\]
and their inverses. Here $\beta_{12}$ is the standard braid generator and $\twist$ is the twist generator of $\RB_1\cong\mathbb Z$, represented pictorially by $\tikz{\fill[black] (3,-0.25) circle[radius=1mm]}$ on a ribbon strand to denote a positive twist. We will denote the inverse of the twist operator $\tau^{-1}$ by replacing $\tikz{\fill[black] (3,-0.25) circle[radius=1mm]}$ in the picture by $\tikz{\draw[fill=white] (3,-0.25) circle[radius=1mm]}$ symbol. We generally write $\beta=\beta_{12}$ and $\alpha=\alpha_{123}$ when the labelling is clear.

\begin{thm}[
\textnormal{\cite[Lemmas~7.1 and~7.4]{Boavida-Horel-Robertson}};
\textnormal{\cite[Theorem~6.2.4]{FresseBook1}}]
\label{thm: maps out of PaRB}
Let $\calQ$ be an operad in groupoids, and denote its operadic unit by $\mathbb{1}_{\calQ}\in\ob(\calQ(1))$. To give an operad morphism $f\colon\PaRB\to\calQ$ is equivalent to choosing an object $\mathtt m\in\ob(\calQ(2))$ and morphisms
\[
\mathtt t\colon
\mathbb{1}_{\calQ}\longrightarrow\mathbb{1}_{\calQ},
\qquad
\mathtt b_{12}\colon
\mathtt m\longrightarrow\mathtt m_{21},
\qquad
\mathtt a_{123}\colon
\mathtt m\circ_1\mathtt m
\longrightarrow
\mathtt m\circ_2\mathtt m,
\]
in arities $1$, $2$, and $3$, respectively, where
$\mathtt m_{21}=(12)^*\mathtt m$. Set
\[
\mathtt b_{21}:=(12)^*\mathtt b_{12}\colon
\mathtt m_{21}\longrightarrow\mathtt m.
\]
These data determine an operad morphism if and only if they satisfy the following relations:
\begin{equation}
\label{twist_rel_morph}\tag{T}
\mathtt t\circ_1 1_{\mathtt m}
=
\mathtt b_{21}\cdot\mathtt b_{12}\cdot
\bigl(1_{\mathtt m}\circ(\mathtt t,\mathtt t)\bigr);
\end{equation}
\begin{equation}
\label{hex_rel_1_morph}\tag{H1}
(1_{\mathtt m}\circ_2\mathtt b_{12})_{213}
\cdot\mathtt a_{213}
\cdot(1_{\mathtt m}\circ_1\mathtt b_{12})
=
\mathtt a_{231}
\cdot(\mathtt b_{12}\circ_2 1_{\mathtt m})
\cdot\mathtt a_{123};
\end{equation}
\begin{equation}
\label{hex_rel_2_morph}\tag{H2}
(1_{\mathtt m}\circ_1\mathtt b_{12})_{132}
\cdot\mathtt a_{132}^{-1}
\cdot(1_{\mathtt m}\circ_2\mathtt b_{12})
=
\mathtt a_{312}^{-1}
\cdot(\mathtt b_{12}\circ_1 1_{\mathtt m})
\cdot\mathtt a_{123}^{-1};
\end{equation}
and
\begin{equation}
\label{pent_rel_1_morph}\tag{P}
(\mathtt a_{123}\circ_3 1_{\mathtt m})
\cdot(\mathtt a_{123}\circ_1 1_{\mathtt m})
=
(1_{\mathtt m}\circ_2\mathtt a_{123})
\cdot(\mathtt a_{123}\circ_2 1_{\mathtt m})
\cdot(1_{\mathtt m}\circ_1\mathtt a_{123}).
\end{equation}

Here $1_{\mathtt m}$ denotes the identity morphism of $\mathtt m$, and $1_{\mathtt m}\circ(\mathtt t,\mathtt t)$ denotes the simultaneous substitution of $\mathtt t$ into its two inputs. Categorical composition is denoted by $\cdot$ and uses the standard convention: if $f\colon x\to y$ and $g\colon y\to z$, then $g\cdot f\colon x\to z$ means first applying $f$ and then $g$. The symbol $\circ_i$ denotes operadic partial composition.

Subscripts indicate ordered relabellings in the sense of Remark~\ref{rem: permutation conventions PaRB}. Thus $\mathtt a_{231}$ is obtained from $\mathtt a_{123}$ by relabelling its inputs, in order, by $2,3,1$.
\end{thm}

\begin{remark}
In $\PaRB$, the ribbon-braid relations allow the individual strand
twists to pass through the pure double braiding. Consequently, the
right-hand side of \eqref{twist_rel_morph} may equivalently be written
\[
\bigl(1_{\mathtt m}\circ(\mathtt t,\mathtt t)\bigr)
\cdot\mathtt b_{21}\cdot\mathtt b_{12} = \mathtt b_{21} \cdot \bigl(1_{\mathtt m_{21}}\circ(\mathtt t,\mathtt t)\bigr)
\cdot\mathtt b_{12}.
\]
The same identity holds after applying any operad morphism
$\PaRB\to\calQ$.
\end{remark}

\begin{remark}
The operad $\PaRB$ may equivalently be obtained from $\CoRB$ by
pulling back its object set along
$\omega\colon\magma\to\ob(\CoRB)$; see
\cite[Section~6.1.5]{FresseBook1} for pullbacks of operads in
groupoids. In other words,
\[
\PaRB\cong\omega^*\CoRB,
\]
where $\omega^*\CoRB$ has object set $\magma(n)$ in arity $n$ and
morphism sets
\[
\Hom_{\omega^*\CoRB(n)}(T_1,T_2)
=
\Hom_{\CoRB(n)}\bigl(\omega_n(T_1),\omega_n(T_2)\bigr).
\]
Forgetting the parenthesisation gives an entrywise equivalence
$\PaRB\to\CoRB$, and in positive arities there is a zigzag of
entrywise weak equivalences
\[
\PaRB\xrightarrow{\simeq}\CoRB
\xrightarrow{\simeq}\Pi_1(\mathsf{FD}_2).
\]
Moreover, $\PaRB$ is cofibrant in $\Op(\Grpd)$; see
\cite[Corollary~6.14]{Boavida-Horel-Robertson}. For the comparison
with the framed little disks operad, see also
\cite[Chapter~1]{Wahl_Thesis}.
\end{remark}

The relations in Theorem~\ref{thm: maps out of PaRB} are precisely the
coherence relations for balanced monoidal categories introduced in
Section~\ref{subsec: tensor categories}. This is explained by the
following theorem of Wahl.

\begin{thm}[\textnormal{\cite[Theorem~II.5.11]{Wahl_Thesis}}]
\label{thm: CoRB_and_PaRB_algebras}
In the category of small categories $\mathbf{Cat}$, algebras over $\CoRB$ are precisely strict balanced monoidal categories, while algebras over $\PaRB$ are balanced monoidal categories.
\end{thm}

\subsection{The cyclic operad of parenthesised ribbon braids}
\label{subsec: cyclic PaRB}
The operad $\PaRB$ admits a natural cyclic enhancement. Formulas for
this structure are given in
\cite[Section~3.1]{campos2019configuration}, motivated by the
inversion map $z\mapsto-1/z$ on normalised configurations of framed
points in $\mathbb C$. This involution exchanges the output with the
first input, placing all boundary labels on an equal footing. The
corresponding cyclic structure on the ribbon-braid groupoid is
compared with the geometric cyclic structure on the framed little
disks operad in \cite[Section~5]{cyclic_ribbon}. In
\cite[Proposition~3.15]{BR_modular}, this cyclic operad is identified
with the genus-zero truncation of the modular operad of surfaces
constructed from mapping class groups.

We first describe the cyclic structure on objects and then specify its action on the generating morphisms of $\PaRB$.

\begin{definition}
\label{def: objects of cyclic PaRB}
The cyclic operad of objects $\ob(\PaRB^{\mathrm{cyc}})$ is the free
cyclic operad generated by the right $\Sigma_2^+$-set of left cosets
\[
A_2^+\backslash\Sigma_2^+
=
\{A_2^+\sigma \mid \sigma\in\Sigma_2^+\},
\qquad
A_2^+=\langle(012)\rangle \subset\Sigma_2^{+},
\]
where $\Sigma_2^+$ acts on the right by $A_2^+\sigma\cdot\tau =
A_2^+(\sigma\tau)$. This coset space has two elements, represented,
after choosing the boundary edge labelled $0$ as the root, by the
rooted binary corollas
\[
\mu=\mu(1,2;0)
\qquad\text{and}\qquad
\mu_{21}:=(12)^*\mu=\mu(2,1;0).
\]
The transpositions $(01)$ and $(12)$ lie in the same left coset of
$A_2^+$, while $(012)\in A_2^+$. Consequently,
\[
(01)^*\mu=(12)^*\mu=\mu_{21},
\qquad
(012)^*\mu=\mu.
\]
Forgetting the cyclic structure recovers the magma operad:
$u^*\ob(\PaRB^{\mathrm{cyc}})\cong\magma$.
\end{definition}

Concretely, for $n \geq 2$, the objects of $\PaRB^{\mathrm{cyc}}(n)$
may be represented by planar unrooted trivalent trees whose boundary
edges are labelled by $\{0,\ldots,n\}$. The planar structure
determines a cyclic order on the three half-edges incident to each
vertex. Since the generator of $A_2^+ = \langle (012) \rangle$ acts
by cyclically rotating these half-edges, the two elements of
$A_2^+ \backslash \Sigma_2^+$ correspond to the two distinct
planar binary corollas $\mu$ and $\mu_{21}$. 

We now specify the extended symmetric-group actions on morphisms.
Since $\Sigma_n^+$ is generated by the subgroup $\Sigma_n$ fixing
$0$ together with the transposition $s_0=(01)$, it suffices to
specify the additional action on the generating morphisms of
$\PaRB$.

\begin{definition}
\label{def: cyclic PaRB}
The cyclic operad $\PaRB^{\mathrm{cyc}}$ is the operad $\PaRB$
equipped with the cyclic structure on objects given in
Definition~\ref{def: objects of cyclic PaRB} and the action on
generating morphisms prescribed by
\begin{equation}
\label{eq: cyclic action PaRB generators}
s_0^*\twist = \twist,
\qquad
s_0^*\beta_{12}
= \bigl(1_\mu\circ_2\twist^{-1}\bigr)\cdot\beta_{12}^{-1},
\qquad
s_0^*\alpha_{123} = \alpha_{231}^{-1}.
\end{equation}
\end{definition}

\begin{figure}[h]
    \centering
    \includegraphics[width=0.8\linewidth]{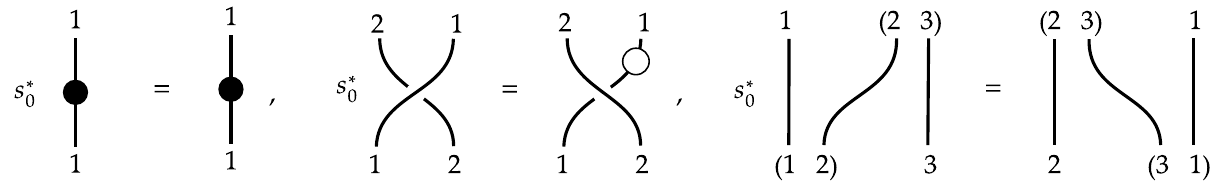}
    \caption{The cyclic action on the generating morphisms of $\PaRB$,
as given in \eqref{eq: cyclic action PaRB generators}.}
    \label{fig:cyclic_action_on_PaRB}
\end{figure}

Here, a filled circle $\tikz{\fill[black] (3,-0.25) circle[radius=1mm]}$ on a ribbon strand represents a positive twist, while an open circle $\tikz{\draw[fill=white] (3,-0.25) circle[radius=1mm]}$ represents a negative twist.

\medskip

These formulas are consistent with the cyclic structures described in
\cite[Section~3.1]{campos2019configuration} and
\cite[Section~5]{cyclic_ribbon}. Here the subscript $231$ denotes the
ordered relabelling of
Remark~\ref{rem: permutation conventions PaRB}, rather than a
permutation written in cycle notation. In
Lemma~\ref{lemma: extended symmetric action on generators}, we check
that the displayed morphisms have the required sources and targets
and satisfy the extended symmetric-group relations.
Proposition~\ref{prop: cyclic structure on PaRB} then establishes
compatibility with operadic composition.

\begin{lemma}
\label{lemma: extended symmetric action on generators}
The assignments in \eqref{eq: cyclic action PaRB generators},
together with the ordinary symmetric-group actions, satisfy the
extended symmetric-group relations on the generating morphisms
$\twist$, $\beta_{12}$, and $\alpha_{123}$.
\end{lemma}

Before giving the proof, we fix notation used throughout the remainder of this section. For $n \geq 1$, write
\[
s_0 = (01)
\qquad\text{and}\qquad
s_i = (i\ i+1), \quad 1 \leq i < n,
\]
for the adjacent transpositions generating $\Sigma_n^+ = \operatorname{Aut}(\{0,\ldots,n\})$. They satisfy the
Coxeter relations
\[
s_i^2 = 1, \qquad
s_i s_j = s_j s_i \quad \text{if } |i-j| \geq 2, \qquad
s_i s_{i+1} s_i = s_{i+1} s_i s_{i+1}.
\]
The generators $s_1, \ldots, s_{n-1}$ generate the subgroup $\Sigma_n \subseteq \Sigma_n^+$ fixing $0$. Since the
symmetric-group actions are right actions, we have $(\sigma\rho)^* = \rho^*\sigma^*$. The ordinary $\Sigma_n$-actions
already satisfy all Coxeter relations among $s_1, \ldots, s_{n-1}$. It therefore remains only to verify the three Coxeter relations involving $s_0$, namely
\[
s_0^2 = 1, \qquad
s_0 s_2 = s_2 s_0 \quad \text{(in arity } \geq 3\text{)},
\qquad
s_0 s_1 s_0 = s_1 s_0 s_1,
\]
on each of the generating morphisms.

\begin{proof}[Proof of Lemma~\ref{lemma: extended symmetric action on generators}]
In arity one, the extended symmetric group $\Sigma_1^+\cong\Sigma_2$
is generated by $s_0$. The only relation to verify is $s_0^2=1$,
which follows immediately from $s_0^*\twist=\twist$.

In arity two, the extended symmetric group $\Sigma_2^+\cong\Sigma_3$
is generated by $s_0=(01)$ and $s_1=(12)$. We must verify that
$(s_0^*)^2\beta_{12}=\beta_{12}$ and
$(s_0s_1s_0)^*\beta_{12}=(s_1s_0s_1)^*\beta_{12}$.

We first check that $s_0^*\beta_{12}=(1_\mu\circ_2\twist^{-1})
\cdot\beta_{12}^{-1}$ has the correct source and target. Since
$\beta_{12}\colon\mu\to\mu_{21}$ and $s_0^*\mu=\mu_{21}$,
$s_0^*\mu_{21}=\mu$ by
Definition~\ref{def: objects of cyclic PaRB}, the morphism
$s_0^*\beta_{12}$ must be a morphism from $\mu_{21}$ to $\mu$.
The morphism $(1_\mu\circ_2\twist^{-1})\cdot\beta_{12}^{-1}$
has source $\mu_{21}$ and target $\mu$, as required.

The two Coxeter relations may be checked directly from the
ribbon-braid relations. Equivalently, set
$z_3=(012)=s_0s_1$. Since the action is on the right,
$z_3^*=s_1^*s_0^*$, and the formula above gives
\begin{equation}
\label{z3beta}
z_3^*\beta_{12}
=\bigl(1_{\mu_{21}}\circ_1\twist^{-1}\bigr)
 \cdot\beta_{21}^{-1}.
\end{equation}
Multiplying the twist equation \eqref{twist_rel_morph} by
$\beta_{21}^{-1}$ on the left and
$\twist^{-1}\circ_1 1_\mu$ on the right gives
\begin{equation}
\label{beta21inverse}
\beta_{21}^{-1} =\beta_{12}  \cdot\bigl(1_\mu\circ_1\twist\bigr)  \cdot\bigl(1_\mu\circ_2\twist\bigr) \cdot\bigl(\twist^{-1}\circ_1 1_\mu\bigr).
\end{equation}

Sliding the twists past the braiding in \eqref{z3beta} and then substituting
\eqref{beta21inverse}, we obtain
\begin{equation}\label{z3squaredbeta}
\begin{aligned}
z_3^*\beta_{12}
&=\beta_{21}^{-1}
  \cdot\bigl(1_\mu\circ_1\twist^{-1}\bigr) =\underbrace{\beta_{12}\cdot\bigl(1_\mu\circ_1\twist\bigr)\cdot\bigl(1_\mu\circ_2\twist\bigr)\cdot\bigl(\twist^{-1}\circ_1 1_\mu\bigr)}_{\eqref{beta21inverse}}
  \cdot\bigl(1_\mu\circ_1\twist^{-1}\bigr) =\beta_{12}
  \cdot\bigl(\twist^{-1}\circ_1 1_\mu\bigr)
  \cdot\bigl(1_\mu\circ_2\twist\bigr)
  \end{aligned}
\end{equation}

The last equality follows by commuting twists cancelling the inverse twists on the first strand.
Applying the cyclic composition rule and substituting the expression
for $z_3^*\beta_{12}$ gives

\begin{equation}\label{z3squarebeta'}
\begin{aligned}
(z_3^*)^2\beta_{12}
&=(z_3^*\beta_{12})
  \cdot\bigl(1_\mu\circ_2\twist^{-1}\bigr)
  \cdot\bigl(1_\mu\circ_1\twist\bigr) =\underbrace{\beta_{12}\cdot\bigl(\twist^{-1}\circ_1 1_\mu\bigr)
  \cdot\bigl(1_\mu\circ_2\twist\bigr)}_{\eqref{z3squaredbeta}}
  \cdot\bigl(1_\mu\circ_2\twist^{-1}\bigr)
  \cdot\bigl(1_\mu\circ_1\twist\bigr) =\beta_{12}\cdot\bigl(\twist^{-1}\circ_1 1_\mu\bigr)
  \cdot\bigl(1_\mu\circ_1\twist\bigr).
\end{aligned}
\end{equation}
Where the last equality comes from cancelling the twists on the second strand. 
Applying the same rule once more, we obtain
\begin{equation}
\begin{aligned}
(z_3^*)^3\beta_{12}
&=(z_3^*\beta_{12})
  \cdot\bigl(1_\mu\circ_2\twist^{-1}\bigr)
  \cdot\bigl(\twist\circ_1 1_\mu\bigr)\\
&=\underbrace{\beta_{12}\cdot\bigl(\twist^{-1}\circ_1 1_\mu\bigr)
  \cdot\bigl(1_\mu\circ_2\twist\bigr)}_{\eqref{z3beta}}
  \cdot\bigl(1_\mu\circ_2\twist^{-1}\bigr)
  \cdot\bigl(\twist\circ_1 1_\mu\bigr)\\
&=\beta_{12}\cdot\bigl(\twist^{-1}\circ_1 1_\mu\bigr)
  \cdot\bigl(\twist\circ_1 1_\mu\bigr)\\
&=\beta_{12}.
\end{aligned}
\end{equation} Here the last two equalities come from cancelling the twists on strand $2$ and then strand $1$. One can similarly check directly that $(z_3s_1)^*\beta_{12}=(s_1z_3^2)^*\beta_{12}$. Together with
$s_1^2=1$, these are equivalent to the two Coxeter relations above.
A similar verification has taken place in \cite[Lemma 5.1]{cyclic_ribbon}. This means that $z_3^*$ extends the $\Sigma_2$-action on $\PaRB(2)$ to a $\Sigma_2^+$-action.
\medskip

We now consider the associator $\alpha_{123}\colon\mu\circ_1\mu\longrightarrow\mu\circ_2\mu$ in $\PaRB(3)$. The extended symmetric group $\Sigma_3^+\cong\Sigma_4$
is generated by $s_0,s_1,s_2$. The relations involving $s_0$ are
\[
s_0^2=1,\qquad s_0s_2=s_2s_0,\qquad s_0s_1s_0=s_1s_0s_1.
\]
For an ordered string $ijk$ obtained by rearranging $123$, let
$\alpha_{ijk}$ denote the associator $\alpha_{123}$ with its input
labels reordered as $i,j,k$, as in
Remark~\ref{rem: permutation conventions PaRB}. Thus
$\alpha_{123}=\alpha$. The action of $s_0=(01)$ on the associator is
\[
s_0^*\alpha_{123}=\alpha_{231}^{-1}.
\]
This is obtained by exchanging the boundary labels $0$ and $1$ and
then regarding the edge newly labelled $0$ as the root. This change
of root reverses the corresponding $A$-move and cyclically reorders
the remaining boundary labels, as shown in
Figure~\ref{fig:cyclic action on associator}.

Applying the same rerooting operation to the symmetric translates of
$\alpha_{123}$ gives
\begin{equation}
\label{eq: rerooting associator translates}
\tag{\(\dagger\)}
\begin{aligned}
s_0^*\alpha_{123}&=\alpha_{231}^{-1}, &
s_0^*\alpha_{231}&=\alpha_{123}^{-1},\\
s_0^*\alpha_{213}&=\alpha_{312}, &
s_0^*\alpha_{312}&=\alpha_{213},\\
s_0^*\alpha_{132}&=\alpha_{321}^{-1}, &
s_0^*\alpha_{321}&=\alpha_{132}^{-1}.
\end{aligned}
\end{equation}
These identities follow uniformly from the action of $s_0$. It
exchanges the external labels $0$ and $1$ and, after the resulting
tree is written with $0$ as its output, cyclically shifts the
remaining labels. This rerooting reverses the direction of the
$A$-move relating $\alpha_{123}$ to $\alpha_{231}$ and
$\alpha_{132}$ to $\alpha_{321}$, accounting for the inverses in
these two cases. It exchanges $\alpha_{213}$ and $\alpha_{312}$
without reversing the corresponding $A$-move. Finally, functoriality
of the relabelling action on each groupoid gives
\[
s_0^*(\alpha_{ijk}^{-1})
=
\bigl(s_0^*\alpha_{ijk}\bigr)^{-1}.
\]

Using the first pair of identities in
\eqref{eq: rerooting associator translates}, we obtain
\[
(s_0^*)^2\alpha_{123}
=s_0^*\bigl(\alpha_{231}^{-1}\bigr)
=\bigl(s_0^*\alpha_{231}\bigr)^{-1}
=\bigl(\alpha_{123}^{-1}\bigr)^{-1}
=\alpha_{123}.
\]
Thus the relation $s_0^2=1$ holds on $\alpha_{123}$.

Next consider the commutation relation $s_0s_2=s_2s_0$. Since the
action is on the right,
\[
(s_0s_2)^*\alpha_{123}
=s_2^*s_0^*\alpha_{123}
=s_2^*\bigl(\alpha_{231}^{-1}\bigr)
=\alpha_{321}^{-1}.
\]
On the other hand,
\[
(s_2s_0)^*\alpha_{123}
=s_0^*s_2^*\alpha_{123}
=s_0^*\alpha_{132}
=\alpha_{321}^{-1},
\]
where the final equality follows from
\eqref{eq: rerooting associator translates}. Hence
$(s_0s_2)^*\alpha_{123}=(s_2s_0)^*\alpha_{123}$.

Finally, consider the braid relation $s_0s_1s_0=s_1s_0s_1$. For the
left-hand side,
\[
(s_0s_1s_0)^*\alpha_{123}
=s_0^*s_1^*s_0^*\alpha_{123}
=s_0^*s_1^*\bigl(\alpha_{231}^{-1}\bigr)
=s_0^*\alpha_{132}^{-1}
=\bigl(s_0^*\alpha_{132}\bigr)^{-1}
=\bigl(\alpha_{321}^{-1}\bigr)^{-1}
=\alpha_{321}.
\]
For the right-hand side,
\[
(s_1s_0s_1)^*\alpha_{123}
=s_1^*s_0^*s_1^*\alpha_{123}
=s_1^*s_0^*\alpha_{213}
=s_1^*\alpha_{312}
=\alpha_{321}.
\]
Hence $(s_0s_1s_0)^*\alpha_{123}=(s_1s_0s_1)^*\alpha_{123}$,
as required.

We have verified all relations involving $s_0$ on each of the
generating morphisms. The remaining relations are already satisfied
by the ordinary symmetric-group actions. Hence the assignments in
\eqref{eq: cyclic action PaRB generators}, together with the ordinary
symmetric-group actions, define right actions of the appropriate
extended symmetric groups on $\twist$, $\beta_{12}$, and
$\alpha_{123}$.
\end{proof}

\begin{figure}[htbp]
\centering
\resizebox{0.7\textwidth}{!}{%
\begin{circuitikz}
\tikzstyle{every node}=[font=\normalsize]
\draw [short] (2,6.75) -- (2,6.25)
  node[pos=0.5,above,fill=white]{0};
\draw [short] (2,6.25) -- (1,5.25)
  node[pos=1,below,fill=white]{1};
\draw [short] (1.5,5.75) -- (2,5.25)
  node[pos=1,below,fill=white]{2};
\draw [short] (2,6.25) -- (3,5.25)
  node[pos=1,below,fill=white]{3};

\draw [short] (2,1.75) -- (2,1.25)
  node[pos=0.5,above,fill=white]{0};
\draw [short] (2,1.25) -- (1,0.25)
  node[pos=1,below,fill=white]{1};
\draw [short] (2,1.25) -- (3,0.25)
  node[pos=1,below,fill=white]{3};
\draw [short] (2.5,0.75) -- (2,0.25)
  node[pos=1,below,fill=white]{2};

\draw [short] (7,6.75) -- (7,6.25)
  node[pos=0.5,above,fill=white]{1};
\draw [short] (7,6.25) -- (6,5.25)
  node[pos=1,below,fill=white]{0};
\draw [short] (6.5,5.75) -- (7,5.25)
  node[pos=1,below,fill=white]{2};
\draw [short] (7,6.25) -- (8,5.25)
  node[pos=1,below,fill=white]{3};

\draw [short] (7,1.75) -- (7,1.25)
  node[pos=0.5,above,fill=white]{1};
\draw [short] (7,1.25) -- (6,0.25)
  node[pos=1,below,fill=white]{0};
\draw [short] (7,1.25) -- (8,0.25)
  node[pos=1,below,fill=white]{3};
\draw [short] (7.5,0.75) -- (7,0.25)
  node[pos=1,below,fill=white]{2};

\draw [short] (11.5,6.75) -- (11.5,6.25)
  node[pos=0.5,above,fill=white]{0};
\draw [short] (11.5,6.25) -- (10.5,5.25)
  node[pos=1,below,fill=white]{2};
\draw [short] (11.5,6.25) -- (12.5,5.25)
  node[pos=1,below,fill=white]{1};
\draw [short] (12,5.75) -- (11.5,5.25)
  node[pos=1,below,fill=white]{3};

\draw [short] (11.5,1.75) -- (11.5,1.25)
  node[pos=0.5,above,fill=white]{0};
\draw [short] (11.5,1.25) -- (10.5,0.25)
  node[pos=1,below,fill=white]{2};
\draw [short] (11,0.75) -- (11.5,0.25)
  node[pos=1,below,fill=white]{3};
\draw [short] (11.5,1.25) -- (12.5,0.25)
  node[pos=1,below,fill=white]{1};

\draw [->,>=Stealth] (2,3.75) -- (2,2.75)
  node[pos=0.5,left,fill=white]{$\alpha_{123}$};
\draw [->,>=Stealth] (11.5,3.75) -- (11.5,2.75)
  node[pos=0.5,right,fill=white]
  {$s_0^*\alpha_{123}=\alpha_{231}^{-1}$};

\draw [->,>=Stealth] (3.5,6.25) -- (5,6.25)
  node[pos=0.5,above,fill=white]{$s_0^*$};
\draw [->,>=Stealth] (3.5,1.25) -- (5,1.25)
  node[pos=0.5,above,fill=white]{$s_0^*$};

\node [font=\normalsize] at (9,6.25) {$=$};
\node [font=\normalsize] at (9,1.25) {$=$};
\end{circuitikz}
}
\caption{The action of $s_0=(01)$ on the associator
$\alpha_{123}\colon\mu\circ_1\mu\to\mu\circ_2\mu$.}
\label{fig:cyclic action on associator}
\end{figure}

The preceding lemma shows that the assignments in
\eqref{eq: cyclic action PaRB generators}, together with the ordinary
symmetric-group actions, define well-defined extended symmetric-group
actions on the generating morphisms of $\PaRB$. It remains to verify
that these actions are compatible with operadic composition.

For $n \geq 1$, let $z_{n+1} = (0\,1\,\cdots\,n) \in \Sigma_n^+$ denote the standard cyclic permutation of the boundary labels, so that $z_{n+1} = s_0 s_1 \cdots s_{n-1}$ in our convention for multiplying permutations. Compatibility with operadic composition requires that
\begin{equation}
\label{eq: cyclic compatibility PaRB}
z_{n+m}^*(x \circ_i y) =
\begin{cases}
z_{n+1}^*(x) \circ_{i-1} y, & 2 \leq i \leq n,\\[3pt]
z_{m+1}^*(y) \circ_m z_{n+1}^*(x), & i = 1,
\end{cases}
\end{equation}
for every $x \in \PaRB(n)$ and $y \in \PaRB(m)$, where $z_{n+m}$ denotes the standard cyclic generator of $\Sigma_{n+m-1}^+$. These are the standard equivariance identities of Definition~\ref{def: cyclic operad} relating the cyclic action to operadic partial composition. 

\begin{prop}
\label{prop: cyclic structure on PaRB}
The cyclic structure on $\ob(\PaRB)$ described in Definition~\ref{def: objects of cyclic PaRB}, together with the assignments in \eqref{eq: cyclic action PaRB generators}, extends uniquely to a cyclic operad structure on $\PaRB$. We denote the resulting cyclic operad by $\PaRB^{\mathrm{cyc}}$.
\end{prop}

\begin{proof}
Since $\PaRB$ is generated by $\mu$, $\twist$, $\beta_{12}$, and
$\alpha_{123}$ under symmetric-group actions and operadic and
categorical composition, the cyclic equivariance identities
\eqref{eq: cyclic compatibility PaRB} determine at most one extension
of the assignments in \eqref{eq: cyclic action PaRB generators} to
all of $\PaRB$. These assignments first extend, using the cyclic
equivariance identities, to the free cyclic operad generated by these
objects and morphisms. It therefore remains to verify that the
congruence generated by the defining relations of $\PaRB$ is
preserved.

By Lemma~\ref{lemma: extended symmetric action on generators}, the
assignments satisfy the extended symmetric-group relations on the
generators. Since the ordinary symmetric-group actions already
preserve the defining relations and $\Sigma_n^+$ is generated by
$\Sigma_n$ together with $z_{n+1}$, it suffices to check that the
twist, hexagon, and pentagon relations are preserved by $z_{n+1}$.

\smallskip
\noindent\emph{The pentagon relation.}
We begin by observing that
\[
z_3^*1_\mu=1_\mu
\qquad\text{and}\qquad
z_4^*\alpha_{123}=\alpha_{123}^{-1}.
\]
Indeed, $z_3=(012)$ cyclically permutes the three boundary labels of
the binary corolla and therefore fixes both $\mu$ and its identity
morphism. Moreover, since $z_4=s_0s_1s_2$ and the action is on the
right, we have
\[
z_4^*\alpha_{123}
=
s_2^*s_1^*s_0^*\alpha_{123}
=
s_2^*s_1^*\alpha_{231}^{-1}
=
\alpha_{123}^{-1}.
\]

We can now compute the effect of the cyclic permutation on each morphism appearing at the vertices of the pentagon diagram. Consider, for example,
$\alpha_{123}\circ_1 1_\mu$. Here $\alpha_{123}\in\PaRB(3)$, $1_\mu\in\PaRB(2)$, and the substitution takes place in the first input. The second case of \eqref{eq: cyclic compatibility PaRB} therefore gives \[z_5^*(\alpha_{123}\circ_1 1_\mu) = z_3^*(1_\mu)\circ_2z_4^*(\alpha_{123}) = 1_\mu\circ_2\alpha_{123}^{-1}. \] The remaining factors are computed in the same way. Substitutions into the second or third input use the first case of  \eqref{eq: cyclic compatibility PaRB}, while substitutions into the first input use the second. We obtain
\begin{equation}
\label{eq: cyclic action on pentagon factors}
\begin{aligned}
z_5^*(\alpha_{123}\circ_1 1_\mu)
  &=1_\mu\circ_2\alpha_{123}^{-1},
&
z_5^*(\alpha_{123}\circ_2 1_\mu)
  &=\alpha_{123}^{-1}\circ_1 1_\mu,\\
z_5^*(\alpha_{123}\circ_3 1_\mu)
  &=\alpha_{123}^{-1}\circ_2 1_\mu,
&
z_5^*(1_\mu\circ_1\alpha_{123})
  &=\alpha_{123}^{-1}\circ_3 1_\mu,\\
z_5^*(1_\mu\circ_2\alpha_{123})
  &=1_\mu\circ_1\alpha_{123}.
\end{aligned}
\end{equation}

Recall that the pentagon relation \eqref{pent_rel_1_morph} in Theorem~\ref{thm: maps out of PaRB}, specialized to the generators of
$\PaRB$, is
\[ (\alpha_{123}\circ_3 1_\mu) \cdot(\alpha_{123}\circ_1 1_\mu) = (1_\mu\circ_2\alpha_{123})
\cdot(\alpha_{123}\circ_2 1_\mu) \cdot(1_\mu\circ_1\alpha_{123}). \] Applying $z_5^*$ and using
\eqref{eq: cyclic action on pentagon factors}, we obtain
\begin{equation}
\label{eq: cyclic image pentagon}
(\alpha_{123}^{-1}\circ_2 1_\mu) \cdot(1_\mu\circ_2\alpha_{123}^{-1}) = (1_\mu\circ_1\alpha_{123}) \cdot(\alpha_{123}^{-1}\circ_1 1_\mu) \cdot(\alpha_{123}^{-1}\circ_3 1_\mu). 
\end{equation}

See Figure~\ref{fig:cyclic action on Pentagon}. It remains to see that this is a consequence of the original pentagon. Inverting \eqref{pent_rel_1_morph}, after specializing to the generators of $\PaRB$, gives
\[
(\alpha_{123}^{-1}\circ_1 1_\mu) \cdot(\alpha_{123}^{-1}\circ_3 1_\mu) = (1_\mu\circ_1\alpha_{123}^{-1}) \cdot(\alpha_{123}^{-1}\circ_2 1_\mu) \cdot(1_\mu\circ_2\alpha_{123}^{-1}).\] Postcomposing both sides with $1_\mu\circ_1\alpha_{123}$ yields
\[
(1_\mu\circ_1\alpha_{123})
\cdot(\alpha_{123}^{-1}\circ_1 1_\mu)
\cdot(\alpha_{123}^{-1}\circ_3 1_\mu)
=
(\alpha_{123}^{-1}\circ_2 1_\mu)
\cdot(1_\mu\circ_2\alpha_{123}^{-1}).
\]
This is precisely \eqref{eq: cyclic image pentagon}, with its two sides interchanged. Hence the cyclic action preserves the pentagon relation.

\smallskip
\noindent\emph{The twist relation.}
Recall that the twist relation in $\PaRB$ is
\begin{equation}
\label{eq: twist relation PaRB proof}
\twist\circ_1 1_\mu
=
\beta_{21}\cdot\beta_{12}\cdot
\bigl(1_\mu\circ(\twist,\twist)\bigr).
\end{equation}

We verify that this relation is preserved by $z_3^*$ (see Figure~\ref{fig:Cyclic action on twist_LHS} and Figure~\ref{fig:Cyclic action on twist_RHS}). The cyclic
compatibility identities give
$z_3^*(\twist\circ_1 1_\mu)=1_\mu\circ_2\twist$,
$z_3^*(1_\mu\circ_1\twist)=\twist\circ_1 1_\mu$, and
$z_3^*(1_\mu\circ_2\twist)=1_\mu\circ_1\twist$.
Moreover,
\[
z_3^*\beta_{12}
=
\bigl(1_{\mu_{21}}\circ_1\twist^{-1}\bigr)
\cdot\beta_{21}^{-1},
\qquad
z_3^*\beta_{21}
=
\bigl(1_\mu\circ_1\twist^{-1}\bigr)
\cdot\beta_{12}^{-1}.
\]
To obtain the second formula, note that
$\beta_{21}=s_1^*\beta_{12}$ and
$s_1z_3=(12)(012)=(02)=z_3^2s_1$. Hence
$z_3^*\beta_{21}=s_1^*z_3^{*2}\beta_{12}$, and the formula follows
by applying $z_3^*$ once more to the expression for
$z_3^*\beta_{12}$.

\begin{figure}[h]
    \centering
    \includegraphics[width=0.8\linewidth]{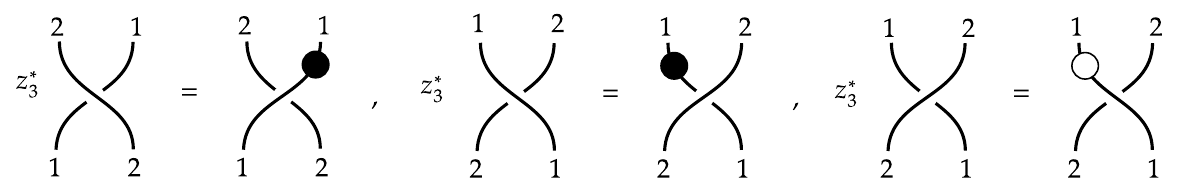}
\caption{The images under $z_3^*$ of the braids
$\beta_{12}^{-1}$, $\beta_{21}^{-1}$, and $\beta_{21}$, from left
to right.}    \label{fig:cyclic_action_on_braids}
\end{figure}

Applying $z_3^*$ to the right-hand side of
\eqref{eq: twist relation PaRB proof} gives
\begin{equation}
\label{eq: expanded cyclic image twist relation rhs}
\begin{split}
&\bigl(1_\mu\circ_1\twist^{-1}\bigr)
\cdot\beta_{12}^{-1}
\cdot\bigl(1_{\mu_{21}}\circ_1\twist^{-1}\bigr)
\cdot\beta_{21}^{-1} \cdot(\twist\circ_1 1_\mu)
\cdot(1_\mu\circ_1\twist).
\end{split}
\end{equation}
In the ribbon braid group, a twist may be moved along its labelled
ribbon through a crossing. In the present case this gives
\begin{equation}
\label{eq: sliding twist past inverse braiding}
\beta_{12}^{-1}
\cdot\bigl(1_{\mu_{21}}\circ_1\twist^{-1}\bigr)
=
\bigl(1_\mu\circ_1\twist^{-1}\bigr)
\cdot\beta_{12}^{-1}.
\end{equation}
Both twists in \eqref{eq: sliding twist past inverse braiding} lie on
the ribbon carrying boundary label $1$, which occupies the second
position at the source of $\beta_{12}^{-1}$ and the first position
at its target.

Rearranging \eqref{eq: twist relation PaRB proof} gives
\begin{equation}
\label{eq: rearranged twist relation PaRB}
\beta_{12}^{-1}\cdot\beta_{21}^{-1}
\cdot(\twist\circ_1 1_\mu)
=
\bigl(1_\mu\circ_1\twist\bigr)
\cdot\bigl(1_\mu\circ_2\twist\bigr).
\end{equation}
Using \eqref{eq: sliding twist past inverse braiding} and then
\eqref{eq: rearranged twist relation PaRB}, the expression in
\eqref{eq: expanded cyclic image twist relation rhs} becomes
\[
\begin{split}
&\bigl(1_\mu\circ_1\twist^{-1}\bigr)^2
\cdot\bigl(1_\mu\circ_1\twist\bigr)
\cdot\bigl(1_\mu\circ_2\twist\bigr)
\cdot\bigl(1_\mu\circ_1\twist\bigr) =1_\mu\circ_2\twist.
\end{split}
\]
Indeed, twists on distinct ribbons commute, so the two positive
twists on strand $1$ may be moved next to the two negative twists and
cancelled. The result is
$1_\mu\circ_2\twist
=z_3^*(\twist\circ_1 1_\mu)$, which is the cyclic image of the
left-hand side of \eqref{eq: twist relation PaRB proof}. Hence the
cyclic action preserves the twist relation.


\begin{figure}[htbp]
    \centering
    \begin{minipage}{0.48\textwidth}
        \centering
        \includegraphics[width=0.4\linewidth]{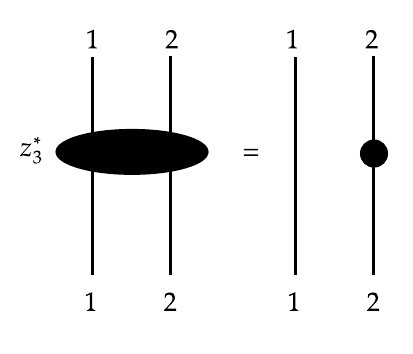}
\caption{The image under $z_3^*$ of the left-hand side of the twist
relation~\eqref{eq: twist relation PaRB proof}.}\label{fig:Cyclic action on twist_LHS}
    \end{minipage}\hfill
    \begin{minipage}{0.48\textwidth}
        \centering
        \includegraphics[width=0.7\linewidth]{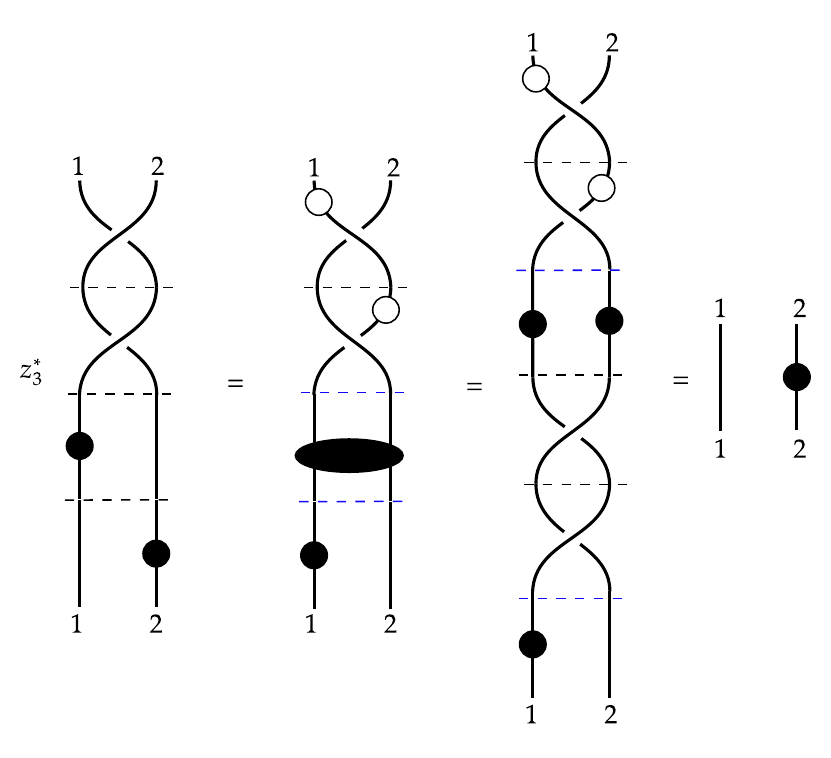}
       \caption{The image under $z_3^*$ of the right-hand side of the twist
relation~\eqref{eq: twist relation PaRB proof}.}
        \label{fig:Cyclic action on twist_RHS}
    \end{minipage}
\end{figure}

\smallskip
\noindent\emph{The second hexagon relation.}
We first verify directly that the cyclic action preserves the second
hexagon \eqref{hex_rel_2_morph}, which we write as an automorphism of
$\mu\circ_2\mu$ in $\PaRB(3)$:
\begin{equation}
\label{eq: second hexagon automorphism}
\alpha_{123} \cdot(\beta_{12}^{-1}\circ_1 1_\mu) \cdot\alpha_{312} \cdot(1_\mu\circ_1\beta_{12})_{132} \cdot\alpha_{132}^{-1} \cdot(1_\mu\circ_2\beta_{12}) =
1_{\mu\circ_2\mu}.
\end{equation}

We apply $z_4^*$ to each factor in this equation. The cyclic images of the associators are computed as before. For example, $\alpha_{123}$ is a morphism from $(12)3$ to $1(23)$, and the cyclic action exchanges these objects; the corresponding calculation gives
$z_4^*(\alpha_{123})=\alpha_{123}^{-1}$. Similarly, the cyclic action sends $(31)2$ to $(13)2$ and $3(12)$ to $1(32)$, and the calculation gives $z_4^*(\alpha_{312})=\alpha_{132}$. Finally, it sends $1(32)$ to $(21)3$ and $(13)2$ to $2(13)$, giving $z_4^*(\alpha_{132}^{-1})=\alpha_{213}$.

We next compute the cyclic images of the braiding factors (see Figure~\ref{fig:Cyclic action on H2}). Applying
the second case of \eqref{eq: cyclic compatibility PaRB}, with
$x=\beta_{12}^{-1}$, $y=1_\mu$, and $i=1$, gives
$z_4^*(\beta_{12}^{-1}\circ_1 1_\mu)
=z_3^*(1_\mu)\circ_2z_3^*(\beta_{12}^{-1})$.
Since $z_3^*1_\mu=1_\mu$ and the cyclic action preserves inverses,
we have
$z_3^*(\beta_{12}^{-1})
=(z_3^*\beta_{12})^{-1}
=\beta_{21}\cdot(1_{\mu_{21}}\circ_1\twist)$.
Consequently,
\begin{equation}
\label{eq: cyclic image inverse braiding first insertion}
z_4^*(\beta_{12}^{-1}\circ_1 1_\mu) = (1_\mu\circ_2\beta_{21}) \cdot \bigl(1_\mu\circ_2 (1_{\mu_{21}}\circ_1\twist)\bigr).
\end{equation}
This is a morphism from
$\mu\circ_2\mu_{21}=1(32)$ to $\mu\circ_2\mu=1(23)$, with a positive twist on the strand labelled $2$.

For the next factor, the case $i=2$ of
\eqref{eq: cyclic compatibility PaRB} gives
\begin{equation}
\label{eq: cyclic image second braiding factor}
z_4^*(1_\mu\circ_2\beta_{12}) = z_3^*(1_\mu)\circ_1\beta_{12} = 1_\mu\circ_1\beta_{12}.
\end{equation}
This is a morphism from $(12)3$ to $(21)3$, since the original morphism has type $1(23)\to1(32)$.

Finally, in the factor
$(1_\mu\circ_1\beta_{12})_{132}$, the subscript $132$ means that the
inputs have been relabelled by the permutation $(23)$; that is,
$(1_\mu\circ_1\beta_{12})_{132}
=(23)^*(1_\mu\circ_1\beta_{12})$.
Using the right-action convention and the identity
$(23)z_4=z_4(12)$, we obtain
\begin{equation}
\label{eq: cyclic image relabelled braiding factor}
z_4^*((1_\mu\circ_1\beta_{12})_{132}) =(12)^*\left[\bigl((1_{\mu_{21}}\circ_2 1_\mu) \circ_1\twist^{-1}\bigr) \cdot(\beta_{21}^{-1}\circ_2 1_\mu) \right].
\end{equation}
Here we have used the case $i=1$ of \eqref{eq: cyclic compatibility PaRB} and the formula
$z_3^*\beta_{12} =(1_{\mu_{21}}\circ_1\twist^{-1}) \cdot\beta_{21}^{-1}$. The morphism in \eqref{eq: cyclic image relabelled braiding factor} has type
$2(13)\to(13)2$, with a negative twist on the strand labelled $2$.

Applying $z_4^*$ to \eqref{eq: second hexagon automorphism}, the positive and negative twists on strand $2$ occur on opposite sides of $\alpha_{132}$. The associator changes only the parenthesisation and has trivial underlying ribbon braid. We may therefore move the twists along strand $2$ by ribbon-braid isotopy until they are adjacent, at which point they cancel. We are left with
\begin{equation}
\label{eq: cyclic image second hexagon automorphism}
\alpha_{123}^{-1} \cdot(1_\mu\circ_2\beta_{21}) \cdot\alpha_{132} \cdot(12)^*(\beta_{21}^{-1}\circ_2 1_\mu) \cdot\alpha_{213} \cdot(1_\mu\circ_1\beta_{12}) =
1_{\mu\circ_1\mu}.
\end{equation}
Thus $z_4^*$ carries the second hexagon to an automorphism of $\mu\circ_1\mu=(12)3$.

Rearranging \eqref{eq: cyclic image second hexagon automorphism}, it remains to prove
\begin{equation}
\label{eq: reduced cyclic second hexagon} 
\alpha_{132}\cdot(12)^*(\beta_{21}^{-1}\circ_2 1_\mu) \cdot\alpha_{213} = (1_\mu\circ_2\beta_{21}^{-1}) \cdot\alpha_{123} \cdot(1_\mu\circ_1\beta_{12}^{-1}).
\end{equation}
Both sides are morphisms from $(21)3$ to $1(32)$. Reading from right to left, the left-hand side is the composite
\begin{equation}
\label{eq: reduced cyclic second hexagon left string}
(21)3\xrightarrow{\alpha_{213}}2(13)\xrightarrow{(12)^*(\beta_{21}^{-1}\circ_2 1_\mu)}(13)2\xrightarrow{\alpha_{132}} 1(32),
\end{equation}
whereas the right-hand side is
\begin{equation}
\label{eq: reduced cyclic second hexagon right string}
(21)3\xrightarrow{1_\mu\circ_1\beta_{12}^{-1}}(12)3\xrightarrow{\alpha_{123}}1(23)\xrightarrow{1_\mu\circ_2\beta_{21}^{-1}}1(32).
\end{equation}

We now identify
\eqref{eq: reduced cyclic second hexagon} with a relabelled copy of
the original second hexagon. Applying the ordinary permutation
$(23)$ to \eqref{hex_rel_2_morph} gives
\begin{equation}
\label{eq: second hexagon relabelled 132}
(1_\mu\circ_1\beta_{12}) \cdot\alpha_{123}^{-1} \cdot(1_\mu\circ_2\beta_{12})_{132} =
\alpha_{213}^{-1} \cdot(\beta_{12}\circ_1 1_\mu)_{132} \cdot\alpha_{132}^{-1}.
\end{equation}
Both sides are morphisms from $1(32)$ to $(21)3$. More explicitly,
the left-hand side is
\begin{equation}
\label{eq: relabelled second hexagon left string}
1(32)\xrightarrow{(1_\mu\circ_2\beta_{12})_{132}}1(23)\xrightarrow{\alpha_{123}^{-1}}(12)3\xrightarrow{1_\mu\circ_1\beta_{12}}(21)3,
\end{equation}
while the right-hand side is
\begin{equation}
\label{eq: relabelled second hexagon right string}
1(32)\xrightarrow{\alpha_{132}^{-1}}(13)2\xrightarrow{(\beta_{12}\circ_1 1_\mu)_{132}}2(13)\xrightarrow{\alpha_{213}^{-1}}(21)3.
\end{equation}

Inverting \eqref{eq: second hexagon relabelled 132} reverses these
two strings of morphisms. The definitions of operadic insertion and
relabeling give
\begin{equation}
\label{eq: inverse relabelled braid identifications}
\bigl((\beta_{12}\circ_1 1_\mu)_{132}\bigr)^{-1} =
(12)^*(\beta_{21}^{-1}\circ_2 1_\mu),
\qquad
\bigl((1_\mu\circ_2\beta_{12})_{132}\bigr)^{-1} = 1_\mu\circ_2\beta_{21}^{-1}.
\end{equation}
Moreover, $(1_\mu\circ_1\beta_{12})^{-1} =1_\mu\circ_1\beta_{12}^{-1}$. After these identifications, the inverse of \eqref{eq: second hexagon relabelled 132} is precisely \eqref{eq: reduced cyclic second hexagon}. Hence the cyclic image of the second hexagon is the inverse of its ordinary $(23)$-relabeling, and therefore follows from the second hexagon itself.

\begin{figure}[htbp]
    \centering
    \begin{minipage}{0.48\textwidth}
        \centering
        \includegraphics[width=\linewidth]{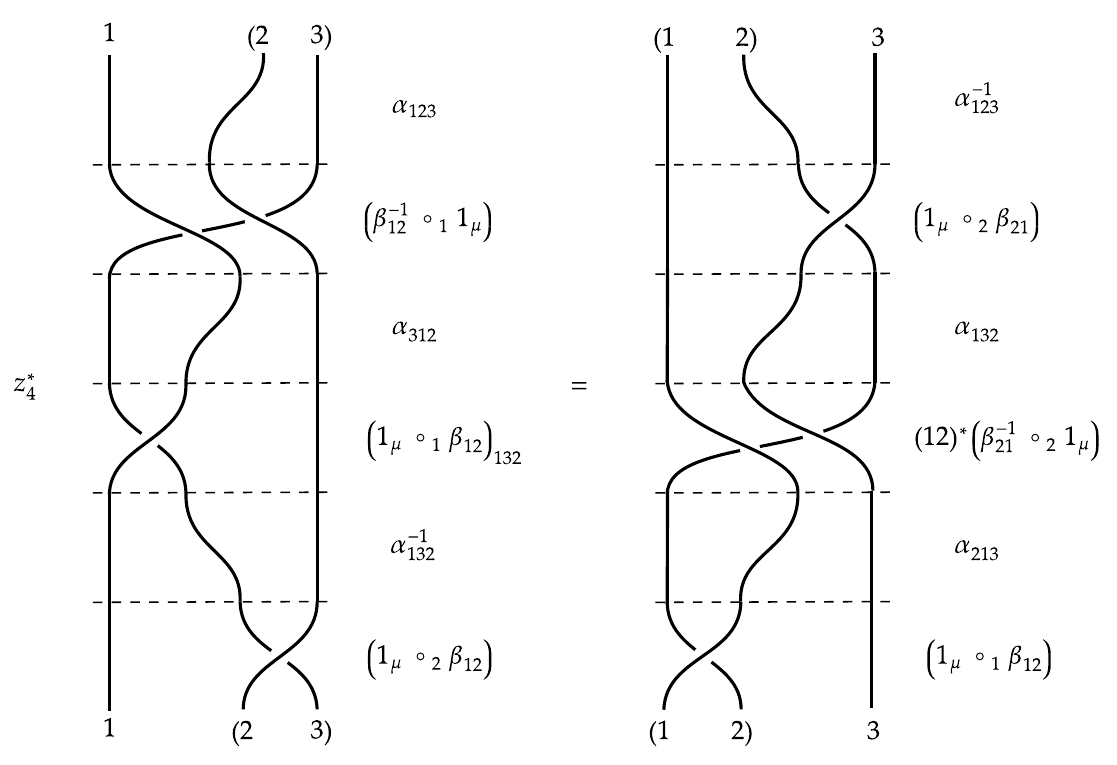}
\caption{The image under $z_4^*$ of the second hexagon
relation~\eqref{eq: cyclic image second hexagon automorphism}. The
right-hand composite is an automorphism of $(12)3$ and reduces to
$1_{\mu\circ_1\mu}\in\Aut_{\PaRB(3)}((12)3)$.}        \label{fig:Cyclic action on H2}
    \end{minipage}\hfill
    \begin{minipage}{0.48\textwidth}
        \centering
        \includegraphics[width=0.7\linewidth]{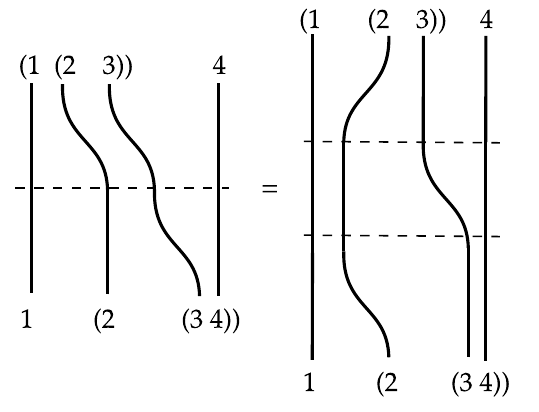}
        \caption{Cyclic action on the pentagon equation~\eqref{eq: cyclic image pentagon} of $\PaRB$}
        \label{fig:cyclic action on Pentagon}
    \end{minipage}
\end{figure}

\smallskip
\noindent\emph{The first hexagon relation.}
We now verify that the cyclic action preserves the first hexagon \eqref{hex_rel_1_morph}. We begin by writing it as an automorphism of
$\mu\circ_1\mu$ in $\PaRB(3)$:
\begin{equation}
\label{eq: first hexagon automorphism}
\alpha_{123}^{-1} \cdot(\beta_{12}^{-1}\circ_2 1_\mu) \cdot\alpha_{231}^{-1} \cdot(1_\mu\circ_2\beta_{12})_{213} \cdot\alpha_{213} \cdot(1_\mu\circ_1\beta_{12}) =
1_{\mu\circ_1\mu}.
\end{equation}

We again apply $z_4^*$ to each factor. The cyclic images of the associators, computed as before, are $z_4^*(\alpha_{123}^{-1})=\alpha_{123}$, $z_4^*(\alpha_{231}^{-1})=\alpha_{312}$, and $z_4^*(\alpha_{213})=\alpha_{231}$.

For the first braiding factor, the case $i=2$ of \eqref{eq: cyclic compatibility PaRB} gives
\[
z_4^*(\beta_{12}^{-1}\circ_2 1_\mu) = z_3^*(\beta_{12}^{-1})\circ_1 1_\mu.
\]
Since $z_3^*(\beta_{12}^{-1}) =\beta_{21}\cdot(1_{\mu_{21}}\circ_1\twist)$, it follows that
\begin{equation}\label{eq: cyclic image inverse braiding second insertion}
z_4^*(\beta_{12}^{-1}\circ_2 1_\mu) = (\beta_{21}\circ_1 1_\mu) \cdot \bigl(1_{\mu_{21}}\circ_1 (\twist\circ_1 1_\mu)\bigr).
\end{equation}
This is a morphism from $3(12)$ to $(12)3$. Here the twist is applied to a strand that has been doubled, so it will produce both individual
twists and a full two-strand twist.

For the rightmost braiding factor, the case $i=1$ of
\eqref{eq: cyclic compatibility PaRB} gives
\begin{equation}
\label{eq: cyclic image first insertion braiding H1}
z_4^*(1_\mu\circ_1\beta_{12})
=
\bigl((1_{\mu_{21}}\circ_2 1_\mu)
\circ_1\twist^{-1}\bigr)
\cdot(\beta_{21}^{-1}\circ_2 1_\mu).
\end{equation}
This is a morphism from $1(23)$ to $(23)1$. Reading from right to left, the cable formed by strands $2$ and $3$ first crosses over
strand $1$, after which a negative twist is applied to strand $1$.

It remains to compute the cyclic image of the relabelled factor $(1_\mu\circ_2\beta_{12})_{213}$. The subscript $213$ means relabelling by $(12)$. The right-action convention and the identity $(12)z_4=z_4^2(132)$ give
\[
z_4^*((1_\mu\circ_2\beta_{12})_{213}) = (132)^*z_4^{*2}(1_\mu\circ_2\beta_{12}).
\]
The first application of $z_4^*$ sends $1_\mu\circ_2\beta_{12}$ to $1_\mu\circ_1\beta_{12}$, so a second application and
\eqref{eq: cyclic image first insertion braiding H1} give
\begin{equation}
\label{eq: cyclic image relabelled braiding H1}
z_4^*((1_\mu\circ_2\beta_{12})_{213}) = (132)^*\left[ \bigl((1_{\mu_{21}}\circ_2 1_\mu) \circ_1\twist^{-1}\bigr) \cdot(\beta_{21}^{-1}\circ_2 1_\mu)
\right].
\end{equation}
This is a morphism from $2(31)$ to $(31)2$, with a negative twist on the strand labelled $2$.

We now expand the doubled positive twist in
\eqref{eq: cyclic image inverse braiding second insertion}. Using the
equivalent form of the ribbon-twist relation noted above, and writing
it in our categorical composition order, we have
\begin{equation}
\label{eq: positive ribbon twist expansion H1}
\twist\circ_1 1_\mu = \bigl(1_\mu\circ(\twist,\twist)\bigr) \cdot\beta_{21}\cdot\beta_{12}.
\end{equation}
Consequently,
\begin{equation}
\label{eq: inserted positive ribbon twist expansion H1}
1_{\mu_{21}}\circ_1(\twist\circ_1 1_\mu) = \bigl(1_{\mu_{21}}\circ_1 (1_\mu\circ(\twist,\twist))\bigr) \cdot(1_{\mu_{21}}\circ_1\beta_{21}) \cdot(1_{\mu_{21}}\circ_1\beta_{12}).
\end{equation}
The first factor on the right supplies positive twists on strands $1$ and $2$, while the remaining two factors form the full twist of these strands.

After substituting \eqref{eq: cyclic image inverse braiding second insertion},
\eqref{eq: cyclic image first insertion braiding H1},
\eqref{eq: cyclic image relabelled braiding H1}, and
\eqref{eq: inserted positive ribbon twist expansion H1} into the
cyclic image of \eqref{eq: first hexagon automorphism}, the positive twist on strand $1$ cancels the negative twist in \eqref{eq: cyclic image first insertion braiding H1}. Similarly, the positive twist on strand $2$ cancels the negative twist in
\eqref{eq: cyclic image relabelled braiding H1}. These cancellations are ribbon-braid isotopies along the corresponding labelled strands.
We are left with
\begin{equation}
\label{eq: reduced cyclic image first hexagon}
\alpha_{123}
\cdot(\beta_{21}\circ_1 1_\mu)
\cdot(1_{\mu_{21}}\circ_1\beta_{21})
\cdot(1_{\mu_{21}}\circ_1\beta_{12})
\cdot\alpha_{312}
\cdot(132)^*(\beta_{21}^{-1}\circ_2 1_\mu)
\cdot\alpha_{231}
\cdot(\beta_{21}^{-1}\circ_2 1_\mu) = 1_{\mu\circ_2\mu}. \end{equation} This is an automorphism of
$\mu\circ_2\mu=1(23)$.

To show that this equation holds in $\PaRB(3)$, it is easiest to identify the underlying braid of
\eqref{eq: reduced cyclic image first hexagon} using the isomorphism $\Aut_{\PaRB(3)}(\mu\circ_2\mu)\cong \PRB_3$ (see Figure~\ref{fig:Cyclic action on H1}). Recall that we write $\beta_i$ for the Artin generator in which the strand in position $i$ crosses over the strand in position $i+1$. The associators have trivial underlying braid, and the remaining factors contribute as follows. The two cabled inverse
braidings
\[ 
\beta_{21}^{-1}\circ_2 1_\mu
\qquad\text{and}\qquad
(132)^*(\beta_{21}^{-1}\circ_2 1_\mu)
\]
each contribute $\beta_2^{-1}\beta_1^{-1}$. The surviving full twist occurs on positions $2$ and $3$ and therefore contributes $\beta_2^2$. Finally,
$\beta_{21}\circ_1 1_\mu$ contributes $\beta_2\beta_1$. Hence the underlying braid of
\eqref{eq: reduced cyclic image first hexagon} is
\begin{equation}
\label{eq: braid word cyclic first hexagon}
(\beta_2\beta_1)
\cdot\beta_2^2
\cdot(\beta_2^{-1}\beta_1^{-1})
\cdot(\beta_2^{-1}\beta_1^{-1}).
\end{equation}
After one cancellation, this becomes
$\beta_2\beta_1\beta_2
\beta_1^{-1}\beta_2^{-1}\beta_1^{-1}$.
Using the braid relation
$\beta_2\beta_1\beta_2=\beta_1\beta_2\beta_1$, we obtain
\begin{equation}
\label{eq: braid reduction cyclic first hexagon}
\beta_2\beta_1\beta_2\beta_1^{-1}\beta_2^{-1}\beta_1^{-1} =
\beta_1\beta_2\beta_1 \beta_1^{-1}\beta_2^{-1}\beta_1^{-1} = 1.
\end{equation}
Geometrically, this is a Reidemeister III move followed by Reidemeister II cancellations. The framing twists have already
cancelled, so the complete ribbon braid in \eqref{eq: reduced cyclic image first hexagon} is the identity.
Thus the cyclic action preserves the first hexagon relation.

\begin{figure}[htbp]
    \centering
    \includegraphics[width=0.8\linewidth]{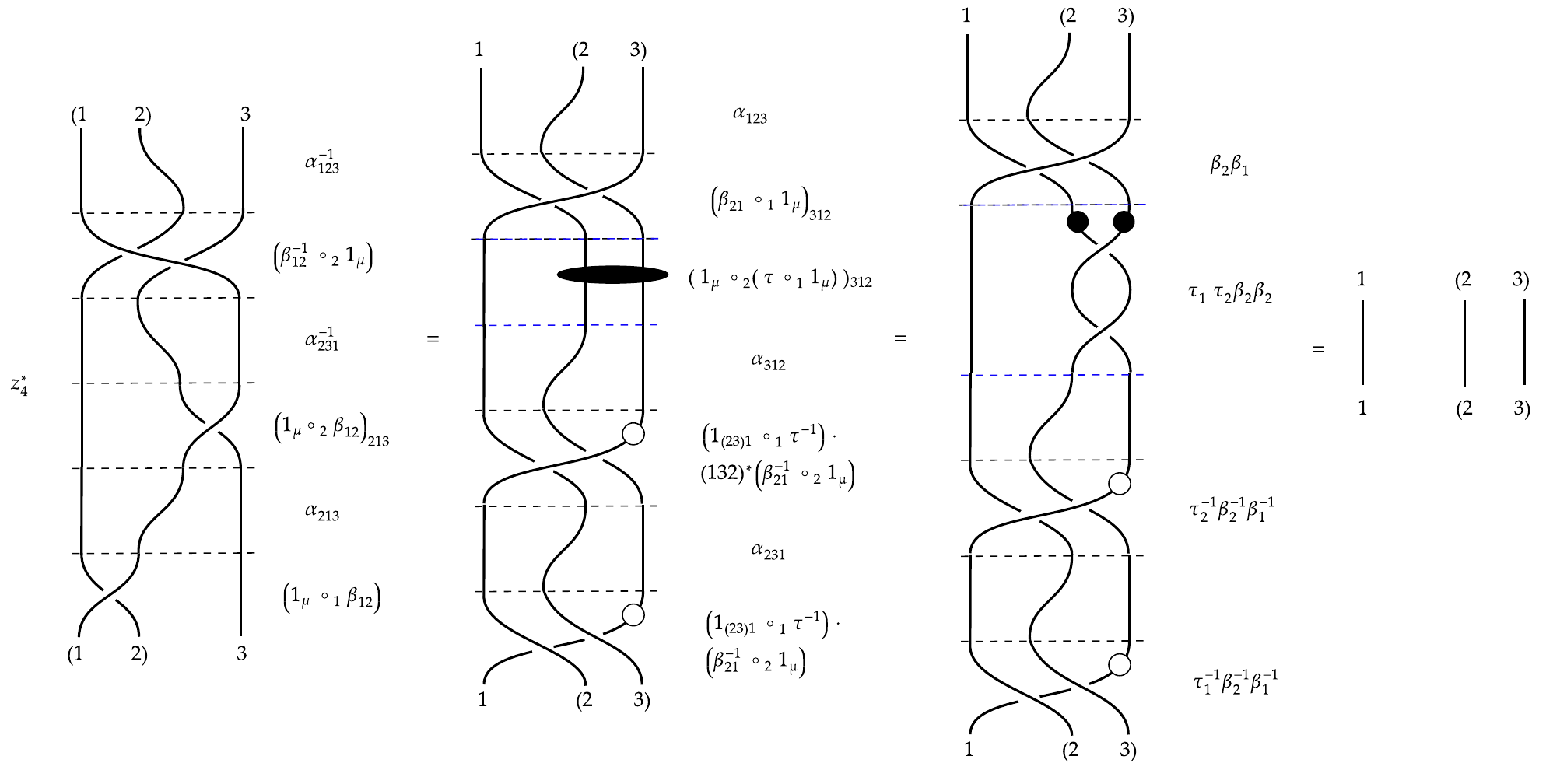}
\caption{The image under $z_4^*$ of the first hexagon
relation~\eqref{eq: first hexagon automorphism}. The second
braid is an automorphism of $1(23)$, and the third is evaluated in
$\Br_3$. The rightmost morphism is
$1_{\mu\circ_2\mu}\in\Aut_{\PaRB(3)}(1(23))$.}    \label{fig:Cyclic action on H1}
\end{figure}

We have now verified that the cyclic action preserves the twist,
pentagon, and hexagon relations in the presentation of $\PaRB$.
Together with
Lemma~\ref{lemma: extended symmetric action on generators}, this shows
that the assignments in
\eqref{eq: cyclic action PaRB generators} extend uniquely to a cyclic
operad structure on $\PaRB$. This completes the proof.
\end{proof}

\subsection{A cofibrant cyclic operad of parenthesised ribbon braids}
\label{subsec: cofibrant cyclic PaRB}

The cyclic operad $\PaRB^{\mathrm{cyc}}$ constructed above is not
cofibrant in $\Cyc(\Grpd)$. Indeed,
$\ob(\PaRB^{\mathrm{cyc}})$ is isomorphic to
$F_{\mathrm{cyc}}(A_2^+\backslash\Sigma_2^+)$, where
$A_2^+=\langle(012)\rangle$. The quotient
$\pi\colon\Sigma_2^+\to A_2^+\backslash\Sigma_2^+$ induces an
entrywise surjective morphism of cyclic operads in sets
\[
F_{\mathrm{cyc}}(\pi)\colon
F_{\mathrm{cyc}}(\Sigma_2^+)
\longrightarrow
F_{\mathrm{cyc}}(A_2^+\backslash\Sigma_2^+),
\]
where both generating extended symmetric sequences are concentrated
in arity two. If $\PaRB^{\mathrm{cyc}}$ were cofibrant, the lifting
criterion of Proposition~\ref{prop:cofibrations-cyc} would give a
section of $F_{\mathrm{cyc}}(\pi)$, and hence a
$\Sigma_2^+$-equivariant section of~$\pi$ in arity two.
No such section exists. The coset $A_2^+$ has stabilizer $A_2^+$,
whereas the regular action of $\Sigma_2^+$ on itself has trivial
stabilizers. Hence $\PaRB^{\mathrm{cyc}}$ is not cofibrant.

We now construct an explicit cofibrant replacement. Its objects are
pairs $(T,r)$, where $T$ is an object of
$\PaRB^{\mathrm{cyc}}$ and $r=(r_v)_{v\in V(T)}$ specifies a
distinguished incident half-edge at each vertex of $T$. These choices
are independent and need not point toward a common boundary edge.

Recall that the forgetful functor
$u^*\colon\Cyc(\Grpd)\to\Op(\Grpd)$ has a left adjoint
$u_!\colon\Op(\Grpd)\to\Cyc(\Grpd)$, called the \emph{cyclic
envelope}. We use the same notation for the corresponding adjunction
in sets. Let $M$ be the symmetric sequence concentrated in arity two,
with $M(2)=\Sigma_2$ carrying its regular right action. Since
$\magma=\ob(\PaRB)\cong F_{\mathrm{op}}(M)$, the free-operad and
cyclic-envelope adjunctions give
\[
u_!\magma
\cong
F_{\mathrm{cyc}}
\bigl(\Sigma_2\times_{\Sigma_2}\Sigma_2^+\bigr)
\cong
F_{\mathrm{cyc}}(\Sigma_2^+).
\]
The first isomorphism follows formally by applying the two
adjunctions to maps into an arbitrary cyclic operad; see also
\cite[Lemma~3.1 and Section~9.11]{Ward}.
Thus $u_!\magma$ is freely generated as a cyclic operad by the
regular right $\Sigma_2^+$-set in arity two.

More explicitly, an object of $(u_!\magma)(n)$ is represented by a
pair $(T,r)$, where $T$ is a planar unrooted tree with boundary
labelling
$\tau_T\colon\{0,\ldots,n\}\xrightarrow{\cong}\partial T$
and $r_v$ is a distinguished incident half-edge at each vertex $v$.
For $n\geq2$, the tree $T$ is trivalent. Its cyclic ordering at each
vertex, together with the choice of $r_v$, determines a linear
ordering of the three incident half-edges with $r_v$ first. The six
linear orderings account for the regular $\Sigma_2^+$-set of binary
generators. Operadic composition is given by grafting and retaining
the distinguished half-edge at each vertex; a distinguished half-edge
may become part of an internal edge. The extended symmetric-group
actions relabel the boundary. In arity one, we include the exceptional
one-edge tree, which has a unique decoration since it has no vertices.
There are no objects in arity zero.

The identification $u^*\PaRB^{\mathrm{cyc}}\cong\PaRB$ determines,
by adjunction, a counit
$\epsilon\colon u_!\PaRB\to\PaRB^{\mathrm{cyc}}$. On objects,
it is the morphism
$\epsilon_0\colon u_!\magma\to\ob(\PaRB^{\mathrm{cyc}})$ given by
$\epsilon_0(T,r)=T$. It forgets the distinguished half-edges while
retaining the planar structure and boundary labelling.

In arity two, $\epsilon_0$ is the quotient
$\Sigma_2^+\to A_2^+\backslash\Sigma_2^+$. It identifies the three
linear orderings related by cyclic rotation. The two elements of the
quotient are represented by the corollas
$\mu(1,2;0)$ and $\mu(2,1;0)$.

\begin{definition}
\label{def: cofibrant cyclic PaRB}
For each $n\geq0$, define
$\PaRB_{\mathrm{cof}}^{\mathrm{cyc}}(n)$ to be the groupoid obtained
by pulling back the morphisms of $\PaRB^{\mathrm{cyc}}(n)$ along
$\epsilon_{0,n}\colon(u_!\magma)(n)\to
\ob(\PaRB^{\mathrm{cyc}}(n))$.

Its objects are the pairs $(T,r)$ described above, so that
$\ob(\PaRB_{\mathrm{cof}}^{\mathrm{cyc}})=u_!\magma$. For two such
objects, set
\[
\Hom_{\PaRB_{\mathrm{cof}}^{\mathrm{cyc}}(n)}
\bigl((T_1,r_1),(T_2,r_2)\bigr)
=
\Hom_{\PaRB^{\mathrm{cyc}}(n)}(T_1,T_2).
\]
Categorical composition is inherited from
$\PaRB^{\mathrm{cyc}}(n)$. Because $\epsilon_0$ is a morphism of
cyclic operads in sets, operadic composition and the extended
symmetric-group actions on $u_!\magma$ and
$\PaRB^{\mathrm{cyc}}$ induce corresponding operations on these
groupoids. Thus $\PaRB_{\mathrm{cof}}^{\mathrm{cyc}}$ is a cyclic
operad in groupoids.
\end{definition}

There is a morphism
$q\colon\PaRB_{\mathrm{cof}}^{\mathrm{cyc}}\to
\PaRB^{\mathrm{cyc}}$ given by $q(T,r)=T$ on objects and by the
identity under the defining identification of morphism sets. The
counit of the cyclic-envelope adjunction factors as
\[
u_!\PaRB\longrightarrow
\PaRB_{\mathrm{cof}}^{\mathrm{cyc}}
\xrightarrow{\,q\,}
\PaRB^{\mathrm{cyc}}.
\]
The first morphism is the identity on objects and sends each morphism
of $u_!\PaRB$ to its image under the counit.

\begin{prop}
\label{prop: cofibrant cyclic PaRB replacement}
The cyclic operad $\PaRB_{\mathrm{cof}}^{\mathrm{cyc}}$ is
cofibrant, and
$q\colon\PaRB_{\mathrm{cof}}^{\mathrm{cyc}}\to
\PaRB^{\mathrm{cyc}}$ is an entrywise trivial fibration of
groupoids. Hence $q$ is a cofibrant replacement of
$\PaRB^{\mathrm{cyc}}$.
\end{prop}

\begin{proof}
The object operad
$\ob(\PaRB_{\mathrm{cof}}^{\mathrm{cyc}})=u_!\magma$
is freely generated by the regular right $\Sigma_2^+$-set in arity
two. Its generating extended symmetric sequence therefore has free
extended symmetric-group actions.
Proposition~\ref{prop:cofibrations-cyc} and
Corollary~\ref{cor:free-cyc-cof} imply that
$\PaRB_{\mathrm{cof}}^{\mathrm{cyc}}$ is cofibrant.

For each $n\geq1$, the object map $\epsilon_{0,n}$ has an aritywise
section obtained by choosing, at every vertex of $T$, the incident
half-edge pointing along the unique path toward the boundary edge
$\tau_T(0)$. For the exceptional one-edge tree, no choices are
required. Hence $q_n$ is surjective on objects; this also holds in
arity zero, where both groupoids are empty. The section does not
respect the extended symmetric-group actions, since it distinguishes
the boundary edge labelled~$0$.

By definition, $q_n$ induces a bijection on the morphism sets between
every pair of objects. It is therefore fully faithful and surjective
on objects. These are precisely the trivial fibrations in the
standard model structure on groupoids; see
\cite[Section~4]{Horel_profinite_groupoids}. Hence every $q_n$ is a
trivial fibration.

Weak equivalences and fibrations in $\Cyc(\Grpd)$ are defined
entrywise, so $q$ is a trivial fibration of cyclic operads. Its source
is cofibrant, and hence $q$ is a cofibrant replacement of
$\PaRB^{\mathrm{cyc}}$.
\end{proof}
\subsection{Profinite completion}
\label{subsec: Profinite completion}

For a small category $\bC$, its pro-category $\Pro(\bC)$ is the free
completion of $\bC$ under small cofiltered limits. Its objects are
small cofiltered diagrams in $\bC$. If
$X=\{X_i\}_{i\in I}$ and $Y=\{Y_j\}_{j\in J}$ are pro-objects, then
\[
\Hom_{\Pro(\bC)}(X,Y)
=
\lim_{j\in J}\operatorname*{colim}_{i\in I}
\Hom_{\bC}(X_i,Y_j).
\]

Let $\bC$ be an essentially small full subcategory of a category
$\bE$ admitting small cofiltered limits. The inclusion
$\bC\hookrightarrow\bE$ extends to a realization functor
$|-|\colon\Pro(\bC)\to\bE$, which sends a pro-object
$\{X_i\}_{i\in I}$ to $\lim_iX_i$. When this functor has a left
adjoint, we write the resulting adjunction as
\begin{equation}
\label{profinite completion adjunction}
\widehat{(-)}\colon\bE
\rightleftarrows
\Pro(\bC):|-|
\end{equation}
and call $\widehat{(-)}$ the \emph{pro-$\bC$ completion} functor.

Suppose that $\bE$ is a simplicial model category in which every
object is cofibrant, and that $\bC\subseteq\bE$ is an essentially
small full subcategory closed under finite limits and cotensors by
finite simplicial sets. Under the hypotheses of
\cite[Theorem~1.1]{BM20}, the category $\Pro(\bC)$ admits a
fibrantly generated simplicial model structure for which
\eqref{profinite completion adjunction} is a simplicial Quillen
adjunction. We use this construction below for profinite groupoids
and profinite spaces.

\begin{example}
\label{ex:profinite-sets}
Let $\mathbf{Fin}\subseteq\mathbf{Set}$ be the full subcategory of
finite sets. The category of \emph{profinite sets} is
$\widehat{\mathbf{Set}}:=\Pro(\mathbf{Fin})$. Taking inverse limits
identifies it with the category of compact, Hausdorff, totally
disconnected spaces and continuous maps. These are also called
\emph{Stone spaces}; see \cite[Chapter~VI]{Johnstone_stone_spaces}.

The category $\widehat{\mathbf{Set}}$ is complete and cocomplete, and
its cartesian product makes it symmetric monoidal. One may therefore
consider categories enriched in profinite sets. Such a category has
a set of objects, a profinite hom-set $\bD(x,y)$ for every pair of
objects, and continuous composition maps
$\bD(y,z)\times\bD(x,y)\to\bD(x,z)$.
\end{example}

\begin{example}
Let $\Grp_{\mathrm{fin}}\subseteq\Grp$ be the full subcategory of
finite groups. The category of \emph{profinite groups} is
$\widehat{\Grp}:=\Pro(\Grp_{\mathrm{fin}})$. It is equivalent to the
category of compact, Hausdorff, totally disconnected topological
groups and continuous homomorphisms.

The realization functor sends a profinite group to its underlying
abstract group. Its left adjoint is the classical profinite
completion functor. For a discrete group $G$, this completion is
\[
\widehat G
=
\lim_{N\triangleleft_fG}G/N,
\]
where the limit runs over the finite-index normal subgroups of~$G$.
\end{example}

\begin{example}
\label{example: etale fundamental groups}
Let $X$ be a connected, locally Noetherian scheme with geometric
point $\bar x\to X$. Its étale fundamental group
$\pi_1^{\mathrm{et}}(X,\bar x)$ is the automorphism group of the
geometric fibre functor from finite étale covers of $X$ to finite
sets. Equivalently, it is the inverse limit of the groups of deck
transformations of a cofinal system of pointed, connected finite
Galois covers.

For $n\geq2$, let
$\mathrm{Conf}_n=\mathbb A_{\mathbb Q}^n
\setminus\bigcup_{i<j}\{z_i=z_j\}$ be the scheme of ordered
configurations. Its complex points form a $K(\PB_n,1)$; see
\cite{FadellNeuwirth}.

Choose an embedding
$\overline{\mathbb Q}\hookrightarrow\mathbb C$ and compatible
geometric and topological basepoints $\bar p$ and $p$. The comparison
between étale and topological fundamental groups gives
\[
\pi_1^{\mathrm{et}}
\bigl(
\mathrm{Conf}_n\otimes_{\mathbb Q}\overline{\mathbb Q},
\bar p
\bigr)
\cong
\widehat{\pi_1(\mathrm{Conf}_n(\mathbb C),p)}
\cong
\widehat{\PB}_n.
\]
See \cite[Exposé~XII, Corollaire~5.2]{SGA1} for this comparison and
\cite[Remark~1.14]{furusho_galois_2017} for its relation to
profinite braid groups and the Grothendieck--Teichmüller action.
\end{example}

Profinite completion also extends from groups to groupoids. A
groupoid is \emph{finite} if it has finitely many morphisms, which
also implies that it has finitely many objects. Let
$\Grpd_{\mathrm{fin}}\subseteq\Grpd$ be the full subcategory of
finite groupoids. The category of \emph{profinite groupoids} is
$\widehat{\Grpd}:=\Pro(\Grpd_{\mathrm{fin}})$.

The completion of a discrete groupoid has a concrete description. A
congruence relation $R$ on a groupoid $\bE$ is \emph{cofinite} if the
quotient groupoid $\bE/R$ is finite. These relations form a
cofiltered poset, and the profinite completion of $\bE$ is represented
by
\[
\widehat{\bE}
=
\{\bE/R\}_{R\text{ cofinite}}.
\]
If $\bE$ is the one-object groupoid associated to a group $G$, this
recovers the one-object profinite groupoid associated to $\widehat G$;
see \cite[Definitions~4.18--4.20]{Horel_profinite_groupoids}.

The category $\widehat{\Grpd}$ carries the cocombinatorial model
structure of
\cite[Theorem~4.12]{Horel_profinite_groupoids}. This model structure
is fibrantly generated, and the adjunction
\[
\widehat{(-)}\colon\Grpd
\rightleftarrows
\widehat{\Grpd}:|-|
\]
is a Quillen adjunction with profinite completion as its left Quillen
functor; see
\cite[Proposition~4.22]{Horel_profinite_groupoids}.

Although profinite completion does not preserve limits in general, it
does preserve products of groupoids with finite object sets.

\begin{prop}[
\textnormal{\cite[Proposition~4.23]{Horel_profinite_groupoids}}]
\label{prop: profinite completion preserves products of groupoids}
Let $\bD$ and $\bE$ be groupoids with finite sets of objects. Then the
canonical map
\[
\widehat{\bD\times\bE}
\longrightarrow
\widehat{\bD}\times\widehat{\bE}
\]
is an isomorphism of profinite groupoids.
\end{prop}

The groupoids themselves need not be finite. In particular, their
automorphism groups may be infinite. Together with preservation of
the terminal groupoid, the proposition says that profinite completion
is strong symmetric monoidal on the full subcategory of groupoids
with finitely many objects.

Let $\calO$ be a cyclic operad in groupoids such that every
$\calO(n)$ has finitely many objects. Its \emph{entrywise profinite
completion} is defined by
$\widehat{\calO}(n)=\widehat{\calO(n)}$. Using the inverse of the
product isomorphism above, its partial compositions are the
composites
\[
\widehat{\calO}(n)\times\widehat{\calO}(m)
\cong
\widehat{\calO(n)\times\calO(m)}
\longrightarrow
\widehat{\calO}(n+m-1),
\]
where the final arrow is the completion of the corresponding
composition in $\calO$. The units and extended symmetric-group
actions are obtained in the same way. Naturality of the product
isomorphisms ensures that the cyclic-operad axioms are preserved.
Thus $\widehat{\calO}$ is a cyclic operad in profinite groupoids.

\subsubsection{Profinite completion in spaces}

A space $X$ is \emph{$\pi$-finite} if it has finitely many connected
components, each homotopy group $\pi_k(X,x)$ is finite, and there
exists an integer $N$ such that $\pi_k(X,x)=0$ for every $k>N$ and
every choice of basepoint~$x$. Let $\mathcal S_\pi$ denote the full
subcategory of the $\infty$-category of spaces spanned by the
$\pi$-finite spaces. This subcategory is closed under finite limits,
and the $\infty$-category of \emph{profinite spaces} is
$\Pro(\mathcal S_\pi)$.

For point-set arguments, we use Quick's model category of simplicial
profinite sets. Its underlying category is
\[
\widehat{\sSet}
:=
\operatorname{Fun}
(\Delta^{\mathrm{op}},\widehat{\mathbf{Set}}).
\]
Quick constructs a fibrantly generated model structure on
$\widehat{\sSet}$ in which the cofibrations are the monomorphisms and
the weak equivalences are the maps inducing isomorphisms on $\pi_0$,
on profinite fundamental groups, and on continuous cohomology with
finite local coefficients; see \cite[Theorem~2.12]{Quick}. This model
category presents the $\infty$-category $\Pro(\mathcal S_\pi)$; see
\cite[Theorem~1.1(4) and Corollary~6.6]{BM20}. The profinite
completion and underlying simplicial-set functors form a simplicial
Quillen adjunction
\begin{equation}
\label{eq: profinite completion of spaces}
\begin{tikzcd}
\widehat{(-)}\colon\sSet
  \arrow[r, shift left=1]
&
  \arrow[l, shift left=1]
\widehat{\sSet}\colon|-|.
\end{tikzcd}
\end{equation}
See \cite[Proposition~2.28]{Quick}.

The fundamental-groupoid and classifying-space adjunction also admits
a profinite analogue:
\[
\begin{tikzcd}
\widehat{\Pi}_1\colon\widehat{\sSet}
  \arrow[r, shift left=1.2]
&
  \arrow[l, shift left=1.2]
\widehat{\Grpd}\colon\widehat{\cs}.
\end{tikzcd}
\]
The classifying space of a finite groupoid is $\pi$-finite, so the
ordinary adjunction extends to pro-objects. With the model structures
above, the resulting adjunction is simplicial Quillen; see
\cite[Example~7.14]{BM20}. Together with its discrete counterpart, it
fits into a square of simplicial Quillen adjunctions that commutes up
to natural isomorphism:
\[
\begin{tikzcd}
\Grpd
  \arrow[d, shift right=1, "\widehat{(-)}"']
  \arrow[r, shift right=1, "\cs"']
&
\sSet
  \arrow[l, shift right=1, "\Pi_1"']
  \arrow[d, shift left=1, "\widehat{(-)}"]
\\
\widehat{\Grpd}
  \arrow[u, shift right=1, "|-|"']
  \arrow[r, shift right=1, "\widehat{\cs}"']
&
\widehat{\sSet}
  \arrow[l, shift right=1, "\widehat{\Pi}_1"']
  \arrow[u, shift left=1, "|-|"].
\end{tikzcd}
\]

For a general groupoid $\bE$, the canonical comparison map
$(\cs\bE)^{\wedge}\to
\widehat{\cs}\bigl(\widehat{\bE}\bigr)$ need not be a weak
equivalence in $\widehat{\sSet}$. A useful sufficient condition is
expressed in terms of Serre's notion of a good group.

\begin{definition}
\label{def: good groupoids}
A discrete group $\mathsf G$ is \emph{good} if, for every finite
discrete $\widehat{\mathsf G}$-module $\mathsf M$, the canonical map
\[
H^k_{\mathrm{cts}}(\widehat{\mathsf G},\mathsf M)
\longrightarrow
H^k(\mathsf G,\mathsf M)
\]
is an isomorphism for every $k\geq0$.

A groupoid $\bE$ with finitely many objects is \emph{good} if every
automorphism group $\Aut_{\bE}(x)$ is good.
\end{definition}

By \cite[Proposition~5.9]{Horel_profinite_groupoids}, if $\bE$ is a
good groupoid with finitely many objects, then this comparison map is
a weak equivalence in $\widehat{\sSet}$.

\begin{example}
\label{example: profinite framed configuration spaces}
For $n\geq1$, let
$\Conf_n=\{(z_1,\ldots,z_n)\in\mathbb A^n_{\mathbb Q}
\mid z_i\neq z_j\text{ for }i\neq j\}$ be the ordered configuration
scheme over $\mathbb Q$. Its \emph{framed configuration scheme} is
$\Conf_n^{\mathrm{fr}}=\Conf_n\times\mathbb G_m^n$, where the
$i$th $\mathbb G_m$-coordinate may be interpreted as a nonzero
tangent vector at $z_i$. Since $\mathbb C^\times\simeq S^1$, its
space of complex points is homotopy equivalent to
$\Conf_n(\mathbb C)\times(S^1)^n$, the configuration space of $n$
framed points in the plane. This is a $K(\PRB_n,1)$, where
$\PRB_n\cong\PB_n\times\mathbb Z^n$.

The framed pure braid group also has a mapping-class-group
interpretation. If $\Gamma_{0,n+1}$ denotes the mapping class group
of an oriented genus-zero surface with $n+1$ boundary components,
with diffeomorphisms and isotopies fixing the boundary pointwise,
then $\Gamma_{0,n+1}\cong\PRB_n$; see
\cite[Chapter~1]{Wahl_Thesis}.

The étale comparison theorem
\cite{Artin_Mazur,Friedlander} identifies the geometric étale
homotopy type of $\Conf_n^{\mathrm{fr}}$ with the profinite
completion of $\cs\PRB_n$. Since $\PRB_n$ is good, the comparison
criterion in
\cite[Proposition~5.9]{Horel_profinite_groupoids} identifies this
completion with the profinite classifying space of
$\widehat{\PRB}_n$. Consequently,
\[
\Pi_\infty^{\mathrm{et}}
\bigl(\Conf_n^{\mathrm{fr}}
\otimes_{\mathbb Q}\overline{\mathbb Q}\bigr)
\simeq
\widehat{\cs}\bigl(\widehat{\PRB}_n\bigr)
\cong
\widehat{\cs}\bigl(\widehat{\Gamma}_{0,n+1}\bigr).
\]
Thus the geometric étale homotopy type of the framed configuration
scheme is modelled by the classifying space of the profinite
completion of the corresponding genus-zero mapping class group.
\end{example}

\begin{example}
\label{example: Oda's theorem}
Let $n\geq2$ and let $\Gamma_0^{n+1}$ denote the pure mapping class
group of a sphere with $n+1$ marked points. In contrast,
$\Gamma_{0,n+1}$ denotes the mapping class group of a sphere with
$n+1$ boundary components, fixing the boundary pointwise.

The complex moduli space $\mathcal M_{0,n+1}(\mathbb C)$ is an
Eilenberg--Mac Lane space with
$\pi_1(\mathcal M_{0,n+1}(\mathbb C))\cong\Gamma_0^{n+1}$.
Oda's comparison theorem identifies its geometric étale homotopy type
with the profinite completion of this classifying space:
\[
\Pi_\infty^{\mathrm{et}}
\bigl(
\mathcal M_{0,n+1}
\otimes_{\mathbb Q}\overline{\mathbb Q}
\bigr)
\simeq
\widehat{\cs\Gamma_0^{n+1}}.
\]
See \cite[Theorem~1]{Oda}, as well as
\cite{Artin_Mazur,Friedlander} for the construction of étale homotopy
types.
\end{example}

Unlike the exact product-preservation result for groupoids in
Proposition~\ref{prop: profinite completion preserves products of groupoids},
profinite completion of spaces does not preserve finite products in
general. The following gives a useful sufficient condition for
product preservation up to weak equivalence.

\begin{prop}[{\cite[Proposition~3.9]{Boavida-Horel-Robertson}}]
\label{prop: profinite spaces products}
Let $X$ and $Y$ be connected simplicial sets whose homotopy groups are
good. Then the canonical map
\[
\widehat{X\times Y}
\longrightarrow
\widehat X\times\widehat Y
\]
is a weak equivalence of profinite spaces.
\end{prop}


We next explain why profinite completion of dendroidal objects may be
computed entrywise. Let $\mathscr K$ be a small category with finite
hom-sets and let $\bC$ be a small category with finite limits. A
functor $A\colon\mathscr K^{\mathrm{op}}\to\bC$ is
\emph{finitely coskeletal} if it is the right Kan extension of its
restriction to a finite full subcategory of
$\mathscr K^{\mathrm{op}}$. We denote the full subcategory of such
functors by
\[
\Fun_{\mathrm{fcosk}}(\mathscr K^{\mathrm{op}},\bC)
\subseteq
\Fun(\mathscr K^{\mathrm{op}},\bC).
\]
The dual of \cite[Theorem~2.3]{BM20} gives an equivalence
\[
\Pro\bigl(
\Fun_{\mathrm{fcosk}}(\mathscr K^{\mathrm{op}},\bC)
\bigr)
\simeq
\Fun\bigl(\mathscr K^{\mathrm{op}},\Pro(\bC)\bigr).
\]

Both $\Omega$ and $\Omega^{\mathrm{cyc}}$ are small categories with
finite hom-sets. For either category, let $\mathscr K_{\leq d}$ be
the full subcategory spanned by trees with at most $d$ edges. Each
$\mathscr K_{\leq d}$ is finite, and these subcategories exhaust
$\mathscr K$.

Restriction to $\mathscr K_{\leq d}^{\mathrm{op}}$ has a right
adjoint, denoted $\operatorname{cosk}_d$. For every dendroidal object
$X$, the canonical map
$X\to\lim_d\operatorname{cosk}_dX$ is an isomorphism. Indeed, if $T$
has at most $d$ edges, then
$X_T\to(\operatorname{cosk}_kX)_T$ is an isomorphism whenever
$k\geq d$. Since each $\operatorname{cosk}_dX$ is finitely
coskeletal, the preceding equivalence shows that pro-$\bC$
completion of dendroidal objects is computed entrywise.

We apply this with $\bE=\sSet$ or $\Grpd$. In the groupoid case,
$\bC$ is the category of finite groupoids. In the spatial case,
$\bC$ is the category of lean simplicial sets, meaning degreewise
finite simplicial sets that are $n$-coskeletal for some $n\geq0$;
see \cite[Section~2]{BM20}. The resulting model structure on
$\Pro(\bC)$ agrees with Quick's model structure on profinite spaces
by \cite[Example~5.5 and Corollary~6.6]{BM20}. Consequently, if
$X\colon(\Omega^{\mathrm{cyc}})^{\mathrm{op}}\to\bE$ is a cyclic
dendroidal object, its profinite completion is represented by the
cyclic dendroidal object $\widehat X$ defined by
$(\widehat X)_T=\widehat{X_T}$.

Entrywise completion therefore represents the pro-$\bC$ completion
of the dendroidal diagram. It need not preserve the Segal condition,
since this also requires completion to preserve the finite products
occurring in the Segal maps.

\begin{lemma}
\label{lem: profinite completion of cyclic infinity operads}
Let
$X\colon(\Omega^{\mathrm{cyc}})^{\mathrm{op}}\to\Grpd$
be a cyclic $\infty$-operad such that $X_{C_{n+1}}$ has finitely many
objects for every corolla $C_{n+1}$. Then its entrywise profinite
completion
$\widehat X\colon(\Omega^{\mathrm{cyc}})^{\mathrm{op}}
\to\widehat{\Grpd}$
is again a cyclic $\infty$-operad.
\end{lemma}

\begin{proof}
Since $X_\eta\simeq *$ and profinite completion preserves weak
equivalences between groupoids, we have $\widehat X_\eta\simeq *$.

For every unrooted tree $T$, the Segal map
\[
X_T\longrightarrow
\prod_{v\in V(T)}X_{C_{|\mathrm{nb}(v)|}}
\]
is a weak equivalence. Every groupoid is cofibrant and profinite
completion is a left Quillen functor, so its profinite completion is
again a weak equivalence. Moreover, since the groupoids
$X_{C_{|\mathrm{nb}(v)|}}$ have finitely many objects,
Proposition~\ref{prop: profinite completion preserves products of groupoids}
gives an isomorphism
\[
\widehat{\prod_{v\in V(T)}X_{C_{|\mathrm{nb}(v)|}}}
\cong
\prod_{v\in V(T)}
\widehat X_{C_{|\mathrm{nb}(v)|}}.
\]
It follows that the Segal map for $\widehat X$ is a weak equivalence,
so $\widehat X$ is a cyclic $\infty$-operad.
\end{proof}

The following lemma shows that the cyclic dendroidal nerve commutes
with entrywise profinite completion.

\begin{lemma}
\label{lemma: nerve commutes with profinite completion}
Let $\calO$ be a cyclic operad in groupoids such that each
$\calO(n)$ has finitely many objects. Then there is a natural
isomorphism of cyclic dendroidal objects
\[
\widehat{\nerve\calO}
\cong
\nerve(\widehat{\calO}).
\]
In particular, both sides are cyclic $\infty$-operads.
\end{lemma}

\begin{proof}
For every unrooted tree $T$, the cyclic dendroidal nerve is given by
\[
(\nerve\calO)_T
\cong
\prod_{v\in V(T)}\calO(\operatorname{nb}(v)).
\]
Since each factor has finitely many objects,
Proposition~\ref{prop: profinite completion preserves products of groupoids}
gives a natural isomorphism
\[
(\widehat{\nerve\calO})_T
=
\widehat{\prod_{v\in V(T)}
\calO(\operatorname{nb}(v))}
\cong
\prod_{v\in V(T)}
\widehat{\calO}(\operatorname{nb}(v))
=
(\nerve\widehat{\calO})_T.
\]
These isomorphisms are natural in $T$ and are compatible with
relabeling and grafting. They therefore assemble into the required
isomorphism of cyclic dendroidal objects.
\end{proof}

We can similarly describe the profinite completion of sufficiently
well-behaved cyclic $\infty$-operads in spaces. A complete
classification would require a cyclic analogue of the theory of
profinite $\infty$-operads developed in \cite{BM22}.

\begin{lemma}
\label{lemma: profinite completion of cyclic infinity operads in spaces}
Let
$X\colon(\Omega^{\mathrm{cyc}})^{\mathrm{op}}\to\sSet$
be a cyclic $\infty$-operad such that every
$X_{C_{n+1}}$ is connected and has good homotopy groups. Then its
entrywise profinite completion
\[
\widehat X\colon
(\Omega^{\mathrm{cyc}})^{\mathrm{op}}
\longrightarrow
\widehat{\sSet}
\]
is again a cyclic $\infty$-operad in profinite spaces.
\end{lemma}

\begin{proof}
For every unrooted tree $T$, the Segal map
\[
X_T\longrightarrow
\prod_{v\in V(T)}X_{C_{|\operatorname{nb}(v)|}}
\]
is a weak equivalence. Every simplicial set is cofibrant and
profinite completion is a left Quillen functor, so completing this
map again gives a weak equivalence. By
Proposition~\ref{prop: profinite spaces products} and induction on the
finite set $V(T)$, there is also a weak equivalence
\[
\widehat{\prod_{v\in V(T)}
X_{C_{|\operatorname{nb}(v)|}}}
\longrightarrow
\prod_{v\in V(T)}
\widehat X_{C_{|\operatorname{nb}(v)|}}.
\]
Their composite is the Segal map for $\widehat X$. Finally,
$\widehat X_\eta\simeq *$, since $X_\eta\simeq *$ and profinite
completion preserves weak equivalences between cofibrant objects.
Thus $\widehat X$ is a cyclic $\infty$-operad.
\end{proof}

\section{The Grothendieck--Teichmüller group
\texorpdfstring{$\GT$}{GT}}
\label{sec: GT}

The goal of this section is to identify $\GT$ with the group of
homotopy automorphisms of the profinite completion of the cyclic
operad of parenthesised ribbon braids. We proceed in two steps. In
Section~\ref{subsec: maps out of cyclic PaRB}, we identify the
profinite Grothendieck--Teichmüller monoid with the monoid of
object-fixing endomorphisms of
$\widehat{\PaRB}^{\mathrm{cyc}}$. In
Section~\ref{subsec: homotopy automorphisms of PaRB}, we pass from
these strict endomorphisms to the derived endomorphisms of the
associated profinite cyclic $\infty$-operad.

We first recall the profinite Grothendieck--Teichmüller monoid. Let
$\widehat{\Free}_2$ denote the profinite completion of the free group
on generators $x$ and $y$. Given $f\in\widehat{\Free}_2$, a
profinite group $G$, and elements $a,b\in G$, we write $f(a,b)$ for
the image of $f$ under the unique continuous homomorphism
$\widehat{\Free}_2\to G$ sending $x$ to $a$ and $y$ to $b$.

\begin{definition}[{\cite{Drin}}]
\label{defn: profinite GT}
The \emph{profinite Grothendieck--Teichmüller monoid}
$\underline{\GT}$ consists of pairs
$(\lambda,f)\in
(1+2\widehat{\mathbb Z})\times\widehat{\Free}_2$
satisfying the following relations:
\begin{enumerate}
\item[\textnormal{(I)}]
$f(x,y)^{-1}=f(y,x)$;

\item[\textnormal{(II)}]
$x^\nu f(x,y)y^\nu f(y,z)z^\nu f(z,x)=1$
in $\widehat{\Free}_2$, where $xyz=1$ and $\nu =\frac{(\lambda-1)}{2}\in\widehat{\mathbb Z}$;

\item[\textnormal{(III)}]
the pentagon relation
\[f(x_{12},x_{23})
\cdot f(x_{12}x_{13},x_{24}x_{34})
\cdot f(x_{23},x_{34})
=
f(x_{13}x_{23},x_{34})
\cdot f(x_{12},x_{23}x_{24})\]
holds in $\widehat{\PB}_4$.
\end{enumerate}
\end{definition}

The group of invertible elements $\GT:=\underline{\GT}^{\times}$ is the
\emph{profinite Grothendieck--Teichmüller group}. Every $(\lambda,f)\in\GT$ satisfies $\lambda\in\widehat{\mathbb Z}^{\times}$, since the projection $(\lambda,f)\mapsto\lambda$ is a monoid homomorphism. Multiplication in $\underline{\GT}$ is given by 
\begin{equation}\label{eqn:GT multiplication}
 (\lambda_1,f_1)*(\lambda_2,f_2)
=
\left(
\lambda_1\lambda_2,\,
f_1\bigl(
f_2(x,y)^{-1}x^{\lambda_2}f_2(x,y),
y^{\lambda_2}
\bigr)\cdot f_2(x,y)
\right).   
\end{equation}

\begin{remark}\label{rem:GT-convention}
Relations (II) and (III) are Drinfeld’s, but in reverse order. We compose morphisms
right to left, $q\cdot p$ meaning ``first $p$'', while products in $\PB_n$ are read left to
right, so that $g(x_{12},x_{23})$ multiplies the associator on the right in
\eqref{eqn:GT-action on ribbons}. 
Our multiplication convention agrees  with the endomorphism of
$\widehat{\Free}_2$ associated to $(\lambda,f)$, which sends
$x$ to $f(x,y)^{-1}x^\lambda f(x,y)$ and $y$ to $y^\lambda$.
\end{remark}

\subsection{Maps out of the cyclic operad of parenthesised ribbon
braids}
\label{subsec: maps out of cyclic PaRB}

Let
$\widehat{(-)}\colon
\Grpd\rightleftarrows\widehat{\Grpd}:|-|$
denote the profinite-completion adjunction
\eqref{profinite completion adjunction}. Since completion is applied
aritywise and preserves the finite object sets of $\PaRB$, forgetting
the cyclic structure gives a natural identification
$u^*\widehat{\PaRB}^{\mathrm{cyc}}\cong\widehat{\PaRB}$.

It follows that every morphism of cyclic operads
$f\colon\PaRB^{\mathrm{cyc}}\to
\lvert\widehat{\PaRB}^{\mathrm{cyc}}\rvert$
has an underlying operad morphism
$u^*(f)\colon\PaRB\to\lvert\widehat{\PaRB}\rvert$. We therefore
obtain a natural map
\[
u^*\colon
\Hom_{\Cyc(\Grpd)}
\bigl(
\PaRB^{\mathrm{cyc}},
\lvert\widehat{\PaRB}^{\mathrm{cyc}}\rvert
\bigr)
\longrightarrow
\Hom_{\Op(\Grpd)}
\bigl(
\PaRB,
\lvert\widehat{\PaRB}\rvert
\bigr).
\]

We will prove that this map is a bijection after restricting to
object-fixing morphisms. Equivalently, every object-fixing morphism
$\PaRB\to\lvert\widehat{\PaRB}\rvert$ extends uniquely to an
object-fixing morphism of cyclic operads
$\PaRB^{\mathrm{cyc}}\to
\lvert\widehat{\PaRB}^{\mathrm{cyc}}\rvert$. In particular, the
action of the profinite Grothendieck--Teichmüller monoid on
$\widehat{\PaRB}$ preserves the cyclic structure constructed in
Proposition~\ref{prop: cyclic structure on PaRB}.

We write $\underline{\Hom}_{\Op(\Grpd)}(\calO,\calP)$ and
$\underline{\Hom}_{\Cyc(\Grpd)}(\calO,\calP)$ for the subsets of operad and
cyclic-operad morphisms, respectively, that induce the identity on
the common operad of objects.

\begin{lemma}
\label{lemma: operad maps lift to cyclic operad maps of PaRB}
The forgetful map
\[
u^*\colon \underline{\Hom}_{\Cyc(\Grpd)} \bigl(\PaRB^{\mathrm{cyc}},
\lvert\widehat{\PaRB}^{\mathrm{cyc}}\rvert\bigr) \longrightarrow \underline{\Hom}_{\Op(\Grpd)} \bigl(\PaRB,\lvert\widehat{\PaRB}\rvert\bigr)
\]
is a bijection. 
\end{lemma}

\begin{proof}
Let $\bar f\colon\PaRB\to\lvert\widehat{\PaRB}\rvert$ be an object-fixing morphism of operads. By the profinite-completion adjunction and
\cite[Propositions~7.3 and~8.1]{Boavida-Horel-Robertson}, the morphism $\bar f$ is determined by a pair $(\lambda,g)\in\underline{\GT}$, where $\lambda=1+2\nu$ for some $\nu\in\widehat{\mathbb Z}$ and $g\in\widehat{\Free}_2$. Its values on the generators are
\begin{equation}\label{eqn:GT-action on ribbons}
   \bar f(\mu)=\mu,\qquad \bar f(\tau)=\tau^\lambda,\qquad \bar f(\beta_{12})=\beta_{12}^\lambda,\qquad \bar f(\alpha_{123}) = \alpha_{123}\cdot g(x_{12},x_{23}). 
\end{equation}
Here $\beta_{12}^{\lambda} =\beta_{12}\cdot(\beta_{21}\cdot\beta_{12})^\nu$. The composite $\beta_{21}\cdot\beta_{12}$ is an automorphism of $\mu$, so its $\nu$th profinite power is well defined.

The morphisms $\tau$, $\beta_{12}$, and $\alpha_{123}$ generate the morphisms of $\PaRB$ under operadic composition, categorical composition, and the ordinary symmetric-group actions. Since $\bar f$ fixes objects, it already commutes with the extended symmetric-group action on objects. It also commutes with the ordinary
symmetric-group actions on morphisms because it is a morphism of symmetric operads. The group $\Sigma_n^+$ is generated by $\Sigma_n$ together with $s_0=(01)$, so it is enough to verify $\bar f(s_0^*\gamma)=s_0^*\bar f(\gamma)$ for $\gamma=\tau,\beta_{12}$, and $\alpha_{123}$.

\smallskip
\noindent\emph{The twist.}
In arity one, $s_0=z_2$ and $s_0^*\tau=\tau$. Since $s_0^*$ acts continuously on $\lvert\widehat{\PaRB}^{\mathrm{cyc}}(1)\rvert$, it commutes with
profinite powers. Therefore $s_0^*\bar f(\tau) = s_0^*(\tau^\lambda) = (s_0^*\tau)^\lambda = \tau^\lambda = \bar f(s_0^*\tau).$

\smallskip
\noindent\emph{The braiding.}
Let $x_{12}=\beta_{21}\cdot\beta_{12} \in\Aut_{\PaRB(2)}(\mu)\cong \PRB_2$ denote the full twist on the two strands. Thus $\beta_{12}^\lambda=\beta_{12}\cdot x_{12}^\nu$. We also write $\tau_i=1_\mu\circ_i\tau$ for the twist on the $i$th strand.

Recall that $s_0^*\beta_{12}=\tau_2^{-1}\cdot\beta_{12}^{-1}$. The corresponding formula for $\beta_{21}=(12)^*\beta_{12}$ is
\[
s_0^*\beta_{21}= (1_{(12)^*\mu}\circ_2\tau^{-1}) \cdot\beta_{21}^{-1}.
\]
Indeed, since $z_3=s_0s_1$ and the action is on the right, one has $s_0^*=s_1^*z_3^*$. Applying this identity to $\beta_{21}$ and using the previously computed formula for $z_3^*\beta_{21}$ gives the displayed expression.

Functoriality gives
$s_0^*x_{12}
=s_0^*\beta_{21}\cdot s_0^*\beta_{12}$.
Substitution of the preceding formulas, followed by simplification
using the ribbon-braid relations, gives
\begin{equation}
\label{eq: s0 action conjugated full twist}
\beta_{12}^{-1}\cdot s_0^*(x_{12})\cdot\beta_{12}
=
\tau_2^{-2}\cdot x_{12}^{-1}.
\end{equation}

Using this identity and
$x_{12}\cdot\tau_2=\tau_2\cdot x_{12}$, we obtain
\[
\begin{aligned}
s_0^*\bar f(\beta_{12})
&=
\tau_2^{-1}\cdot\beta_{12}^{-1}
\cdot(s_0^*x_{12})^\nu\\
&=
\tau_2^{-1}\cdot
\bigl(
\beta_{12}^{-1}\cdot s_0^*(x_{12})\cdot\beta_{12}
\bigr)^\nu
\cdot\beta_{12}^{-1}\\
&=
\tau_2^{-\lambda}\cdot
(\beta_{12}^\lambda)^{-1}.
\end{aligned}
\]
On the other hand,
$\bar f(s_0^*\beta_{12})
=\bar f(\tau_2^{-1}\cdot\beta_{12}^{-1})
=\tau_2^{-\lambda}\cdot(\beta_{12}^\lambda)^{-1}$.
Thus
$\bar f(s_0^*\beta_{12})=s_0^*\bar f(\beta_{12})$.

\smallskip
\noindent\emph{The associator.}
It remains to prove
$\bar f(s_0^*\alpha_{123})=s_0^*\bar f(\alpha_{123})$.
Recall that
$s_0^*\alpha_{123}=(231)^*\alpha_{123}^{-1}$.

We use the natural identification
$\Aut_{\PaRB(3)}(\mu\circ_1\mu)\cong\PRB_3$. Under the decomposition
$\PRB_3\cong\PB_3\times\mathbb Z^3$, the subgroup $\PB_3$ consists
of the pure ribbon braids with trivial strand twists. Its standard
generators are
$x_{12}=\beta_1^2$, $x_{23}=\beta_2^2$, and
$x_{13}=\beta_2\cdot\beta_1^2\cdot\beta_2^{-1}$.
More explicitly,
\[
x_{12}
=
(1_\mu\circ_1\beta_{21})
\cdot(1_\mu\circ_1\beta_{12}),
\]
whereas
\[
x_{23}
=
\alpha_{123}^{-1}
\cdot(1_\mu\circ_2\beta_{21})
\cdot(1_\mu\circ_2\beta_{12})
\cdot\alpha_{123}.
\]
The element $x_{13}$ is the element of
$\PB_3\subseteq\PRB_3$ represented by the braid
$\beta_2\cdot\beta_1^2\cdot\beta_2^{-1}$.

We first compute $s_0^*x_{12}$. Since
$z_4=s_0s_1s_2$, we have $s_0=z_4s_2s_1$ and hence
$s_0^*=(231)^*z_4^*$ on $\PaRB(3)$. In particular,
$s_0^*(\mu\circ_1\mu)=(231)^*(\mu\circ_2\mu)$.
Functoriality and cyclic compatibility give
\[
s_0^*x_{12}
=
(231)^*\Bigl(
(z_3^*\beta_{21}\circ_2 1_\mu)
\cdot
(z_3^*\beta_{12}\circ_2 1_\mu)
\Bigr).
\]
The composite inside the parentheses is an automorphism of
$\mu\circ_2\mu$.

Recall that
$z_3^*\beta_{12}
=(1_{(12)^*\mu}\circ_1\tau^{-1})\cdot\beta_{21}^{-1}$
and
$z_3^*\beta_{21}
=(1_\mu\circ_1\tau^{-1})\cdot\beta_{12}^{-1}$.
Substitution shows that the underlying ribbon braid is
\[
\tau_1^{-1}\cdot\beta_1^{-1}\cdot\beta_2^{-2}
\cdot\beta_1^{-1}\cdot\tau_1^{-1}.
\]

Our standard braid words are written as automorphisms of
$\mu\circ_1\mu$. The corresponding automorphism of
$\mu\circ_2\mu$ is obtained by transport along $\alpha_{123}$.
Consequently,
\[
s_0^*x_{12}
=
(231)^*\left(
\alpha_{123}\cdot\tau_1^{-1}
\cdot\beta_1^{-1}\cdot\beta_2^{-2}\cdot\beta_1^{-1}
\cdot\tau_1^{-1}\cdot\alpha_{123}^{-1}
\right).
\]
The braid between the two twist factors is pure. Under
$\PRB_3\cong\PB_3\times\mathbb Z^3$, every pure braid commutes with
the individual strand twists. We may therefore combine the two
copies of $\tau_1^{-1}$, obtaining
\begin{equation}
\label{eq: s0 action on x12 braid word}
s_0^*x_{12}
=
(231)^*\left(
\alpha_{123}\cdot\tau_1^{-2}\cdot
(\beta_1\cdot\beta_2^2\cdot\beta_1)^{-1}
\cdot\alpha_{123}^{-1}
\right).
\end{equation}

Let
$w=x_{12}\cdot x_{13}\cdot x_{23}
=(\beta_1\cdot\beta_2)^3$
be the standard generator of the centre of $\PB_3$. The braid
relation gives
\[
w\cdot x_{23}^{-1}
=
x_{12}\cdot x_{13}
=
\beta_1\cdot\beta_2^2\cdot\beta_1.
\]
Thus
$(\beta_1\cdot\beta_2^2\cdot\beta_1)^{-1}
=x_{23}\cdot w^{-1}$.
Substituting this into
\eqref{eq: s0 action on x12 braid word} gives
\begin{equation}
\label{eq: s0 action on x12}
s_0^*x_{12}
=
(231)^*\left(
\alpha_{123}\cdot
(\tau_1^{-2}w^{-1})\cdot x_{23}
\cdot\alpha_{123}^{-1}
\right).
\end{equation}

We similarly apply $s_0^*$ to the displayed expression for $x_{23}$.
Cyclic compatibility carries the two middle factors to the
relabelled composite
\[
(1_\mu\circ_1\beta_{21})
\cdot(1_\mu\circ_1\beta_{12})
=
x_{12}.
\]
The cyclic images of the two associators transport this automorphism
from $\mu\circ_1\mu$ to $\mu\circ_2\mu$. Applying the relabelling
$(231)^*$ therefore gives
\begin{equation}
\label{eq: s0 action on x23}
s_0^*x_{23}
=
(231)^*\left(
\alpha_{123}\cdot x_{12}\cdot\alpha_{123}^{-1}
\right).
\end{equation}
The ribbon-braid identities used in deriving
\eqref{eq: s0 action on x12} and
\eqref{eq: s0 action on x23} remain valid after profinite completion.

We now compare the two sides of the required identity. Since
$\bar f(\alpha_{123})
=\alpha_{123}\cdot g(x_{12},x_{23})$, relation~\textnormal{(I)}
gives
\[
\bar f(s_0^*\alpha_{123})
=
(231)^*\left(
g(x_{23},x_{12})\cdot\alpha_{123}^{-1}
\right).
\]
On the other hand,
\[
s_0^*\bar f(\alpha_{123})
=
s_0^*\alpha_{123}\cdot
g(s_0^*x_{12},s_0^*x_{23}).
\]
Substituting
\eqref{eq: s0 action on x12} and
\eqref{eq: s0 action on x23}, and using the compatibility of
free-group word evaluation with simultaneous conjugation, gives
\[
s_0^*\bar f(\alpha_{123})
=
(231)^*\left(
g\bigl((\tau_1^{-2}w^{-1})\cdot x_{23},x_{12}\bigr)
\cdot\alpha_{123}^{-1}
\right).
\]

Relation~\textnormal{(III)} in the definition of $\underline{\GT}$ implies that $g$ belongs to the closed commutator subgroup of $\widehat{\Free}_2$; see the discussion preceding \cite[Proposition~7.3]{Boavida-Horel-Robertson}, and \cite[Proposition 1.5]{ihara}. Equivalently, the exponent sums of both free generators in $g$ vanish.

Set $r=\tau_1^{-2}w^{-1}$. The element $w$ is central in $\PB_3$,
and the individual strand twists commute with pure braids. Hence
$r$ commutes with both $x_{12}$ and $x_{23}$. Evaluating the
commutator word $g$ at $(r\cdot x_{23},x_{12})$ introduces the
central factor $r$ with total exponent zero. Consequently,
\[
g(r\cdot x_{23},x_{12})
=
g(x_{23},x_{12}).
\]
It follows that
$s_0^*\bar f(\alpha_{123})
=\bar f(s_0^*\alpha_{123})$.

We have shown that $\bar f$ commutes with $s_0^*$ on the generating
morphisms, and hence with all the extended symmetric-group actions.
It therefore defines a morphism of cyclic operads
$\PaRB^{\mathrm{cyc}}\to
\lvert\widehat{\PaRB}^{\mathrm{cyc}}\rvert$.
The refinement is unique because the forgetful functor
$u^*\colon\Cyc(\Grpd)\to\Op(\Grpd)$ is faithful. Thus the
forgetful map in the statement is a bijection.
\end{proof}

The following proposition combines Proposition~7.3 and Proposition~8.1 of 
\cite{Boavida-Horel-Robertson}.

\begin{prop}
\label{End of PaRB is profinite GT}
Let $\underline{\End}_{\Op}(-)$ denote the monoid of object-fixing
endomorphisms in $\Op(\widehat{\Grpd})$. There are natural
isomorphisms of profinite monoids
\[
\underline{\End}_{\Op}(\widehat{\PaB})
\cong
\underline{\End}_{\Op}(\widehat{\PaRB})
\cong
\underline{\GT}.
\]
\end{prop}

Explicitly, the object-fixing endomorphism of
$\widehat{\PaRB}$ corresponding to
$(\lambda,f)\in\underline{\GT}$ is determined on the generating
morphisms by
\[
\mu\longmapsto\mu,\qquad
\beta_{12}\longmapsto\beta_{12}^{\lambda},\qquad
\tau\longmapsto\tau^{\lambda},\qquad
\alpha_{123}\longmapsto
\alpha_{123}\cdot f(x_{12},x_{23}).
\]
Here $\lambda=1+2\nu$ for some $\nu\in\widehat{\mathbb Z}$ and $\beta_{12}^{\lambda} =\beta_{12}\cdot x_{12}^{\nu}$. Conversely, every pair satisfying the defining relations of $\underline{\GT}$ determines a unique object-fixing endomorphism. Such an endomorphism is invertible precisely when
the corresponding pair belongs to $\GT$.

We now extend this identification to the cyclic setting.
Lemma~\ref{lemma: operad maps lift to cyclic operad maps of PaRB}
shows that the cyclic structure on
$\widehat{\PaRB}^{\mathrm{cyc}}$ imposes no additional conditions
on object-fixing endomorphisms.

\begin{thm}
\label{thm: GT monoid iso to cyclic endomorphisms}
Let $\underline{\End}_{\Cyc}(-)$ denote the monoid of object-fixing
endomorphisms in $\Cyc(\widehat{\Grpd})$. There is a canonical
isomorphism of profinite monoids $\underline{\End}_{\Cyc}(\widehat{\PaRB}^{\mathrm{cyc}}) \cong
\underline{\GT}.$ 
\end{thm}

\begin{proof}
The profinite-completion adjunction gives a bijection
\[
\underline{\End}_{\Cyc}(\widehat{\PaRB}^{\mathrm{cyc}})\cong\underline{\Hom}_{\Cyc(\Grpd)}\left(\PaRB^{\mathrm{cyc}},\lvert\widehat{\PaRB}^{\mathrm{cyc}}\rvert\right).
\]
It sends a continuous endomorphism of $\widehat{\PaRB}^{\mathrm{cyc}}$ to its restriction along the completion map $\PaRB^{\mathrm{cyc}}\to \lvert\widehat{\PaRB}^{\mathrm{cyc}}\rvert$. Its inverse is given by continuous extension. Since completion preserves the finite object sets of $\PaRB$, this correspondence restricts to object-fixing morphisms.

By
Lemma~\ref{lemma: operad maps lift to cyclic operad maps of PaRB},
forgetting the cyclic structure induces a bijection
\[
\underline{\Hom}_{\Cyc(\Grpd)}
\left(
\PaRB^{\mathrm{cyc}},
\lvert\widehat{\PaRB}^{\mathrm{cyc}}\rvert
\right)
\cong
\underline{\Hom}_{\Op(\Grpd)}
\left(
\PaRB,
\lvert\widehat{\PaRB}\rvert
\right).
\]
A second application of the profinite-completion adjunction identifies the latter set with
$\underline{\End}_{\Op}(\widehat{\PaRB})$. Proposition~\ref{End of PaRB is profinite GT} therefore gives
\[
\underline{\End}_{\Cyc}(\widehat{\PaRB}^{\mathrm{cyc}})\cong \underline{\GT}.
\]

Restriction and continuous extension preserve composition, so this is an isomorphism of monoids. The explicit description by
generators identifies both sides as the same closed submonoid of a finite product of profinite morphism groups; hence the isomorphism and its inverse are continuous.
\end{proof}
\subsection{Homotopy automorphisms of
\texorpdfstring{$\widehat{\PaRB}^{\mathrm{cyc}}$}{PaRBcyc}}
\label{subsec: homotopy automorphisms of PaRB}

Recall from
Proposition~\ref{prop: cofibrant cyclic PaRB replacement} that
$\PaRB_{\mathrm{cof}}^{\mathrm{cyc}}$ is cofibrant in
$\Cyc(\Grpd)$ and that
\[
q\colon
\PaRB_{\mathrm{cof}}^{\mathrm{cyc}}
\longrightarrow
\PaRB^{\mathrm{cyc}}
\]
is an entrywise weak equivalence. Since every groupoid is cofibrant
and profinite completion is a left Quillen functor, completion
preserves weak equivalences between groupoids. Applying it in every
arity therefore gives an entrywise weak equivalence
\begin{equation}
\label{eq: completed cyclic PaRB resolution}
\widehat q\colon
\widehat{\PaRB}_{\mathrm{cof}}^{\mathrm{cyc}}
\overset{\sim}{\longrightarrow}
\widehat{\PaRB}^{\mathrm{cyc}}.
\end{equation}
We use this map to compare the strict endomorphisms studied above
with the derived endomorphisms obtained after applying the cyclic
dendroidal nerve.

Let $I$ be the groupoid with objects $0$ and $1$ and a unique
isomorphism between them, regarded as a constant profinite
groupoid. For a cyclic operad $\calQ$ in profinite groupoids, the
entrywise cotensor $\calQ^I$ is again a cyclic operad. Evaluation at
the two objects of $I$ gives morphisms
$d_0,d_1\colon\calQ^I\to\calQ$.

A \emph{homotopy} between cyclic-operad morphisms
$f_0,f_1\colon\calP\to\calQ$ is a morphism
$H\colon\calP\to\calQ^I$ satisfying $d_0H=f_0$ and $d_1H=f_1$.
Concretely, such a homotopy assigns to every object
$x\in\calP(n)$ a morphism
$\eta_x\colon f_0(x)\to f_1(x)$ in $\calQ(n)$ such that
\[
f_1(\gamma)\cdot\eta_x
=
\eta_y\cdot f_0(\gamma)
\]
for every morphism $\gamma\colon x\to y$. These morphisms must also
satisfy
\[
\eta_{x\circ_i y}
=
\eta_x\circ_i\eta_y,
\qquad
\eta_{\sigma^*x}
=
\sigma^*\eta_x
\quad\text{for every }\sigma\in\Sigma_n^+,
\]
and the component at the operadic unit must be the identity. This defines an equivalence relation on cyclic-operad morphisms. We
write $[\calP,\calQ]_{\Cyc(\widehat{\Grpd})}$ for the resulting set of homotopy classes.

Recall that an object of
$\PaRB_{\mathrm{cof}}^{\mathrm{cyc}}(n)$ is a pair $(T,r)$, where
$T$ is an object of $\PaRB^{\mathrm{cyc}}(n)$ and
$r=(r_v)_{v\in V(T)}$ chooses a half-edge incident to each
vertex $v$. These choices are independent and need not point
towards a common boundary edge. The morphism
$q\colon\PaRB_{\mathrm{cof}}^{\mathrm{cyc}}
\to\PaRB^{\mathrm{cyc}}$
forgets this decoration, so $q(T,r)=T$. On morphisms, $q$ is
the identity under the defining identification
\[
\Hom_{\PaRB_{\mathrm{cof}}^{\mathrm{cyc}}(n)}
\bigl((T,r),(S,s)\bigr)
=
\Hom_{\PaRB^{\mathrm{cyc}}(n)}(T,S).
\]

Choosing, at every vertex $v$, the half-edge $r_v$ pointing
towards the boundary edge $\tau_T(0)$ gives a distinguished
decoration of $T$. The full subgroupoids spanned by trees with
these decorations are preserved by the ordinary symmetric-group
actions and operadic compositions. They therefore assemble into
a copy of $\PaRB$ inside
$u^*\PaRB_{\mathrm{cof}}^{\mathrm{cyc}}$.
After profinite completion, this gives an inclusion
\[
j\colon
\widehat{\PaRB}
\longrightarrow
u^*\widehat{\PaRB}_{\mathrm{cof}}^{\mathrm{cyc}}.
\]
The composite
\[
\widehat{\PaRB}
\xrightarrow{\,j\,}
u^*\widehat{\PaRB}_{\mathrm{cof}}^{\mathrm{cyc}}
\xrightarrow{\,u^*\widehat q\,}
u^*\widehat{\PaRB}^{\mathrm{cyc}}
\cong
\widehat{\PaRB}
\]
is the identity.

\begin{prop}
\label{prop: rooted and cyclic homotopy classes PaRB}
Restriction along $j$ induces a bijection
\[
\left[
\widehat{\PaRB}_{\mathrm{cof}}^{\mathrm{cyc}},
\widehat{\PaRB}^{\mathrm{cyc}}
\right]_{\Cyc(\widehat{\Grpd})}
\longrightarrow
\left[
\widehat{\PaRB},
\widehat{\PaRB}
\right]_{\Op(\widehat{\Grpd})}.
\]
\end{prop}
\begin{proof}
Let $F\colon\widehat{\PaRB}_{\mathrm{cof}}^{\mathrm{cyc}}
\to\widehat{\PaRB}^{\mathrm{cyc}}$ be a cyclic-operad morphism,
and set $g:=u^*F\circ j$. Here $j$ selects, at every vertex, the
incident half-edge pointing toward the boundary edge labelled~$0$.
By \cite[Proposition~8.1]{Boavida-Horel-Robertson}, $g$ is homotopic
to a unique object-fixing morphism $g_0$. By
Lemma~\ref{lemma: operad maps lift to cyclic operad maps of PaRB}
and the universal property of profinite completion, $g_0$ has a
unique object-fixing cyclic refinement $g_0^{\mathrm{cyc}}$.
We construct a cyclic homotopy
$F\simeq g_0^{\mathrm{cyc}}\circ\widehat q$.

Write $j(\mu)=(\mu,r)$, where $r$ is the distinguished half-edge incident to the boundary edge labelled~$0$, and let $z_3=(012)$. For two decorations $r',r''$ of the same tree $T$, let $\iota_{r',r''}\colon(T,r')\to(T,r'')$ be the unique morphism sent by $\widehat q$ to $1_T$. Since $z_3^*\mu=\mu$, the morphism $\iota_{r,z_3^*r}$ has source $j(\mu)$ and target $z_3^*j(\mu)$. Full faithfulness of $\widehat q$ gives $\iota_{z_3^{2*}r,r}\cdot\iota_{z_3^*r,z_3^{2*}r} \cdot\iota_{r,z_3^*r}=1_{j(\mu)}$.

Choose an isomorphism $\rho\in\Hom_{\widehat{\PaRB}^{\mathrm{cyc}}(2)}(F(j(\mu)),\mu)$, which exists because the target groupoid is
connected. We first adjust $\rho$ so that $(z_3^*\rho)\cdot F(\iota_{r,z_3^*r})=\rho$. The discrepancy $d=(z_3^*\rho)\cdot F(\iota_{r,z_3^*r})\cdot\rho^{-1}$
belongs to
\[
\Aut_{\widehat{\PaRB}^{\mathrm{cyc}}(2)}(\mu) \cong\widehat{\PB}_2\times\widehat{\mathbb Z}^{\,2} \cong\widehat{\mathbb Z}^{\,3}.
\]
This group has a basis given by the boundary twists $t_0=\tau\circ_1 1_\mu$, $t_1=1_\mu\circ_1\tau$, and $t_2=1_\mu\circ_2\tau$: the latter two are the individual ribbon twists, while $t_0$ is their product with the full twist of the two strands. 

The chosen $\rho$ need not yet respect cyclic rotation. The two morphisms $\rho$ and $(z_3^*\rho)\cdot F(\iota_{r,z_3^*r})$ have the same source $F(j(\mu))$ and target $\mu$. Their discrepancy is the automorphism
\[
d=(z_3^*\rho)\cdot F(\iota_{r,z_3^*r})\cdot\rho^{-1}
\]
of $\mu$. We will modify $\rho$ so that $d$ becomes the identity. The cyclic axioms give $z_3^*t_0=t_2$, $z_3^*t_2=t_1$, and $z_3^*t_1=t_0$. Writing this group additively, the preceding identity for the comparison morphisms implies $(1+z_3^*+z_3^{2*})d=0$. Since $z_3^*$ cyclically permutes the
basis,
\[
\ker(1+z_3^*+z_3^{2*})=\operatorname{im}(1-z_3^*),
\]
both being the subgroup of triples whose coordinates sum to zero. Choose $b\in\widehat{\mathbb Z}^{\,3}$ with $d=b-z_3^*b$. Replacing $\rho$ by $b\cdot\rho$ replaces $d$ by $z_3^*b+d-b=0$, giving the required compatibility.

The object operad of
$\widehat{\PaRB}_{\mathrm{cof}}^{\mathrm{cyc}}$ is freely
generated by the regular $\Sigma_2^+$-set containing $j(\mu)$.
Thus $\rho$ extends uniquely, by the extended symmetric-group
actions and operadic composition, to a compatible family
$\eta_y\colon F(y)\to\widehat q(y)$. Define $H$ on objects by
$H(y)=\widehat q(y)$ and, for $\gamma\colon y\to y'$, set
\[
H(\gamma)=\eta_{y'}\cdot F(\gamma)\cdot\eta_y^{-1}.
\]
The compatibility of the family makes $H$ a cyclic-operad morphism
and gives a homotopy $F\simeq H$. Our choice of $\rho$ ensures that
$H(\iota_{r,z_3^*r})=1_\mu$.

For any two decorations $r',r''$ of a tree $T$, changing one
distinguished half-edge at a time expresses $\iota_{r',r''}$ as
a composite of elementary changes. Each is obtained from an
extended symmetric translate of $\iota_{r,z_3^*r}$ or its inverse
and identity morphisms at the remaining vertices by operadic
composition. These identities follow from full faithfulness of
$\widehat q$. Since $H$ sends the binary comparison to an identity,
it follows that $H(\iota_{r',r''})=1_T$. Hence $H$ descends uniquely
through $\widehat q$ to an object-fixing cyclic-operad morphism
$h\colon\widehat{\PaRB}^{\mathrm{cyc}}
\to\widehat{\PaRB}^{\mathrm{cyc}}$. The descended maps on morphism
sets are continuous because $\widehat q$ identifies the corresponding
profinite morphism sets.

Restricting $F\simeq h\circ\widehat q$ along $j$ gives
$g\simeq u^*h$. Uniqueness of the object-fixing representative
implies $u^*h=g_0$, and uniqueness of the cyclic refinement gives
$h=g_0^{\mathrm{cyc}}$. Therefore
\begin{equation}
\label{eq: cyclic homotopy to object fixing representative}
F\simeq g_0^{\mathrm{cyc}}\circ\widehat q.
\end{equation}

For surjectivity, given an operad morphism
$g\colon\widehat{\PaRB}\to\widehat{\PaRB}$, choose its unique
object-fixing representative $g_0$. Since
$\widehat q\circ j=\id_{\widehat{\PaRB}}$, the cyclic morphism
$g_0^{\mathrm{cyc}}\circ\widehat q$ restricts to $g_0$ and hence
maps to the homotopy class of $g$.

For injectivity, suppose that $u^*F_0\circ j$ and $u^*F_1\circ j$
are homotopic. By
\cite[Proposition~8.1]{Boavida-Horel-Robertson}, their common
homotopy class has a unique object-fixing representative $g_0$.
Applying \eqref{eq: cyclic homotopy to object fixing representative}
to $F_0$ and $F_1$ gives
$F_0\simeq g_0^{\mathrm{cyc}}\circ\widehat q\simeq F_1$.
Thus restriction along $j$ is injective and hence induces the
asserted bijection.
\end{proof}

We now compare these explicit homotopy classes with the derived
endomorphisms of the corresponding profinite cyclic
$\infty$-operad. The following is the cyclic analogue of
\cite[Proposition~8.2]{Boavida-Horel-Robertson}.

\begin{prop}
\label{prop: strict and derived cyclic endomorphisms PaRB}
There is a natural bijection
\[
\left[\widehat{\PaRB}_{\mathrm{cof}}^{\mathrm{cyc}},\widehat{\PaRB}^{\mathrm{cyc}}\right]_{\Cyc(\widehat{\Grpd})}\cong\pi_0\mathbb R\End_{\Cyc_\infty(\widehat{\Grpd})}\left(\nerve\widehat{\PaRB}^{\mathrm{cyc}}\right).
\]
\end{prop}

\begin{proof}
The profinite-completion adjunction identifies morphisms $\widehat{\PaRB}_{\mathrm{cof}}^{\mathrm{cyc}} \to\widehat{\PaRB}^{\mathrm{cyc}}$ of cyclic operads in profinite groupoids with morphisms $\PaRB_{\mathrm{cof}}^{\mathrm{cyc}}\to\lvert\widehat{\PaRB}^{\mathrm{cyc}}\rvert$ of cyclic operads in groupoids. Since $I$ is finite, the underlying-groupoid functor satisfies $\lvert\calQ^I\rvert\cong\lvert\calQ\rvert^I$. The adjunction therefore preserves $I$-homotopies and induces a natural bijection
\[
\left[\widehat{\PaRB}_{\mathrm{cof}}^{\mathrm{cyc}},\widehat{\PaRB}^{\mathrm{cyc}}\right]_{\Cyc(\widehat{\Grpd})}\cong\left[\PaRB_{\mathrm{cof}}^{\mathrm{cyc}},\left|\widehat{\PaRB}^{\mathrm{cyc}}\right|\right]_{\Cyc(\Grpd)}.
\]

The cyclic operad $\PaRB_{\mathrm{cof}}^{\mathrm{cyc}}$ is cofibrant in $\Cyc(\Grpd)$, while every cyclic operad in groupoids is fibrant. Hence
\[
\left[\PaRB_{\mathrm{cof}}^{\mathrm{cyc}},\left|\widehat{\PaRB}^{\mathrm{cyc}}\right|\right]_{\Cyc(\Grpd)}
\cong
\pi_0\mathbb R\Map_{\Cyc(\Grpd)}\left(\PaRB_{\mathrm{cof}}^{\mathrm{cyc}},\left|\widehat{\PaRB}^{\mathrm{cyc}}\right|\right).
\]

By the homotopical full faithfulness of the cyclic dendroidal nerve, Theorem~\ref{thm: dendroidal nerve is homotopically fully faithful}, and the derived profinite-completion adjunction for cyclic dendroidal objects, there are natural weak equivalences
\[
\begin{aligned}
\mathbb R\Map_{\Cyc(\Grpd)}\left(\PaRB_{\mathrm{cof}}^{\mathrm{cyc}},\left|\widehat{\PaRB}^{\mathrm{cyc}}\right|\right)
&\simeq
\mathbb R\Map_{\Cyc_\infty(\Grpd)}
\left(
\nerve\PaRB_{\mathrm{cof}}^{\mathrm{cyc}},
\nerve\left|\widehat{\PaRB}^{\mathrm{cyc}}\right|
\right)
\\
&\simeq
\mathbb R\Map_{\Cyc_\infty(\widehat{\Grpd})}
\left(
\widehat{\nerve\bigl(
\PaRB_{\mathrm{cof}}^{\mathrm{cyc}}
\bigr)},
\nerve\widehat{\PaRB}^{\mathrm{cyc}}
\right).
\end{aligned}
\]
For the second equivalence, we use the natural identification
\[
\nerve\left|\widehat{\PaRB}^{\mathrm{cyc}}\right|
\cong
\left|\nerve\widehat{\PaRB}^{\mathrm{cyc}}\right|.
\]
Indeed, the cyclic dendroidal nerve is evaluated on a tree by taking
a finite product of arity groupoids, and the underlying-groupoid
functor preserves finite products.

The groupoid
$\PaRB_{\mathrm{cof}}^{\mathrm{cyc}}(n)$ has finitely many objects
in every arity. Lemma~\ref{lemma: nerve commutes with profinite completion}
therefore gives a natural identification
\[
\widehat{\nerve\bigl(
\PaRB_{\mathrm{cof}}^{\mathrm{cyc}}
\bigr)}
\cong
\nerve\widehat{\PaRB}_{\mathrm{cof}}^{\mathrm{cyc}}.
\]
Moreover, the entrywise weak equivalence
\[
\widehat q\colon
\widehat{\PaRB}_{\mathrm{cof}}^{\mathrm{cyc}}
\longrightarrow
\widehat{\PaRB}^{\mathrm{cyc}}
\]
induces a weak equivalence after applying the cyclic dendroidal
nerve, since the nerve is computed on each tree by a finite product
of arity groupoids. Precomposition with $\nerve\widehat q$ therefore
identifies the components of the preceding mapping space with
\[
\pi_0\mathbb R\End_{\Cyc_\infty(\widehat{\Grpd})}
\left(
\nerve\widehat{\PaRB}^{\mathrm{cyc}}
\right),
\]
as required.
\end{proof}

\begin{thm}
\label{thm: GT as homotopy cyclic endomorphisms}
There is an isomorphism of monoids
\[
\pi_0\mathbb R\End_{\Cyc_\infty(\widehat{\Grpd})}
\bigl(\nerve\widehat{\PaRB}^{\mathrm{cyc}}\bigr)
\cong
\underline{\GT}.
\]
\end{thm}

\begin{proof}
The action of $\underline{\GT}$ on
$\widehat{\PaRB}^{\mathrm{cyc}}$ induces a monoid homomorphism
\[
\underline{\GT}
\longrightarrow
\pi_0\mathbb R\End_{\Cyc_\infty(\widehat{\Grpd})}
\bigl(\nerve\widehat{\PaRB}^{\mathrm{cyc}}\bigr).
\]
Propositions~\ref{prop: strict and derived cyclic endomorphisms PaRB}
and~\ref{prop: rooted and cyclic homotopy classes PaRB}, together
with \cite[Proposition~8.1]{Boavida-Horel-Robertson}, show that its
underlying map of sets identifies each homotopy class with its unique
object-fixing representative. By
Proposition~\ref{End of PaRB is profinite GT}, these
representatives are parametrised by $\underline{\GT}$. The displayed
monoid homomorphism is therefore bijective.
\end{proof}

\begin{cor}
\label{cor: GT as cyclic operad}
There is an isomorphism of groups
\[
\pi_0\mathbb R\Aut_{\Cyc_\infty(\widehat{\Grpd})}
\bigl(\nerve\widehat{\PaRB}^{\mathrm{cyc}}\bigr)
\cong
\GT.
\]
\end{cor}

\begin{proof}
This follows by passing to the invertible elements in
Theorem~\ref{thm: GT as homotopy cyclic endomorphisms}.
\end{proof}

In the next section, we transport this action to a profinite cyclic
$\infty$-operad modelling framed configuration spaces.

\subsection{Galois actions on a cyclic
\texorpdfstring{$\infty$}{infinity}-operad of configuration spaces}
\label{subsec: Galois action on cyclic infinity conf space}

We now transport the action of
Theorem~\ref{thm: GT as homotopy cyclic endomorphisms} from the
parenthesised ribbon-braid model to the cyclic $\infty$-operad of
framed configuration spaces constructed in
Theorem~\ref{thm: cyclic infinity operad of framed configurations}.

\begin{lemma}
\label{lem: csNPaRB is profinite cyclic infty operad}
The entrywise profinite completion
$\bigl(\cs\nerve\PaRB^{\mathrm{cyc}}\bigr)^{\wedge}$ is a cyclic
$\infty$-operad in profinite spaces.
\end{lemma}

\begin{proof}
In arity zero, $\PaRB^{\mathrm{cyc}}(0)$ is empty. For every
$n\geq1$, the groupoid $\PaRB^{\mathrm{cyc}}(n)$ is connected, has
finitely many objects, and has automorphism groups isomorphic to
$\PRB_n\cong\PB_n\times\mathbb Z^n$. The groups $\PB_n$ are good by
\cite[Corollary~5.11]{Horel_profinite_groupoids}. Since $\PB_n$ and
$\mathbb Z^n$ are finitely presented and good,
\cite[Corollary~5.12]{Horel_profinite_groupoids} implies that
$\PRB_n$ is good. The result now follows from
Lemma~\ref{lemma: profinite completion of cyclic infinity operads in spaces}.
\end{proof}

There are two ways to obtain a profinite cyclic dendroidal object
from $\PaRB^{\mathrm{cyc}}$: one may apply the cyclic dendroidal
nerve and classifying-space functor before completion, or complete
the cyclic operad first. These constructions agree up to weak
equivalence.

\begin{lemma}
\label{lem:comparison cyclic PaRB nerves}
There is a weak equivalence of cyclic $\infty$-operads in profinite
spaces
\[
\bigl(\cs\nerve\PaRB^{\mathrm{cyc}}\bigr)^{\wedge}
\overset{\sim}{\longrightarrow}
\widehat{\cs}\nerve\widehat{\PaRB}^{\mathrm{cyc}}.
\]
\end{lemma}

\begin{proof}
Both sides satisfy the cyclic Segal condition, so it is enough to
check the comparison on corollas and on the exceptional edge. On the
exceptional edge, both sides are terminal. At $C_1$, both sides are
empty, so the comparison is an isomorphism. For $n\geq1$, the
comparison at $C_{n+1}$ is the map
$\bigl(\cs\PaRB^{\mathrm{cyc}}(n)\bigr)^{\wedge}
\to
\widehat{\cs}\bigl(\widehat{\PaRB}^{\mathrm{cyc}}(n)\bigr)$.
It is a weak equivalence by
\cite[Proposition~5.9]{Horel_profinite_groupoids}, since
$\PaRB^{\mathrm{cyc}}(n)$ has finitely many objects and good
automorphism groups. The cyclic Segal conditions then give the
result on every tree.
\end{proof}

The preceding lemmas allow us to pass from homotopy endomorphisms of
the completed ribbon-braid model in profinite groupoids to homotopy
endomorphisms of its completed classifying space in profinite spaces.

\begin{thm}
\label{thm: GT acts on widehat csN PaRB}
There is an isomorphism of monoids
\[
\underline{\GT}
\cong
\pi_0\mathbb R\End_{\Cyc_\infty(\widehat{\sSet})}
\left(
\bigl(\nerve\mathsf S_0\bigr)^{\wedge}
\right).
\]
\end{thm}

\begin{proof}
The profinite classifying-space functor is homotopically fully
faithful by
\cite[Proposition~5.7]{Horel_profinite_groupoids}. Applying it
entrywise and passing to the cyclic Segal localisations therefore
gives an isomorphism
\[
\pi_0\mathbb R\End_{\Cyc_\infty(\widehat{\Grpd})}
\bigl(\nerve\widehat{\PaRB}^{\mathrm{cyc}}\bigr)
\cong
\pi_0\mathbb R\End_{\Cyc_\infty(\widehat{\sSet})}
\bigl(\widehat{\cs}\nerve\widehat{\PaRB}^{\mathrm{cyc}}\bigr).
\]
Combining Theorem~\ref{thm: GT as homotopy cyclic endomorphisms}
with Lemma~\ref{lem:comparison cyclic PaRB nerves} therefore
identifies this monoid with the homotopy endomorphisms of
$\bigl(\cs\nerve\PaRB^{\mathrm{cyc}}\bigr)^{\wedge}$.

The comparison established in
Section~\ref{subsec:framed_conf_operad} gives a zigzag of entrywise
weak equivalences of cyclic $\infty$-operads
$\cs\nerve\PaRB^{\mathrm{cyc}}\simeq\nerve\mathsf S_0$.
Since every simplicial set is cofibrant and profinite completion is
left Quillen, entrywise completion preserves each weak equivalence
in this zigzag. It follows in particular that
$\bigl(\nerve\mathsf S_0\bigr)^{\wedge}$ is a cyclic
$\infty$-operad in profinite spaces. Homotopy invariance of derived
endomorphism spaces now gives the result.
\end{proof}

Passing to units gives the corresponding statement for homotopy
automorphisms.

\begin{cor}
\label{cor:Galois_action_CB}
There is an isomorphism of groups
\[
\GT
\cong
\pi_0\mathbb R\Aut_{\Cyc_\infty(\widehat{\sSet})}
\left(
\bigl(\nerve\mathsf S_0\bigr)^{\wedge}
\right).
\]
Consequently, the injection
$\operatorname{Gal}(\overline{\mathbb Q}/\mathbb Q)\hookrightarrow
\GT$ induces a faithful action of the absolute Galois group by
homotopy automorphisms of
$\bigl(\nerve\mathsf S_0\bigr)^{\wedge}$.
\end{cor}

\section{Galois actions on tangles}
\label{sec: Galois action on tangles}

In this section, we construct a $\GT$-action on parenthesised profinite tangles through topologically enriched functors. The action transports their metric-prop algebra structure, allowing the ribbon operations to vary.

We first recall the categories of framed tangles (Section~\ref{subsec:tangles}) and extend a $\Env(\PaRB)$-action on $q\bT'$ to an action of
$\Pi(\PaRB^{\mathrm{cyc}})$ through its endomorphism prop. Using Furusho's construction of framed profinite tangles \cite{furusho_galois_2017}, we then obtain the topological category $q\widehat{\bT}'$ and the corresponding metric-prop action (Section~\ref{subsec: profinite tangles}). The $\GT$-action on $\widehat{\PaRB}^{\mathrm{cyc}}$ extends
to its metric prop and acts by precomposition on this algebra structure (Section~\ref{subsec:GT_action_profinite_tangles}).

In Section~\ref{subsec:GT action on turaev category}, we will explain how in characteristic zero, this $\GT$-action is
realized by ribbon monoidal equivalences of unoriented tangles.

\subsection{The ribbon category of tangles}
\label{subsec:tangles}

A framed oriented tangle is an isotopy class, relative to the
boundary, of proper embeddings of compact oriented one-manifolds
in the cube $[-1,1]^3$, equipped with framings of their normal
bundles. The boundary points of the interval components lie on
$[-1,1]\times\{0\}\times\{\pm1\}$ and are arranged in a fixed order.

Let $\mathbb M(\{+,-\})$ be the free monoid on the symbols $+$ and
$-$. Reading the orientations of the boundary strands in their
fixed order assigns source and target words to a tangle. The
category $\bT$ has these signed words as objects and framed
oriented tangles as morphisms. Composition is given by vertical
stacking, tensor product by horizontal juxtaposition, and the
monoidal unit is the empty word~$\emptyset$.

The dual word $w^*$ is obtained by reversing the order of the
symbols in $w$ and interchanging $+$ and $-$. For example,
$(+-++-)^*=+--+-$, and in particular $+^*=-$ and $-^*=+$.
With the left-duality convention used throughout this paper,
the basic coevaluation and evaluation are
$\coev_+\colon\emptyset\to(+,-)$ and
$\ev_+\colon(-,+)\to\emptyset$. They are represented by a cup
and a cap and satisfy
\[
(\id_+\otimes\ev_+)\cdot(\coev_+\otimes\id_+)=\id_+,
\qquad
(\ev_+\otimes\id_-)\cdot(\id_-\otimes\coev_+)=\id_-.
\]

The elementary positive crossing defines the braiding
$c_{+,+}\colon(+,+)\to(+,+)$, and a positive full twist
of a framed strand defines $\Theta_+\colon+\to+$.
Tensor products and duality determine the corresponding
braidings and twists for arbitrary signed words.
The twist satisfies the balancing identity
\[
\Theta_{X\otimes Y}
=
c_{Y,X}\cdot c_{X,Y}\cdot(\Theta_X\otimes\Theta_Y).
\]
These structures make $\bT$ the strict ribbon monoidal category
freely generated by one object. See \cite{shum_tortile_1994},
\cite[Chapter~XII.2]{kassel_quantum_1995}, and
\cite[Section~3.5]{KT_Tangles}.

The category $q\bT$ of \emph{$q$-tangles}, or
\emph{parenthesised tangles}, is the parenthesised form of $\bT$
used in the theory of universal finite-type invariants. See
\cite{BN_weak_tangles,Le_Murakami_1996,KT_Tangles,qtanglesbook}.
Its objects are fully parenthesised words in $+$ and $-$,
together with the empty word. A morphism $s\to t$ is a framed
oriented tangle whose source and target are obtained by
forgetting the parentheses in $s$ and~$t$.

Composition and tensor product are again given by vertical
stacking and horizontal juxtaposition. Tensor product of
nonempty words adjoins an outer pair of parentheses, while
the empty word is a strict unit. The associator is
\[
\alpha_{s,t,u}\colon
(s\otimes t)\otimes u
\longrightarrow
s\otimes(t\otimes u).
\]
The underlying tangle of $\alpha_{s,t,u}$ consists of parallel
identity strands. Together with the braiding, twist, evaluation,
and coevaluation inherited from $\bT$, this makes $q\bT$ a
ribbon monoidal category. Forgetting parentheses defines a
ribbon monoidal equivalence $q\bT\to\bT$.

We will also use the unoriented framed tangle category. See
\cite[Remark~4.1 and Section~8.2]{mousaaid2021affinization}.

\begin{definition}
\label{def: self-dual tangle category}
Let $\bT'$ be the strict ribbon category obtained from $\bT$
by forgetting strand orientations. Its objects are words in
a single symbol $*$, including the empty word, and its
morphisms are framed unoriented tangles with prescribed source
and target. The object $*$ is strictly self-dual, with
evaluation and coevaluation represented by the cap and cup.
\end{definition}

Equivalently, $\bT'$ can be viewed as the quotient category $\bT$ equipped with the strict self-duality relation $(+)^* = +$.

Define $q\bT'$ from $\bT'$ as $q\bT$ was defined from $\bT$.
Its objects are fully parenthesised words in $*$, together
with the empty word, and
\[
\Hom_{q\bT'}(p,q)
=
\Hom_{\bT'}(\overline p,\overline q),
\]
where $\overline p$ and $\overline q$ are obtained by
forgetting parentheses. Forgetting parentheses defines a
ribbon monoidal equivalence $q\bT'\to\bT'$.

\subsubsection{Relating tangles to cyclic ribbon braids}

By Theorem~\ref{thm: CoRB_and_PaRB_algebras}, the underlying
balanced monoidal structure on $q\bT'$ makes it a $\PaRB$-algebra,
equivalently an $\Env(\PaRB)$-algebra. We extend this action to
$\Pi(\PaRB^{\mathrm{cyc}})$ using the symmetric nondegenerate
pairing supplied by ribbon duality.

Let $\calC$ be a small ribbon monoidal category with chosen left
duality $X\mapsto X^*$. The functor
$(X,Y)\mapsto\Hom(X,Y^*)$ is contravariant in both variables,
by precomposition in the first variable and duality in the second.
We use ribbon duality to express categorical modules
(profunctors) with all variables contravariant, defining
endomorphism categories with $n$ inputs and $m$ outputs.
This convention allows the cyclic action to permute the input
and output variables and represents evaluation and coevaluation
by $\Hom(X,Y^*)$ in their respective profiles.

\begin{definition}
\label{def: categorical endomorphism prop}
For $n,m\geq0$, define
\[
\End_{\calC}(n,m)
=
\Fun\bigl((\calC^{\times(m+n)})^{\mathrm{op}},\Set\bigr).
\]
We write its objects as
$H(Y_1,\ldots,Y_m;X_1,\ldots,X_n)$, listing outputs before
inputs. Its morphisms are natural transformations.

For $H\in\End_{\calC}(n,m)$ and
$K\in\End_{\calC}(m,\ell)$, composition glues the outputs
of $H$ to the corresponding inputs of $K$, using duality:
\[
\begin{aligned}
&(K\cdot H)(Z_1,\ldots,Z_\ell;X_1,\ldots,X_n)
&=
\int^{Y_1,\ldots,Y_m\in\calC}
K(Z_1,\ldots,Z_\ell;Y_1,\ldots,Y_m)
\times
H(Y_1^*,\ldots,Y_m^*;X_1,\ldots,X_n).
\end{aligned}
\]
The coend ranges over the intermediate objects and identifies
applying a morphism on one side of a gluing with applying its
dual on the other. Explicitly, for $f\colon Y\to Y'$, the
identification is
\[
\bigl(Y,K(f)(k),h\bigr)
\sim
\bigl(Y',k,H(f^*)(h)\bigr),
\]
where $k\in K(\ldots;Y')$ and $h\in H(Y^*;\ldots)$.
These identifications are imposed in each intermediate variable.

Tensor product is the Cartesian product of functors on separate
lists of variables. The identity in $\End_{\calC}(n,n)$ is
the functor
$(Y_1,\ldots,Y_n;X_1,\ldots,X_n)
\mapsto\prod_{i=1}^n\Hom(Y_i,X_i^*)$.
The symmetric groups act by permuting the input and output
variables.
\end{definition}

The co-Yoneda isomorphisms give the unit comparisons, and
interchange of coends gives the associativity comparisons.
Together with the double-dual identifications supplied by the
ribbon structure, these make $\End_{\calC}$ a prop in categories
up to coherent natural isomorphism. We use these double-dual
identifications throughout.

We first describe the given $\PaRB$-action through
$\End_{q\bT'}$. An object $\omega\in\ob(\PaRB(n))$ acts
by the tensor-product functor
$(w_1,\ldots,w_n)\mapsto\omega(w_1,\ldots,w_n)$.
The corresponding object of $\End_{q\bT'}(n,1)$ is
\[
(w_0;w_1,\ldots,w_n)\longmapsto
\Hom\bigl(w_0,\omega(w_1^*,\ldots,w_n^*)\bigr).
\]
Here each $w_i$ is a fully parenthesised word in $*$, or the
empty word. The target word is obtained by substituting $w_i^*$
for the letter $i$ in $\omega$.

\begin{example}
The binary corolla $\mu$ determines the functor
$(w_0;w_1,w_2)\mapsto\Hom(w_0,w_1^*\otimes w_2^*)$.
For $w_0=*$, $w_1=(**)$, and $w_2=*$, its value is
$\Hom\bigl(*,(**)*\bigr)$,
since both $(**)$ and $*$ are self-dual. This is the set of
framed unoriented tangles from one strand to three strands,
with target parenthesised as $(**)*$. For example,
$\coev_*\otimes\id_*\colon *\to(**)*$ is an element of this
set, represented by a cup beside an identity strand.
\end{example}

The map of props $\Env(\PaRB)\to\End_{q\bT'}$ sends the
associator, braiding, and twist of $\PaRB$ to natural
transformations given by postcomposition:
\[
\begin{aligned}
\alpha_{123}
&\longmapsto
\bigl[u\mapsto\alpha_{w_1^*,w_2^*,w_3^*}\cdot u\bigr],\\
\beta_{12}
&\longmapsto
\bigl[u\mapsto c_{w_1^*,w_2^*}\cdot u\bigr],\\
\twist
&\longmapsto
\bigl[u\mapsto\Theta_{w_1^*}\cdot u\bigr].
\end{aligned}
\]
For example, fix the input words $w_1,w_2$. The transformation
assigned to $\beta_{12}$ exchanges the entire blocks $w_1^*$
and $w_2^*$, preserving the parenthesisation within each block.
It sends a tangle $u\colon w_0\to w_1^*\otimes w_2^*$ to
the composite
\[
w_0\xrightarrow{u}w_1^*\otimes w_2^*
\xrightarrow{c_{w_1^*,w_2^*}}w_2^*\otimes w_1^*.
\]
The source word $w_0$ varies independently of the input words.
For instance, when $w_0=*$, $w_1=(**)$, and $w_2=*$, this
sends a tangle $*\to(**)*$ to a tangle $*\to*(**)$ by
exchanging the two-strand block with the one-strand block.

By the Yoneda lemma, the natural transformation is determined
by its value on the identity at
$w_0=w_1^*\otimes w_2^*$. This value is precisely
$c_{w_1^*,w_2^*}$, so the transformation recovers the block
braiding. The same argument recovers the associator and twist.

\begin{prop}
\label{prop: metric prop action through cyclic endomorphisms}
The $\Env(\PaRB)$-action on $q\bT'$ extends to an action of
$\Pi(\PaRB^{\mathrm{cyc}})$ through $\End_{q\bT'}$.
\end{prop}

\begin{proof}
We first verify that the functors and natural transformations
described above define the $\Env(\PaRB)$-action.
For $\omega\in\ob(\PaRB(n))$ and
$\omega'\in\ob(\PaRB(m))$, an element of the coend describing
substitution is represented by morphisms
\[
u\colon w_0\longrightarrow\omega(\ldots,w^*,\ldots),
\qquad
v\colon w^*\longrightarrow\omega'(\ldots).
\]
Send this representative to $\omega(\ldots,v,\ldots)\cdot u$.
The coend relation identifies the two ways of applying a morphism
of the intermediate word. Conversely, a morphism
$w_0\to\omega(\ldots,\omega'(\ldots),\ldots)$ is represented
by taking $w=\omega'(\ldots)^*$ and
$v=\id_{\omega'(\ldots)}$, using the double-dual identification.
The coend relations identify every representative with one of
this form, giving a bijection.
Thus composition agrees with operadic substitution.

The relations on generating morphisms hold by the given
$\PaRB$-action. Taking Cartesian products on separate lists of
variables and permuting those variables gives the corresponding
$\Env(\PaRB)$-action. Evaluating at dual input words and using
$w_i^{**}\cong w_i$ recovers the original $\PaRB$-action.

To extend this action to the metric prop, assign both $\ev$
and $\coev$ the pairing
$(w,w')\mapsto\Hom(w,(w')^*)$, with its variables regarded
as inputs for $\ev$ and outputs for $\coev$.
We verify the symmetry, snake, and cyclic-invariance relations.

For $h\colon w\to(w')^*$, its transpose, defined in
\eqref{eq: transpose of a morphism}, is
$h^*\colon(w')^{**}\to w^*$.
The symmetry $\Hom(w,(w')^*)\cong\Hom(w',w^*)$ sends $h$
to the composite
\[
w'\xrightarrow{\cong}(w')^{**}\xrightarrow{h^*}w^*.
\]
Applying this construction twice returns $h$ under the
double-dual identifications.

For the snake relations, fix words $w_0,w_1$. An element of
the first snake composite evaluated at $(w_0;w_1)$ is
represented by a chain
\[
w_0\xrightarrow{h_1}w_2^*
\xrightarrow{h_2}w_3
\xrightarrow{h_3}w_4^*
\xrightarrow{h_4}w_1^*.
\]
Composition sends its class to
$h_4\cdot h_3\cdot h_2\cdot h_1\in\Hom(w_0,w_1^*)$.
An inverse represents $u\colon w_0\to w_1^*$ by taking
$w_2=w_4=w_0^*$, $w_3=w_0$, the first three morphisms to
be identities, and the last to be $u$.
The coend relations identify every chain with this representative.
The second snake composite gives the same bijection after
transposing the factors. Hence both composites are isomorphic
to the identity operation
$(w_0;w_1)\mapsto\Hom(w_0,w_1^*)$
in $\End_{q\bT'}(1,1)$.

For cyclic invariance, let $H$ be a contravariant functor of
$n+1$ variables, with $n\geq1$. Contracting the identity
factors in the image of the cyclic-invariance relation gives
\[
\begin{aligned}
&\int^{w,w'}
\Hom(w_0,w)\times
H(w',w,w_1,\ldots,w_{n-1})\times
\Hom(w_n,w')\\
&\qquad\cong
H(w_n,w_0,\ldots,w_{n-1})
=(z_{n+1}^*H)(w_0,\ldots,w_n).
\end{aligned}
\]
A representative consists of $a\colon w_0\to w$,
$b\colon w_n\to w'$, and
$h\in H(w',w,w_1,\ldots,w_{n-1})$.
The displayed map sends its class to
$H(b,a,\id,\ldots,\id)(h)$.
Its inverse takes $w=w_0$, $w'=w_n$, and $a,b$ to be
identities. Thus the images of $\ev$ and $\coev$ implement
the required cyclic permutation.

For the functor assigned to $\omega$, duality gives
\[
\Hom\bigl(w_0,\omega(w_1^*,\ldots,w_n^*)\bigr)
\cong
\Hom\bigl(
\emptyset,w_0^*\otimes\omega(w_1^*,\ldots,w_n^*)
\bigr).
\]
Under this identification, changing the output is implemented
by bending strands through cups and caps. On objects, this
reroots the planar tree representing $\omega$. On generating
morphisms, the bending comparisons give
\[
s_0^*\twist=\twist,\qquad
s_0^*\beta_{12}
=(1_\mu\circ_2\twist^{-1})\cdot\beta_{12}^{-1},
\qquad
s_0^*\alpha_{123}=\alpha_{231}^{-1},
\]
as in \eqref{eq: cyclic action PaRB generators} and
Figure~\ref{fig:cyclic_action_on_PaRB}.
These comparisons respect gluing and relabelling by cup--cap
cancellation and ribbon coherence. Since the displayed morphisms
generate $\PaRB$, they establish compatibility with its cyclic
structure.

The assignments therefore preserve the defining relations of
$\Pi(\PaRB^{\mathrm{cyc}})$. The coend and duality comparisons
are compatible with iterated composition and tensor product,
giving the asserted extension.
\end{proof}

\subsection{Profinite tangles}
\label{subsec: profinite tangles}

Furusho defines spaces of profinite oriented tangles and knots in
\cite{furusho_galois_2017}. The framed version is introduced in
Appendix~B of the same paper. We organize the framed profinite
tangles according to their source and target and thereby obtain a
category enriched in compactly generated spaces.

Following
\cite[Definitions~2.2 and~2.3]{furusho_galois_2017}, let $A$ and $C$
be the sets of oriented annihilation and creation symbols
\[
A=
\left\{
a^\epsilon_{k,l}
\ \middle|\
k,l\geq0,\ 
\epsilon\in
\{\uparrow,\downarrow\}^{k}
\times
\{\underset{\rightarrow}{\cap},
  \underset{\leftarrow}{\cap}\}
\times
\{\uparrow,\downarrow\}^{l}
\right\}
\]
and
\[
C=
\left\{
c^\epsilon_{k,l}
\ \middle|\
k,l\geq0,\ 
\epsilon\in
\{\uparrow,\downarrow\}^{k}
\times
\{\underset{\rightarrow}{\cup},
  \underset{\leftarrow}{\cup}\}
\times
\{\uparrow,\downarrow\}^{l}
\right\}.
\]
The set of decorated profinite braid pieces is
\[
\widehat{\mathcal B}
=
\coprod_{n\geq1}
\widehat{\Br}_n\times\{\uparrow,\downarrow\}^n.
\]
For $(\beta,\epsilon)\in\widehat{\mathcal B}$, the incoming
orientations are specified by $\epsilon$, and the outgoing
orientations are obtained from $\epsilon$ using the permutation
underlying~$\beta$.

A \emph{profinite tangle diagram} is a composable word
$T=\gamma_m\cdot\ldots\cdot\gamma_1$, where each
$\gamma_i$ belongs to
$A\amalg\widehat{\mathcal B}\amalg C$ and the target of
$\gamma_i$ agrees with the source of $\gamma_{i+1}$. Its source is
the source of $\gamma_1$, and its target is the target of
$\gamma_m$. We also allow the empty composable word at each
object, representing its identity morphism.

A \emph{framed profinite tangle} is an equivalence class of such
diagrams under the moves
$\mathrm{(FT1)}$--$\mathrm{(FT6)}$ of
\cite[Definition~B.1]{furusho_galois_2017}. The moves
$\mathrm{(FT1)}$--$\mathrm{(FT5)}$ are the corresponding moves
$\mathrm{(T1)}$--$\mathrm{(T5)}$ of
\cite[Definition~2.8]{furusho_galois_2017}. They give the unit and
composition relations for braid pieces, interchange of independent
tangles, compatibility with cabling, and the creation--annihilation
relations. The move $\mathrm{(FT6)}$ is the framed first
Reidemeister move.

The framed theory also contains distinguished twists
$\varphi^\uparrow\colon+\to+$ and
$\varphi^\downarrow\colon-\to-$, defined in
\cite[Definition~B.3]{furusho_galois_2017}. Lemma~B.9 of that paper
extends their powers continuously from $\mathbb N$ to
$\widehat{\mathbb Z}$, while Lemma~B.10 describes their compatibility
with cups and caps.

\begin{definition}
The \emph{framed profinite tangle category} $\widehat{\bT}$ has the
signed words in $\{+,-\}$ as objects. Its morphisms from $w$ to $w'$
are the framed profinite tangles with source $w$ and target~$w'$.
Composition is induced by concatenation of diagrams, and tensor
product is induced by horizontal juxtaposition.

Each morphism space is given the compactly generated refinement of
Furusho's quotient topology. This topology is induced by the
profinite topologies on the braid pieces and the discrete topologies
on the creation and annihilation pieces, as in
\cite[Note~2.10 and Appendix~B]{furusho_galois_2017}.
\end{definition}

Concatenation and juxtaposition are continuous on tangle diagrams
and preserve the framed tangle relations. They therefore induce
continuous composition and tensor product on $\widehat{\bT}$.
The braid pieces, cups, caps, and twists make $\widehat{\bT}$
a strict ribbon monoidal category. The braid relations hold in
the profinite braid groups, the two instances of $\mathrm{(FT5)}$
give the snake identities, and $\mathrm{(FT6)}$ governs the
framing twists.

Let $\widehat{\bT}'$ be the corresponding category of framed
unoriented profinite tangles. Its objects are the natural numbers,
composition is vertical stacking, and tensor product is horizontal
juxtaposition. Its generating object is self-dual, with evaluation
and coevaluation represented by the cap and cup.

Define $q\widehat{\bT}'$ from $\widehat{\bT}'$ as $q\bT'$ was
defined from $\bT'$. Its objects are fully parenthesised words in
$*$, together with the empty word. If $|p|$ denotes the length of
the word obtained from $p$ by forgetting parentheses, then
\[
\Hom_{q\widehat{\bT}'}(p,q)
=
\Hom_{\widehat{\bT}'}(|p|,|q|).
\]
The associator is the identity underlying tangle between the
corresponding parenthesisations. Forgetting parentheses defines
a topological ribbon monoidal equivalence
$q\widehat{\bT}'\to\widehat{\bT}'$.

The $\PaRB$-action on parenthesised tangles extends to an
action of $\widehat{\PaRB}$ on $q\widehat{\bT}'$.
On objects, a parenthesised word acts by the corresponding
tensor-product functor. On morphisms, a profinite ribbon braid
acts by cabling, with the framing parameters interpreted using
the continuous powers of twists from
\cite[Lemma~B.9]{furusho_galois_2017}.
The composition and cabling relations ensure that these
assignments define an operad action, continuous in the
profinite ribbon-braid parameters.

We apply Definition~\ref{def: categorical endomorphism prop}
to the topological category $q\widehat{\bT}'$, using
topologically enriched functors and enriched coends.
We denote the resulting endomorphism prop by
$\End_{q\widehat{\bT}'}$.

\begin{prop}
\label{prop: profinite metric prop action through cyclic endomorphisms}
The action of $\widehat{\PaRB}$ on $q\widehat{\bT}'$
extends to an action of
$\Pi(\widehat{\PaRB}^{\mathrm{cyc}})$ through
$\End_{q\widehat{\bT}'}$.
Its restriction to $\Env(\widehat{\PaRB})$ is induced by the
ribbon structure on $q\widehat{\bT}'$. Both $\ev\colon 2\to 0$
and $\coev\colon 0\to 2$ are sent to the continuous functor
\[
(w,w')\longmapsto
\Hom_{q\widehat{\bT}'}\bigl(w,(w')^*\bigr),
\]
with $w,w'$ regarded as input variables for $\ev$ and as output
variables for $\coev$.
\end{prop}

\begin{proof}
The construction of
Proposition~\ref{prop: metric prop action through cyclic endomorphisms}
applies with hom-sets replaced by the corresponding mapping spaces.
Composition and tensor product in $q\widehat{\bT}'$, as well as
duality and the action of profinite ribbon braids, are continuous.
The enriched co-Yoneda isomorphisms therefore identify composition
in $\End_{q\widehat{\bT}'}$ with operadic substitution in
$\widehat{\PaRB}$.

Assign $\ev$ and $\coev$ the functor displayed in the statement.
Transposition gives a continuous isomorphism
\[
\Hom_{q\widehat{\bT}'}\bigl(w,(w')^*\bigr)
\cong
\Hom_{q\widehat{\bT}'}(w',w^*),
\]
which satisfies the symmetry relations. The enriched co-Yoneda
isomorphisms identify both snake composites with the identity
operation, as in
Proposition~\ref{prop: metric prop action through cyclic endomorphisms}.
The duality used in these identifications is supplied by the
cup--cap relations $\mathrm{(FT5)}$.

The cyclic-invariance relation is verified as in
Proposition~\ref{prop: metric prop action through cyclic endomorphisms}.
Contraction with the functors assigned to $\coev$ and $\ev$
changes the choice of output, and the resulting action on
$\twist$, $\beta_{12}$, and $\alpha_{123}$ is given by
\eqref{eq: cyclic action PaRB generators}. For profinite
ribbon-braid morphisms, compatibility with composition and
cabling follows from $\mathrm{(FT2)}$ and $\mathrm{(FT4)}$,
while compatibility of profinite twist powers with cups and
caps follows from
\cite[Lemma~B.10]{furusho_galois_2017}.

Consequently, the assignments preserve the symmetry, snake, and
cyclic-invariance relations defining
$\Pi(\widehat{\PaRB}^{\mathrm{cyc}})$. All the structure maps are
continuous, so these assignments give the asserted action.
\end{proof}

\subsection{Galois actions on profinite tangles}
\label{subsec:GT_action_profinite_tangles}

We identify the object-fixing automorphisms of the metric
prop $\Pi(\widehat{\PaRB}^{\mathrm{cyc}})$ and use them to
transport the metric-prop algebra structure on
$q\widehat{\bT}'$ constructed in
Proposition~\ref{prop: profinite metric prop action through cyclic endomorphisms}.

By Theorem~\ref{thm: GT monoid iso to cyclic endomorphisms}, each
$(\lambda,f)\in\GT$ determines an object-fixing cyclic operad automorphism
$F_{\lambda,f}$ of $\widehat{\PaRB}^{\mathrm{cyc}}$. On the
generators, it is given by
\begin{equation}
\label{eq: GT action generators profinite PaRB}
F_{\lambda,f}(\beta_{12})
=
\beta_{12}\cdot
(\beta_{21}\cdot\beta_{12})^{(\lambda-1)/2},
\qquad
F_{\lambda,f}(\twist)=\twist^\lambda,
\qquad \text{and} \qquad
F_{\lambda,f}(\alpha_{123})
=
\alpha_{123}\cdot f(x_{12},x_{23}).
\end{equation}
For a prop in groupoids, an object-fixing automorphism means an
automorphism which is the identity on the objects of each of its
entry groupoids.

\begin{prop}
\label{prop:GT-metric-prop-profinite}
Restriction induces natural isomorphisms
\[
\underline{\Aut}_{\Prop}\bigl(\Pi(\widehat{\PaRB}^{\mathrm{cyc}})\bigr)
\cong
\underline{\Aut}_{\Cyc}\bigl(\widehat{\PaRB}^{\mathrm{cyc}}\bigr)
\cong
\GT.
\]
\end{prop}

\begin{proof}
Let $F$ be an object-fixing cyclic operad automorphism of
$\widehat{\PaRB}^{\mathrm{cyc}}$. Apply $\Env(F)$ to
$\Env(\widehat{\PaRB})$ and fix $\ev$ and $\coev$. The relations
in $\Env(\widehat{\PaRB})$ are preserved because $F$ is an operad
automorphism, while the snake and symmetry relations are preserved
because $\ev$, $\coev$, and the structural symmetries are fixed.

For $h\in\widehat{\PaRB}^{\mathrm{cyc}}(n)$, cyclic invariance is
the relation
\[
z_{n+1}^*h
=
(\id_1\otimes\ev)
\cdot(\id_1\otimes h\otimes\id_1)
\cdot(\coev\otimes\id_n).
\]
The image of its left-hand side is
$F(z_{n+1}^*h)=z_{n+1}^*F(h)$, while the image of its right-hand
side is
\[
(\id_1\otimes\ev)
\cdot(\id_1\otimes F(h)\otimes\id_1)
\cdot(\coev\otimes\id_n).
\]
These are equal by cyclic invariance for $F(h)$. Thus $F$ extends
to an object-fixing prop endomorphism $\Pi(F)$.
The extension of $F^{-1}$ is its inverse, so $\Pi(F)$ is an
automorphism.

Conversely, let $\Phi$ be an object-fixing automorphism of
$\Pi(\widehat{\PaRB}^{\mathrm{cyc}})$. Since
$\Env(\widehat{\PaRB})$ is canonically an entrywise full
subprop, $\Phi$ restricts to an object-fixing automorphism
of $\widehat{\PaRB}$. It also fixes $\ev$ and $\coev$.
Applying $\Phi$ to the cyclic-invariance relation therefore
shows that its restriction commutes with the cyclic action.

The defining universal property of the metric prop makes
an extension fixing $\ev$ and $\coev$ unique. Thus extension
and restriction are inverse and respect composition, proving
the first isomorphism. The second follows from
Theorem~\ref{thm: GT monoid iso to cyclic endomorphisms}
by passing to invertible elements.
\end{proof}

Let $\rho\colon
\Pi(\widehat{\PaRB}^{\mathrm{cyc}})
\to\End_{q\widehat{\bT}'}$
be the action constructed in
Proposition~\ref{prop: profinite metric prop action through cyclic endomorphisms}.
For each $(\lambda,f)\in\GT$, precomposition gives another action
\[
\Pi(\widehat{\PaRB}^{\mathrm{cyc}})
\xrightarrow{\ \Pi(F_{\lambda,f})\ }
\Pi(\widehat{\PaRB}^{\mathrm{cyc}})
\xrightarrow{\ \rho\ }
\End_{q\widehat{\bT}'}.
\]
On $\Env(\widehat{\PaRB})$, this replaces each parenthesised
ribbon-braid operation by its image under $F_{\lambda,f}$.
The functors assigned to $\ev$ and $\coev$ remain unchanged because
$\Pi(F_{\lambda,f})$ fixes these operations.

Compatibility with multiplication in $\GT$, that is $\Pi(F_{\lambda_1,f_1}) \cdot \Pi(F_{\lambda_2,f_2}) = \Pi(F_{(\lambda_1,f_1)\ast (\lambda_2,f_2)})$,  makes precomposition
a right action on the metric-prop algebra structures carried by
$q\widehat{\bT}'$. We refer to this transport as the $\GT$-action
on $q\widehat{\bT}'$.

\begin{cor}
\label{cor:Galois_action_qT}
The group $\GT$ acts on parenthesised profinite tangles
$q\widehat{\bT}'$. Consequently, the absolute Galois group
$\operatorname{Gal}(\overline{\mathbb Q}/\mathbb Q)$
also acts on $q\widehat{\bT}'$.
\end{cor}

\begin{proof}
The action is given by transport of the metric-prop algebra
structure through $\End_{q\widehat{\bT}'}$.
Precomposition with $\Pi(F_{\lambda,f})$ preserves the algebra
axioms, the structure comparisons, and continuity of the
structure maps. The multiplication law for these automorphisms
gives the right $\GT$-action, with the identity acting trivially.
Restriction along Ihara's homomorphism
$\operatorname{Gal}(\overline{\mathbb Q}/\mathbb Q)\to\GT$
gives the Galois action.
\end{proof}

\begin{remark}
\label{remark: comparison with Kassel Turaev}
Kassel and Turaev construct an action of the proalgebraic
Grothendieck--Teichmüller group on completed categories of tangles in
\cite[Appendix~D]{KT_Tangles}. Their formulas modify the associator,
braiding, and twist and introduce corresponding corrections to the
duality morphisms. See
\cite[Equations~\textnormal{(D.2)}--\textnormal{(D.7)}]{KT_Tangles}. These corrections involve inverses which exist in the prounipotent linear completion but need not exist in the profinite tangle category. A related issue is noted by Furusho in
\cite[Remark~2.43]{furusho_galois_2017}: the inverse of the
profinite knot $\Lambda_f$ is naturally available after passing
to the fraction group of the monoid of profinite knots.

The action of Corollary~\ref{cor:Galois_action_qT} is defined by
precomposition with $\Pi(F_{\lambda,f})$. This automorphism fixes
the formal operations $\ev$ and $\coev$ and preserves the
cyclic-invariance relations. Transport of the metric-prop algebra
structure therefore requires no inverse of a profinite knot.

To realize this transport by ribbon monoidal equivalences between
the original and transported ribbon structures, the construction
requires invertible correction morphisms for the duality.
These are available in the prounipotent linear completion.
In Section~\ref{subsec:GT action on turaev category}, we identify
them with the Kassel--Turaev corrections and obtain the
corresponding ribbon monoidal equivalences.
\end{remark}

\section{Turaev's tangle category and rational formality of the
framed little disks operad}
\label{sec: Turaev's tangle category and rational formality of the little disks operad}

Let $\kk$ be a field of characteristic zero. In this section, we
study the action of the proalgebraic Grothendieck--Teichmüller group
$\gt(\kk)$ on the completed $\kk$-linear category of framed tangles
introduced by Kassel and Turaev. We first compare their completion
with the Malcev completion of the braid and ribbon-braid groups
(Section~\ref{sec: completion of categories, operads and props}).
We then describe the construction of
\cite[Appendix~D]{KT_Tangles}, including the correction to the
duality maps, and explain its relation to the action on the cyclic
operad $\PaRB_{\kk}^{\mathrm{cyc}}$
(Section~\ref{subsec:GT action on turaev category}).
Finally, we use the $\gt(\mathbb{Q})$-action on this cyclic operad
to prove the cyclic formality of the framed little disks operad
(Section~\ref{subsec: formality}).

\subsection{Prounipotent completion of categories, operads, and props}
\label{sec: completion of categories, operads and props}

We recall the two characteristic-zero completions used below. Malcev
completion applies aritywise to the braid and ribbon-braid operads,
whereas the completed tangle category of Kassel and Turaev is defined
using the augmentation ideal of a braided linear category. On pure
braids, the two constructions induce the same filtration. For ribbon
braids, the Kassel--Turaev completion also retains the parities of the
strand framings.

\subsubsection{Malcev completion of groups}

Let $G$ be a group and let $I_G$ be the augmentation ideal of its group algebra $\kk[G]$. Its completed group algebra is
\[\widehat{\kk[G]}:=\varprojlim_r\kk[G]/I_G^r.\] The \emph{Malcev completion} of $G$ over $\kk$ is the group $G_{\kk}:=\mathcal G(\widehat{\kk[G]})$ of grouplike elements; see \cite[Chapter~8]{FresseBook1}.

The primitive elements of $\widehat{\kk[G]}$ form a pronilpotent Lie algebra $\mathfrak g_{\kk}$. In characteristic zero, exponential and logarithm give mutually inverse bijections between $\mathfrak g_{\kk}$ and $G_{\kk}$; see \cite[Proposition~8.1.5]{FresseBook1}. For $g\in G_{\kk}$ and $a\in\kk$, we therefore write $g^a:=\exp(a\log g)$.

\subsubsection{Completion of braided linear categories}

Let $\calC$ be a $\kk$-linear braided monoidal category.
A \emph{monoidal ideal} $\mathcal I$ in $\calC$ consists of
subspaces $\mathcal I(X,Y)\subseteq\Hom_{\calC}(X,Y)$ closed
under precomposition, postcomposition, and tensoring with
arbitrary morphisms.

The \emph{augmentation ideal} of $\calC$ is the monoidal
ideal generated by the braiding defects
$c_{W,V}\cdot c_{V,W}-\id_{V\otimes W}$.
Following \cite[Section~3.4]{KT_Tangles}, its
\emph{prounipotent completion} is
$\widehat{\calC}:=\varprojlim_r\calC/\mathcal I^{r+1}$.
It has the same objects as $\calC$, with morphism spaces
\[
\Hom_{\widehat{\calC}}(X,Y)
=
\varprojlim_r
\Hom_{\calC}(X,Y)/\mathcal I^{r+1}(X,Y).
\]

\begin{example}
\label{example: completion of braid category}
Let $\bB$ be the braid category, with objects the natural
numbers, $\End_{\bB}(n)=\Br_n$, and no morphisms between
distinct objects. Tensor product is addition on objects
and juxtaposition on morphisms.

In its $\kk$-linearisation, the braided augmentation ideal
on $\kk[\Br_n]$ is generated by $\beta_i-\beta_i^{-1}$,
or equivalently by $\beta_i^2-1$, for $1\leq i<n$.
Its powers induce the ordinary augmentation filtration
on $\kk[\PB_n]$; see
\cite[Lemma~C.1 and Theorem~C.2]{KT_Tangles}.
The resulting completion of $\kk[\PB_n]$ is therefore
$\widehat{\kk[\PB_n]}$, whose grouplike elements form
the Malcev completion $(\PB_n)_{\kk}$.
\end{example}

For the $\kk$-linear tangle category $\bT(\kk)$,
\cite[Theorem~3.6]{KT_Tangles} identifies completion
at the braided augmentation ideal $\mathcal I$ with
the inverse limit of the truncated tangle categories
of \cite[Section~2.2]{KT_Tangles}.
This completion can retain framing parity, reflected
in the $\mathbb Z/2$-residues on components in the
corresponding chord-diagram category;
see \cite[Section~2.3]{KT_Tangles}.

We additionally require the framings to be unipotent.
For a $\kk$-linear ribbon category $\calC$, let
$\mathcal J$ be the monoidal ideal generated by
$\mathcal I$ and the morphisms $\Theta_V-\id_V$,
for all objects $V$. Below, we complete ribbon
categories with respect to the powers of $\mathcal J$.
In particular, we write
$\bT_{\kk}\coloneqq
\varprojlim_{r\geq 0}\bT(\kk)/\mathcal J^{r+1}$.

For the ribbon braid category, the corresponding ideal
on $n$ strands is generated by $\beta_i^2-1$ and
$\tau_j-1$, where $\tau_j$ is the framing twist on
the $j$th strand. Its powers induce the ordinary
augmentation filtration on the group algebra of
the pure ribbon braid group $\PB_n\times\mathbb Z^n$.
Thus its completed group algebra has the Malcev
completion of the pure ribbon braid group as its
group of grouplike elements.

Since $\Theta_V-\id_V\in\mathcal J$, the series
$\log(\Theta_V)$ converges in the completed
endomorphism algebra, and we define
$\Theta_V^a=\exp\bigl(a\log(\Theta_V)\bigr)$
for $a\in\kk$. In this completion,
$\log(\Theta_V^2)=2\log(\Theta_V)$, so
$(\Theta_V^2)^{a/2}=\Theta_V^a$.

\subsubsection{Completion of operads and the tangle category}

Fresse extends Malcev completion from groups to groupoids and
proves that it preserves finite products; see
\cite[Sections~9.1--9.2 and Proposition~9.2.2]{FresseBook1}.
It may therefore be applied aritywise to operads in groupoids.
We write $\calO_{\kk}$ for the resulting Malcev completion
of an operad $\calO$. When $\calO$ is cyclic, its cyclic
structure also extends to $\calO_{\kk}$. In particular,
this gives the cyclic operad $\PaRB_{\kk}^{\mathrm{cyc}}$.

Let $q\bT'_{\kk}$ denote the augmentation-adic completion
of the $\kk$-linearisation of $q\bT'$, using the
Kassel--Turaev augmentation ideal defined above.
It has the same objects as $q\bT'$, and its morphism
spaces are the inverse limits of the corresponding
quotients by powers of this ideal.

The construction of
Proposition~\ref{prop: metric prop action through cyclic endomorphisms}
is compatible with $\kk$-linearisation and $\mathcal I$-adic
completion. It therefore gives an action
\[
\Pi(\PaRB_{\kk}^{\mathrm{cyc}})
\longrightarrow
\End_{q\bT'_{\kk}},
\]
where the endomorphism prop is formed using complete
$\kk$-linear spaces.

\subsection{The proalgebraic Grothendieck--Teichmüller action on tangles}
\label{subsec:GT action on turaev category}

We identify the object-fixing automorphisms of
$\Pi(\PaRB_{\kk}^{\mathrm{cyc}})$ and use them to transport
the metric-prop action on $q\bT'_{\kk}$. We then compare the
resulting ribbon structures with those constructed by Kassel
and Turaev in \cite[Appendix~D]{KT_Tangles}.

Let $\mathsf F(x,y)_{\kk}$ be the Malcev completion over $\kk$
of the free group on generators $x$ and $y$. If
$f\in\mathsf F(x,y)_{\kk}$ and $a,b$ are elements of a
prounipotent group over $\kk$, we write $f(a,b)$ for the image
of $f$ under the homomorphism sending $x$ to $a$ and $y$ to $b$.

\begin{definition}[{\cite{Drin}}]
\label{defn: prounipotent GT}
The \emph{proalgebraic Grothendieck--Teichmüller group}
$\gt(\kk)$ consists of pairs
$(\lambda,f)\in\kk^\times\times\mathsf F(x,y)_{\kk}$ satisfying
\begin{enumerate}
\item
$f(x,y)^{-1}=f(y,x)$;

\item
$x^\nu f(x,y)y^\nu f(y,z)z^\nu f(z,x)=1$, where
$z=(xy)^{-1}$ and $\nu=(\lambda-1)/2$;

\item
$f(x_{12},x_{23}) f(x_{12}x_{13},x_{24}x_{34}) f(x_{23},x_{34}) = f(x_{13}x_{23},x_{34}) f(x_{12},x_{23}x_{24})$
in $(\PB_4)_{\kk}$.
\end{enumerate}
\end{definition}

The group multiplication in $\gt(\kk)$ is given by
\begin{equation}\label{eqn:GT(K) multiplication}
 (\lambda_1,f_1)*(\lambda_2,f_2)
=
\left(
\lambda_1\lambda_2,\,
f_1\bigl(
f_2(x,y)^{-1}x^{\lambda_2}f_2(x,y),
y^{\lambda_2}
\bigr)\cdot f_2(x,y)
\right).   
\end{equation}
The projection $\chi\colon\gt(\kk)\to\kk^\times$, given by $\chi(\lambda,f)=\lambda$, is the cyclotomic character. It is a surjective homomorphism; see \cite[Section~5]{Drin}.

The operadic description of $\gt(\kk)$ through parenthesised
braids originates with Bar--Natan \cite{BN}. Fresse proves the
Malcev-completed identification
\[
\underline{\Aut}_{\Op}(\PaB_{\kk})
\cong
\gt(\kk);
\]
see \cite[Theorem~11.1.7]{FresseBook1}.

The argument of
\cite[Proposition~7.11]{Boavida-Horel-Robertson} applies after
Malcev completion. It first shows that every object-fixing
automorphism of $\PaRB_{\kk}$ preserves the suboperad
$\PaB_{\kk}$. Such an automorphism is determined by its restriction
to $\PaB_{\kk}$ together with the image of the framing generator
$\tau$. The ribbon relation then forces an automorphism with
cyclotomic parameter $\lambda$ to send $\tau$ to $\tau^\lambda$.
Consequently, restriction induces an isomorphism between the
object-fixing automorphism groups of $\PaRB_{\kk}$ and
$\PaB_{\kk}$.

The proof of
Lemma~\ref{lemma: operad maps lift to cyclic operad maps of PaRB}
also applies after Malcev completion. Hence every object-fixing
automorphism of $\PaRB_{\kk}$ preserves its cyclic structure.

\begin{thm}
\label{thm:prounipotent_GT_identifications}
There are natural isomorphisms of groups
\[
\underline{\Aut}_{\Op}(\PaB_{\kk})
\cong
\underline{\Aut}_{\Op}(\PaRB_{\kk})
\cong
\underline{\Aut}_{\Cyc}(\PaRB_{\kk}^{\mathrm{cyc}})
\cong
\gt(\kk),
\]
where all automorphisms fix objects.
\end{thm}

\begin{remark}
Willwacher obtains closely related results for rational
dg Hopf cooperads. He proves that the $\mathsf{GRT}$-action
on the parenthesised ribbon-chord-diagram operad preserves
its cyclic structure
\cite[Proposition~5.2]{willwacher2024cyclic}.
The induced action on the Batalin--Vilkovisky Hopf
cooperad factors through its cyclic homotopy
automorphisms, and the resulting homotopy automorphism
space is computed in
\cite[Corollaries~5.3 and~1.5]{willwacher2024cyclic}.
\end{remark}

For $(\lambda,f)\in\gt(\kk)$, let $F_{\lambda,f}$ denote the
corresponding cyclic-operad automorphism of
$\PaRB_{\kk}^{\mathrm{cyc}}$. On the standard generators,
\[
F_{\lambda,f}(\mu)=\mu,\qquad
F_{\lambda,f}(\alpha)=\alpha\cdot f(x_{12},x_{23}),\qquad
F_{\lambda,f}(\beta)
=\beta\cdot(\beta_{21}\cdot\beta)^{(\lambda-1)/2},
\qquad\text{and}\qquad
F_{\lambda,f}(\tau)=\tau^\lambda.
\]

The restriction-extension argument of
Proposition~\ref{prop:GT-metric-prop-profinite} applies without
change after Malcev completion.

\begin{prop}
\label{prop: GT action on metric prop}
Restriction induces a natural isomorphism
\[
\underline{\Aut}_{\Prop}
\bigl(\Pi(\PaRB_{\kk}^{\mathrm{cyc}})\bigr)
\cong
\gt(\kk),
\]
where the automorphisms fix objects.
\end{prop}

\subsubsection{Comparison with the Kassel--Turaev construction}

In \cite[Appendix~D]{KT_Tangles}, Kassel and Turaev associate to
each $(\lambda,f)\in\gt(\kk)$ a new ribbon monoidal structure on the
same underlying complete $\kk$-linear category. In our notation, its
associator, braiding, and twist are
\begin{equation}
\label{eq: Kassel Turaev transported ribbon operations}
\alpha^{\lambda,f} = \alpha\cdot f(x_{12},x_{23}), \qquad
\beta^{\lambda,f} = \beta\cdot(\beta_{21}\cdot\beta)^{(\lambda-1)/2},
\qquad \text{and} \qquad
\theta^{\lambda,f} 
= \theta^{\lambda}.
\end{equation}
These are \cite[Equations~\textnormal{(D.2)}, \textnormal{(D.3)},
and~\textnormal{(D.5)}]{KT_Tangles}.

Their duality formulas are
\begin{equation}
\label{eq: Kassel Turaev transported duality}
\coev_X^{\lambda,f}
=
\coev_X,
\qquad
\ev_X^{\lambda,f}
=
\ev_X\cdot(\nu'_X\otimes\id_X)
=
\ev_X\cdot(\id_{X^*}\otimes\nu_X),
\end{equation}
where
\begin{equation}
\label{eq: Kassel Turaev duality corrections}
\begin{aligned}
\nu'_X
&=
\left(
(\ev_X\otimes\id_{X^*})
\cdot f(x_{12},x_{23})^{-1}
\cdot\alpha_{X^*,X,X^*}^{-1}
\cdot(\id_{X^*}\otimes\coev_X)
\right)^{-1},\\
\nu_X
&=
\left(
(\id_X\otimes\ev_X)
\cdot\alpha_{X,X^*,X}
\cdot f(x_{12},x_{23})
\cdot(\coev_X\otimes\id_X)
\right)^{-1}.
\end{aligned}
\end{equation}
We have inserted the associators in
\cite[Equations~\textnormal{(D.6)} and~\textnormal{(D.7)}]
{KT_Tangles} to specify the parenthesisations.
Here $\nu'_X$ is an automorphism of $X^*$ and $\nu_X$ is an
automorphism of $X$. Each composite inside parentheses is
congruent to the identity modulo $\mathcal I$, so its inverse
exists in the completed category. The equality of the two
expressions for $\ev_X^{\lambda,f}$ is
\cite[Equation~\textnormal{(D.4)}]{KT_Tangles}.

\begin{prop}
\label{prop: comparison with KT action}
For every $(\lambda,f)\in\gt(\kk)$, the ribbon structure on
$q\bT'_{\kk}$ determined by precomposing the metric-prop action
with $\Pi(F_{\lambda,f})$ agrees with the ribbon structure
constructed by Kassel and Turaev in
\cite[Appendix~D]{KT_Tangles}.
\end{prop}

\begin{proof}
On the operations coming from $\PaRB_{\kk}$, precomposition
with $\Pi(F_{\lambda,f})$ gives the associator, braiding, and
twist in
\eqref{eq: Kassel Turaev transported ribbon operations}. The functor $(X,Y)\mapsto\Hom(X,Y^*)$ assigned to both
$\ev$ and $\coev$ is unchanged, since $\Pi(F_{\lambda,f})$
fixes these operations. To represent this pairing by duality
morphisms for $\alpha^{\lambda,f}$, retain $\coev_X$ and
take $\ev_X^{\lambda,f}$ as in
\eqref{eq: Kassel Turaev transported duality}.

We verify the snake identities using the two expressions for
$\ev_X^{\lambda,f}$. Moving $\nu_X$ and $\nu'_X$ to the input
of the respective composite by naturality of
$\alpha^{\lambda,f}$ gives
\[
\begin{aligned}
(\id_X\otimes\ev_X^{\lambda,f})
\cdot\alpha^{\lambda,f}_{X,X^*,X}
\cdot(\coev_X\otimes\id_X)
&=\nu_X^{-1}\cdot\nu_X\cong\id_X,\\
(\ev_X^{\lambda,f}\otimes\id_{X^*})
\cdot(\alpha^{\lambda,f}_{X^*,X,X^*})^{-1}
\cdot(\id_{X^*}\otimes\coev_X)
&=(\nu'_X)^{-1}\cdot\nu'_X\cong\id_{X^*}.
\end{aligned}
\]
The first calculation uses
$\ev_X^{\lambda,f}=\ev_X\cdot(\id_{X^*}\otimes\nu_X)$,
and the second uses
$\ev_X^{\lambda,f}=\ev_X\cdot(\nu'_X\otimes\id_X)$.

These identities give the duality bijections representing the
unchanged pairing. Explicitly, a morphism $h\colon X\to Y^*$
corresponds to
$\ev_Y^{\lambda,f}\cdot(h\otimes\id_Y)\colon
X\otimes Y\to\emptyset$.
The inverse uses $\coev_Y$ and the transported associator.

Since the evaluation satisfying the snake identities is uniquely
determined by the fixed coevaluation, this identifies the duality
with \eqref{eq: Kassel Turaev transported duality}.
Together with
\eqref{eq: Kassel Turaev transported ribbon operations}, this
gives the required comparison.
\end{proof}

Write $(q\bT'_{\kk})^{\lambda,f}$ for the category equipped
with the transported ribbon structure.

\begin{cor}
The group $\gt(\kk)$ acts by object-fixing automorphisms of
the underlying category $q\bT'_{\kk}$. Each automorphism is
a ribbon monoidal equivalence from the original ribbon
structure to the corresponding transported ribbon structure.
\end{cor}

\begin{proof}
The tangle universal property gives a ribbon monoidal functor
sending the generating strand to itself and the generating
morphisms to their transported counterparts. In particular,
the cup is sent to $\coev_*$ and the cap to
$\ev_*^{\lambda,f}$. The transported symmetric pairing supplies
the required self-duality. Since the assignment preserves the
augmentation filtration, it extends to a $\kk$-linear functor
\[
\Phi_{\lambda,f}\colon
q\bT'_{\kk}\longrightarrow(q\bT'_{\kk})^{\lambda,f}.
\]
The passage to the completion follows
\cite[Theorem~3.7]{KT_Tangles}.

Transport is natural in ribbon monoidal functors and respects
multiplication in $\gt(\kk)$. The corrected evaluations also
respect multiplication, since the coevaluation remains fixed
and the evaluation satisfying the snake identities is unique.
The uniqueness in the tangle universal property therefore
gives the multiplication law for the underlying functors
$\Phi_{\lambda,f}$  \[\Phi_{(\lambda_1,f_1)*(\lambda_2,f_2)}=\Phi_{\lambda_2,f_2}\cdot\Phi_{\lambda_1,f_1}.\] Therefore the assignment $(\lambda,f)\mapsto\Phi_{\lambda,f}$ is a (right) action of $\gt(\kk)$. The functor $\Phi_{\lambda,f}^{-1}$ associated to $(\lambda,f)^{-1}$ is the inverse.
\end{proof}

For every prime $\ell$, restriction along Drinfeld's
homomorphism
$\operatorname{Gal}(\overline{\mathbb Q}/\mathbb Q)
\to\gt(\mathbb Q_\ell)$
gives the corresponding Galois action on
$q\bT'_{\mathbb Q_\ell}$.

\begin{remark}
The Lie-algebraic counterpart to the topological story presented
here is developed by the second author in
\cite{singh2026cyclic}. There, the graded version of
Theorem~\ref{thm:prounipotent_GT_identifications} is established
by identifying the graded proalgebraic
Grothendieck--Teichm\"uller group $\mathsf{GRT}(\kk)$ with the
group of object-fixing automorphisms of the prounipotent cyclic
operad of parenthesised ribbon chord diagrams.
\end{remark}
\subsection{Cyclic formality of the framed little disks operad}
\label{subsec: formality}

We conclude with an alternative proof of the rational formality of the cyclic operad $\mathsf{FD}_2^{\mathrm{cyc}}$. The proof uses the cyclic structure on $\PaRB$ described in Section~\ref{subsec: cyclic PaRB}, together with the action of the proalgebraic Grothendieck--Teichmüller group $\gt(\mathbb Q)$. Recall that a cyclic dg operad is \emph{formal} if it is connected to its homology by a zigzag of quasi-isomorphisms of cyclic dg operads.

The formality of the underlying framed little disks operad was proved by \v{S}evera~\cite{severa2009formality}, whose proof uses $\PaRB$. Giansiracusa and Salvatore subsequently gave an argument for cyclic formality over $\mathbb R$ \cite[Theorem~A]{Salvatore_2012_cyclic_formality}. Campos, Idrissi, and Willwacher later established the rational dg Hopf version needed here \cite[Section~3.1 and Theorem~3.1]{campos2019configuration}. Willwacher subsequently used the cyclic $\mathsf{GRT}(\kk)$-action on the Batalin-Vilkovisky dg Hopf cooperad to show the formality of $\mathsf{FD}_2^{\mathrm{cyc}}$ \cite[Section 5.3]{willwacher2024cyclic}.

Our proof is the cyclic analogue of Petersen's argument deriving formality from the Grothendieck--Teichmüller action \cite{petersen_formality}. We use the cyclic form of the grading-automorphism criterion of Guillén Santos, Navarro, Pascual, and Roig; see \cite[Corollary~5.2.2]{GNPR}.

We recall the following from \cite[Section~5]{GNPR}.

\begin{definition}\label{def:locally_finite}
An automorphism $\sigma$ of a $\kk$-module $M$ is \emph{locally finite} if every element of
$M$ lies in a finite-dimensional $\sigma$-invariant subspace.  For $\lambda\in\kk$ we write
\[W_n(M,\sigma) :=\bigcup_{s\ge1}\ker(\sigma-\lambda^{n})^{s},\] for the primary component of
$\sigma$ at $\lambda^{n}$; if $\sigma$ is locally finite then
$M=C\oplus \big(\bigoplus_{n}W_n(M,\sigma) \big)$, where $C$ is the sum of the primary components for the
remaining irreducible factors of the minimal polynomials.
\end{definition}

\begin{prop}\label{prop: cyclic grading formality criterion}
Let $\kk$ be a field of characteristic zero and $\lambda\in\kk^{\times}$ not be a root of
unity. Let $\calO$ be a cyclic dg operad and suppose that the grading automorphism of
$H_{\ast}(\calO)$, acting on $H_{j}(\calO)$ by $\lambda^{j}$, lifts to a locally finite
automorphism $\sigma$ of $\calO$ as a cyclic dg operad. Then $\calO$ is formal as a
cyclic dg operad.
\end{prop}

\begin{proof}
Let $W_n\calO=W_n(M,\sigma)$ be the primary components of $\sigma$, then we have by Definition~\ref{def:locally_finite} and the local finiteness $\calO=C\oplus\bigoplus_{n}W_n\calO$. This decomposition is
preserved by $d$ and by the $\Sigma^{+}_{l}$-actions, and
$W_m\calO\circ_{i}W_n\calO\subseteq W_{m+n}\calO$ since $\sigma$ acts on
$\calO(l)\otimes\calO(m)$ as $\sigma\otimes\sigma$, the cyclic unit lies in
$W_0\calO$. Thus $\bigoplus_{n}W_n\calO$ is a cyclic dg suboperad. Because $H_{\ast}(\sigma)$ is the grading automorphism and $\lambda$ is not a root of unity, we have $H_{\ast}(C)=0$ and $H_{\ast}(W_n\calO)$ is concentrated in degree $n$, where it is $H_{n}(\calO)$, therefore the inclusion $\bigoplus_{n}W_n\calO\hookrightarrow\calO$ is a quasi-isomorphism. 

Now we set $T\calO=\bigoplus_{n}\tau_{\geq n}W_n\calO$, with
$\tau_{\geq n}$ the canonical truncation: given in degree $n$ by $\tau_{\geq n}W_n\calO := Z_n W_n\calO \oplus \bigoplus_{i>n}W_i\calO $. It is a cyclic dg suboperad because for
$x\in(\tau_{\geq m}W_m\calO)_{p}$ and $y\in(\tau_{\geq n}W_n\calO)_{q}$ the
composite $x\circ_{i}y$ has degree $p+q\geq m+n$. Note that $p+q= m+n$ only if $x$ and $y$ are cycles, in which case $x\circ_{i}y$ is also a cycle. By the above the inclusion $T\calO\hookrightarrow\bigoplus_{n}W_n\calO$ and the projection
$T\calO\to H_{\ast}(\calO)$, given in weight $n$ by $Z_{n}\to H_{n}(\calO)$ and by zero
in degrees $>n$, are quasi-isomorphisms of cyclic dg operads.
\end{proof}

\begin{remark}\label{rem:relation_to_GNPR_Petersen}
Local finiteness holds in either of two conditions: first when $\calO$ is of finite type since $\sigma$ preserves each degree and arity, and second when $\sigma$ lies in an action of $\bG_{m}$ by automorphisms of cyclic dg operads, in which case $C=0$ and the primary decomposition is the weight decomposition. 

Proposition~\ref{prop: cyclic grading formality criterion} is the cyclic analogue of \cite[Corollary~5.2.2]{GNPR}, which is the operadic form \cite[p.~819]{petersen_formality} of \cite[Theorem~12.7]{sullivan}. The transfer to cyclic operads is asserted without proof in \cite[Section 7.1]{GNPR}, therefore we provide a short proof following \cite[Section 5]{GNPR} for completeness. Alternatively, we also note that when $H_{\ast}(\calO(1))\cong\kk$, it follows from the cyclic minimal models of
\cite[Section 7.1]{GNPR} by the argument of \cite{petersen_formality}, the eigenvalue estimate needed there is inherited from the underlying dg operad.
\end{remark}

The following theorem serves as the cyclic, prounipotent analogue to the formality result of
\cite[Theorem~9.1]{Boavida-Horel-Robertson}.

\begin{thm}
\label{thm: cyclic formality framed little disks}
The cyclic operad $\mathsf{FD}_2^{\mathrm{cyc}}$ is rationally
formal.
\end{thm}

\begin{proof}
The aritywise classifying space
$\cs(\PaRB^{\mathrm{cyc}})$ is weakly equivalent, as a cyclic
topological operad, to $\mathsf{FD}_2^{\mathrm{cyc}}$. The
rational-completion argument used in
\cite[Theorem~9.1]{Boavida-Horel-Robertson} applies aritywise to
$\PaRB^{\mathrm{cyc}}$. Indeed, its arity groupoids $\PaRB^{\mathrm{cyc}}$ have finitely
many objects and automorphism groups
$\PRB_n\cong\PB_n\times\mathbb Z^n$ are $\mathbb{Q}$-good, so the Malcev
completion induces an isomorphism on rational homology. Consequently,
the cyclic dg operad
$C_\ast(\cs(\PaRB_{\mathbb Q}^{\mathrm{cyc}});\mathbb Q)$ is a
rational chain model for $\mathsf{FD}_2^{\mathrm{cyc}}$.

The group $\gt(\mathbb Q)$ acts on
$\PaRB_{\mathbb Q}^{\mathrm{cyc}}$ by cyclic-operad automorphisms and
hence acts on this chain model. Its homology is the cyclic
Batalin--Vilkovisky operad
$H_\ast(\mathsf{FD}_2^{\mathrm{cyc}};\mathbb Q)\cong\BV$; see
\cite{getzler_framed_discs}. In homological grading, $\BV$ is
generated by the binary commutative product $m$ in degree zero and
the unary BV operator $\Delta$ in degree one.

Let $(\lambda,f)\in\gt(\mathbb Q)$. The corresponding automorphism $\sigma$ fixes the multiplication operation $\mu$ and therefore fixes the
class $m$. In arity one, $\PaRB(1)$ is a one-object groupoid whose automorphism group $\Aut_{\PaRB(1)}(1,1)\cong \PRB_1 \cong \mathbb Z$ is generated by the framing twist $\tau$, so $\PaRB_{\mathbb Q}(1)$ has automorphism group $\mathbb Q$ and $H_1(\mathrm{B}\PaRB_{\mathbb Q}(1);\mathbb Q)\cong\mathbb Q$, since the automorphism $\sigma$ sends $\tau$ to
$\tau^\lambda$, so it acts by multiplication by $\lambda$ on the
class corresponding to $\Delta$. Since $\BV$ is generated by $m$
and $\Delta$, it acts by multiplication by $\lambda^j$ on the
homological degree-$j$ part of $\BV$.

Note that $\sigma$ lies in $\mathbb G_m$-action due to the map $\gt(\mathbb Q) \rightarrow \mathbb G_m$ sending $(\lambda,f)$ to $\lambda$ where $\mathbb G_m$ acts on homology via the grading action as described above. Therefore $\sigma$ is locally finite by Remark~\ref{def:locally_finite}.

Finally, the cyclotomic character
$\gt(\mathbb Q)\to\mathbb Q^\times$ is surjective; see
\cite[Section~5]{Drin} and \cite{petersen_formality}. We may
therefore choose an element whose cyclotomic parameter $\lambda$ is not a root
of unity. Its action on the chain model lifts the corresponding
grading automorphism of $\BV$. Proposition~\ref{prop: cyclic grading
formality criterion} now proves that
$\mathsf{FD}_2^{\mathrm{cyc}}$ is rationally formal.
\end{proof}

\bibliographystyle{alphaurl}
\bibliography{Version_2/bib1.bib}

\end{document}